\documentclass[11pt]{article}

\usepackage{graphicx,caption}
\usepackage{subfig}
\usepackage{amssymb,amsmath}
\usepackage{color}
\usepackage[colorlinks,filecolor=blue,citecolor=cyan]{hyperref}
\usepackage{setspace}
\usepackage{eufrak}

\newcommand{\sss}{\scriptscriptstyle}	
\newcommand{\refEqn}[1]{(\ref{#1})}
\newcommand{\refEqns}[2]{\mbox{(\ref{#1})--(\ref{#2})}}
\newcommand{\refFigsSecs}[2]{\ref{#1}--\ref{#2}}
\newcommand{\tablesize}{\footnotesize}

\newcommand{\rS}{\mathrm S}

\newcommand{\torus}{{\mathbb T}}

\newcommand{\reals}{{\mathbb R}}
\newcommand{\integers}{{\mathbb Z}}
\newcommand{\naturals}{{\mathbb N}}
\newcommand{\diag}{\hbox{\rm diag}}
\newcommand{\mmbf}[1]{\mbox{\boldmath ${#1}$}}
\newcommand{\F}{\mathbf{F}}
\newcommand{\FL}{\F_{\mathrm{{\sss{L}}}}}
\newcommand{\FNL}{\F_{\mathrm{{\sss{NL}}}}}
\newcommand{\Feqzero}{\mbox{$\F = \mmbf{0}$}}
\newcommand{\nleargs}{\mmbf{\hat{a}}, \varphi, S, T}
\newcommand{\Fargs}{\mbox{$\F(\nleargs)$}}
\newcommand{\FLargs}{\mbox{$\FL(\nleargs)$}}
\newcommand{\FNLargs}{\mbox{$\FNL(\mmbf{\hat{a}})$}}
\newcommand{\Jac}{\mathbf{DF}}
\newcommand{\JL}{\Jac_{\mathrm{{\sss{L}}}}}
\newcommand{\JNL}{\Jac_{\mathrm{{\sss{NL}}}}}

\newcommand{\Ja}{J_{\hat{\mmbf{{\sss{a}}}}}}
\newcommand{\productJa}[1]{\Ja \hat{\mmbf{#1}}}
\newcommand{\ii}{\mathrm{i}}
\newcommand{\e}{\mathrm{e}}
\newcommand{\dd}{\mathrm{d}}

\newcommand{\jacobian}{Jacobian}
\newcommand{\cglegroup}{G}
\newcommand{\cgleisotropy}[1]{\cglegroup_{\!{\sss{{#1}}}}}
\newcommand{\Ap}[1]{A_{\!{\sss #1}}}
\newcommand{\Arot}{\theta}                       
\newcommand{\spacetrans}{\sigma}          
\newcommand{\timetrans}{\tau}
\newcommand{\Lx}{L_{\sss x}}                          
\newcommand{\km}{k_{\sss m}}                        
\newcommand{\omegan}{\omega_{\sss n}}     
\newcommand{\ra}{$\, \longrightarrow \ \ $}
\newcommand{\lsym}{l}
\newcommand{\ce}{c_{\sss{1}}}
\newcommand{\co}{c_{\sss{2}}}

\newcommand{\am}[1]{a_{\sss{{#1}}}}
\newcommand{\amn}[2]{\hat{a}_{\sss{{#1},{#2}}}}
\newcommand{\M}{J}
\newcommand{\U}{\M_s^{-1} \M_r}
\newcommand{\m}{p}
\newcommand{\n}{q}
\newcommand{\x}{z}
\newcommand{\Ai}[1]{\mathcal{A}^{{\sss ({#1})}}}
\newcommand{\pt}[2]{{#1}_{{\sss{{#2}}}}}
\newcommand{\pts}[3]{{#1}_{{\sss{{#2}}}}^{{\sss{{#3}}}}}
\newcommand{\Ni}{\mathrm{N}^{\sss (i)}}

\newcommand{\baseSpace}{\mathcal{B}}
\newcommand{\moduliSpace}{\EuFrak{M}}

\newcommand{\moduliSpaceI}{\mathcal{I}}
\newcommand{\solSpace}{\mathcal{S}}
\newcommand{\fiberM}[2]{\mathcal{#1}_{\sss{#2}}}
\newcommand{\dSpace}{\Sigma}
\newcommand{\diSpace}[1]{\dSpace_{\sss{#1}}}

\newcommand{\contStepVar}{\mathrm{j}}
\newcommand{\contStep}[1]{\contStepVar_{\sss{\mathrm{{#1}}}}}

\newcommand{\contParamVar}{\Lambda}  
\newcommand{\contParam}[1]{\contParamVar_{\sss{\mathrm{{#1}}}}}
\newcommand{\stepChange}[1]{\delta_{\sss{\mathrm{{#1}}}}}
\newcommand{\uAlgo}[1]{\mmbf{u}_{\sss{\mathrm{{#1}}}}}
\newcommand{\fAlgo}{\mathcal{F}}
\newenvironment{myenumerate}
{ \begin{enumerate}
    \setlength{\itemsep}{0pt}
    \setlength{\parskip}{0pt}
    \setlength{\parsep}{0pt}     }
{ \end{enumerate}                  }

\begin{document}
                     
\pagestyle{plain}
\abovedisplayskip=5pt
\belowdisplayskip=5pt
\abovedisplayshortskip=5pt
\belowdisplayshortskip=5pt
\belowcaptionskip=-5pt

\title{Numerical Continuation of Invariant Solutions of the \\Complex Ginzburg-Landau Equation}
\author{Vanessa L\'{o}pez\footnotemark[2]}

\renewcommand{\thefootnote}{\fnsymbol{footnote}}
\footnotetext[2]{{IBM} Research, T.~J. Watson Research Center, 1101 Kitchawan Road, Route 134,
    Yorktown Heights, NY, 10598 USA
    (lopezva@us.ibm.com).}
\renewcommand{\thefootnote}{\arabic{footnote}}

\maketitle

\begin{abstract}

We consider the problem of computation and deformation of group orbits of solutions of the 
complex Ginzburg-Landau equation ({CGLE}) with cubic nonlinearity in $1\!+\!1$ space-time dimension
invariant under the action of the three-dimensional Lie group of symmetries 
$A(x,t) \rightarrow \e^{\ii\Arot}A(x+\spacetrans,t+\timetrans)$.
From an initial set of group orbits of invariant solutions,
for a particular point in the parameter space of the {CGLE}, we obtain new sets of group orbits of invariant solutions
via numerical continuation along paths in the moduli space.
The computed solutions along the continuation paths are unstable, and have multiple modes and frequencies
active in their spatial and temporal spectra, respectively.
Structural changes in the moduli space resulting in symmetry gaining / breaking  associated often with the spatial
reflection symmetry
$A(x,t) \rightarrow A(-x,t)$ of the {CGLE} were frequently uncovered
in the parameter regions traversed.
\\ \textbf{Key Words:} invariant solutions, complex Ginzburg-Landau equation, 
                                        continuous symmetries, numerical continuation
\end{abstract}

\section{Introduction}
\label{sec:introduction}

We consider the problem of numerical computation and deformation of solutions of evolutionary partial differential 
equations ({PDEs}) fixed by the action of a subgroup of a Lie group $\Gamma = \Gamma_1 \times \reals$ of continuous 
symmetries of the {PDEs}, where $\reals$ is the group of time translations and $\Gamma_1$ is non-trivial.  
Within this context, such \emph{invariant solutions} are also known as \emph{relative periodic orbits} or 
\emph{relative time-periodic solutions} of an (autonomous) \emph{equivariant dynamical system}.  
In this paper, we work with the complex Ginzburg-Landau equation with cubic nonlinearity in $1\!+\!1$ space-time dimension, 
with $\Gamma_1 = \torus^2 \ (= \rS^1 \times \rS^1)$ -- the two-torus.  
We note, however, that it should be straightforward 
to apply the methodology described in this paper to other evolutionary parameter-dependent {PDEs} invariant under 
the action of a group of continuous transformations.   
  
The complex Ginzburg-Landau equation ({CGLE}) is a widely studied PDE which has become a model problem
for the study of nonlinear evolution equations exhibiting chaotic spatio-temporal dynamics, as well as being of 
interest in the context of pattern formation.  It has applications in various fields, including fluid dynamics and 
superconductivity.  (For details see, for example, \cite{aranson02,levermore96,mielke02,vanSaarloos94} and 
references therein.) Following \cite{lopez05}, we consider here the {CGLE} with cubic nonlinearity in 
one spatial dimension,
\begin{equation}   \label{eqn:cgle_pde}
    \frac {\partial A}{\partial t}  = 
    R A + (1 + \ii \nu) \frac{\partial^2 A}{\partial x^2} - (1 + \ii \mu) A|A|^2,
\end{equation}
with periodic boundary conditions
\begin{equation}   \label{eqn:cgle_bcs}
    A(x,t)  =  A(x+\Lx,t),    
\end{equation}
and spatial period $L_x = 2\pi$.
The {CGLE} also appears in the literature in the form
\begin{equation}    \label{eqn:cgle_pde2}
    \frac {\partial A}{\partial t}   = 
    A + (1 + i \nu) \frac{\partial^2 A}{\partial x^2} - (1 + i \mu) A|A|^2,
    \qquad A(x,t) = A(x+L,t),
\end{equation}
but note that with a change of variables
$  \ x \rightarrow (L_x/L) \, x, \  
    t \rightarrow (L_x/L)^2 \, t, \
    A \rightarrow (L_x/L) \, A  \ $
one obtains equation~\refEqn{eqn:cgle_pde}, with $R = (L / L_x)^2$.
Thus we adopt the formulation~\refEqn{eqn:cgle_pde} without loss of generality and, henceforth,
when we refer to \emph{the {CGLE}} we mean equation \refEqn{eqn:cgle_pde} with the boundary 
conditions \refEqn{eqn:cgle_bcs} unless otherwise noted.

Equation \refEqn{eqn:cgle_pde} describes the time evolution of a complex-valued field $A(x,t)$.  The parameters 
$R$, $\nu$, and $\mu$ in the equation are real.  
When $R > 0$ there is, in general, nontrivial spatio-temporal behavior and this is therefore the region of interest.  
The parameters $\nu$ and $\mu$ are measures of the linear and nonlinear dispersion, 
respectively \cite{aranson02,levermore96}.  

As will be discussed in detail in Section~\ref{sec:problem_statement}, the {CGLE} has a three-parameter group 
$\cglegroup = \torus^2 \times \reals$ of continuous symmetries generated by space-time translations and a rotation 
of the complex field $A$.  Thus, we focus our study on invariant solutions of the {CGLE}, namely, the ones that in
addition to \refEqn{eqn:cgle_pde} and \refEqn{eqn:cgle_bcs} satisfy
\begin{equation}  \label{eq:Ainvariant}
    A(x,t)  =  \e^{\ii\varphi}A(x+S, t+T)
\end{equation} 
for some $(\varphi, S) \in \torus^2$ and $T > 0$.   
The interest here is on invariant solutions of the 
{CGLE} having multiple frequencies active in their temporal spectrum, not on 
single-frequency solutions $A(x,t) = B(x) \e^{\ii\omega t}$ \cite{doelman89,holmes86,kapitula96} or 
generalized traveling waves $A(x,t) = \rho(x-vt) \e^{\ii\phi(x-vt)} \e^{\ii\omega t}$, 
where $\rho$ and $\phi$ are
real-valued functions and $\omega$ is some frequency \cite{aranson02,brusch01,mielke02,vanHecke02},
which have been considered more extensively than the multiple-frequency class.
The {CGLE} is also invariant under the action of the discrete group of transformations 
$A(x,t) \rightarrow A(-x,t)$ and thus solutions of the {CGLE} may also be fixed by this $\integers_2$-symmetry.
While it is not uncommon in studies to center on
solutions fixed by the $\integers_2$-symmetry (for example, even solutions), 
we make no such restriction here in order to work with a broader solution space.

\begin{figure}[!t]    
    \centering
    {\includegraphics[trim = 30mm 20mm 20mm 15mm, clip, width=0.44\textwidth]{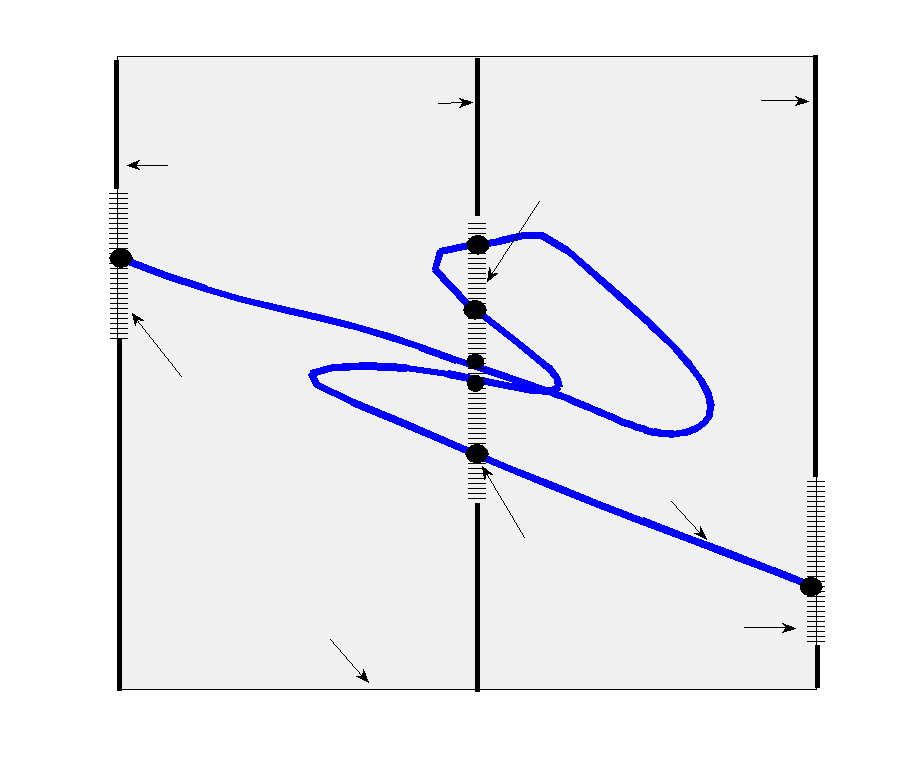}}
    \put(-140,18){$\baseSpace$}
    \put(-175,40){$\moduliSpace$}
    \put(-175,140){$\fiberM{M}{R_{\mathrm{o}},\nu_{\mathrm{o}},\mu_{\mathrm{o}}}$}
    \put(-178,78){$\fiberM{I}{R_{\mathrm{o}},\nu_{\mathrm{o}},\mu_{\mathrm{o}}}$}
    \put(-55,157){$\fiberM{M}{R_{\mathrm{n}},\nu_{\mathrm{n}},\mu_{\mathrm{n}}}$}
    \put(-45,15){$\fiberM{I}{R_{\mathrm{n}},\nu_{\mathrm{n}},\mu_{\mathrm{n}}}$}
    \put(-130,155){$\fiberM{M}{R,\nu,\mu}$}
    \put(-77,133){$\fiberM{I}{R,\nu,\mu}$}
    \put(-82,37){{\small{$\cglegroup$-orbit}}}
    \put(-80,65){{\small{continuation}}}
    \put(-60,57){{\small{path}}}
    \caption{{Fibered space $(\moduliSpace,\baseSpace,\pi)$}}
    \label{fig:intro_plot}
\end{figure}

Since we are actually working with a 3-parameter family \refEqn{eqn:cgle_pde} of equations, 
this family defines {\it implicitly} a fibered space $\solSpace\stackrel{p}{\to}\baseSpace$
over the space of parameters $\baseSpace=\{(R,\nu,\mu)\}\subset\reals^{3}$, where $\solSpace$ is the total space of
solutions of \refEqn{eqn:cgle_pde} and $\baseSpace$ forms the base of the fibered space. Moreover, the group $\cglegroup$
acts on the total space of solutions $\solSpace$. Therefore, we consider the quotient fibered space 
$\moduliSpace\stackrel{\pi}{\to}\baseSpace$ modulo this action. Here $\moduliSpace=\solSpace/\sim$
is the total moduli space,  where $\sim$ is a
relation between the points of $\solSpace$ established by the group action which is compatible with $p$, that is,  for any 
$s',s\in\solSpace$, $s'\sim s$ if and only if $p(s')=p(s)$ and there exists a $g\in\cglegroup$ such that $s'=g \!\cdot\! s$. 
Then $\pi$ is the map induced by $p$ after taking the quotient, and the points of $\moduliSpace$ are in one-to-one
correspondence with $\cglegroup$-orbits whose elements are all mapped by $p$ to the same point in the base $\baseSpace$.
Thus, geometrically we have a fibered space, that is, a triple $(\moduliSpace,\baseSpace,\pi)$, depicted in 
Figure~\ref{fig:intro_plot},  whose fibers $\fiberM{M}{R,\nu,\mu}=\pi^{-1}(R,\nu,\mu)$ over each point 
$(R,\nu,\mu)\in\baseSpace$ of the base are moduli spaces of solutions of the {CGLE}.  Note that we do not know the 
explicit form of the map $\pi$. It is defined implicitly by equation \refEqn{eqn:cgle_pde}.
In essence, it is our goal to understand and reveal its properties.
Therefore, the aim of the present
study is to acquire a more global view of (a part of) the fibered space  of
$\cglegroup$-orbits of the {CGLE} and its structure as we move around the point
$(R,\nu,\mu)$ in the base space $\baseSpace$.

More precisely, and referring again to Figure~\ref{fig:intro_plot}, here we are interested in the (sub)fibered space  
$(\moduliSpaceI,\baseSpace,\pi|\moduliSpaceI)\subset(\moduliSpace,\baseSpace,\pi)$, where the points of the subspace 
$\moduliSpaceI\subset\moduliSpace$ are $\cglegroup$-orbits of {\it invariant}
solutions of the {CGLE}. Namely, these are solutions that satisfy, in addition to \refEqns{eqn:cgle_pde}{eqn:cgle_bcs}, the
functional equation \refEqn{eq:Ainvariant}.  Then the fiber of $\pi|\moduliSpaceI$, 
$\fiberM{I}{R,\nu,\mu}=(\pi|\moduliSpaceI)^{-1}(R,\nu,\mu)\subset\fiberM{M}{R,\nu,\mu}$, over each point  
$(R,\nu,\mu)\in\baseSpace$ of the base  is a moduli space of such invariant solutions of the {CGLE}.
Note that a $\cglegroup$-orbit in  $\fiberM{I}{R,\nu,\mu}$ is determined uniquely by a quadruple 
$\big(A(x,t,R,\nu,\mu),\varphi(R,\nu,\mu),S(R,\nu,\mu),T(R,\nu,\mu)\big)$, where $A(x,t,R,\nu,\mu)$ is an element of the 
orbit (that is, an invariant solution) over the point $(R,\nu,\mu)\in\baseSpace$.  Further, for each point 
$(R,\nu,\mu)\in\baseSpace$ the space $\fiberM{I}{R,\nu,\mu}$ is a union 
$\fiberM{I}{R,\nu,\mu}=\bigcup_{\alpha}\diSpace{R,\nu,\mu}^{(\alpha)}$ 
of symmetry classes $\diSpace{R,\nu,\mu}^{(\alpha)}\subset\fiberM{I}{R,\nu,\mu}$ of $\cglegroup$-orbits. 
A number of such $\cglegroup$-orbits and their symmetry classes were found in \cite{lopez05} at a particular point 
$(R_{\mathrm{o}},\nu_{\mathrm{o}},\mu_{\mathrm{o}})\in\baseSpace$. 
The main goal here is to understand structural changes in the
spaces $\fiberM{I}{R,\nu,\mu}$ as we trace paths in the fibered space  $(\moduliSpaceI,\baseSpace,\pi|\moduliSpaceI)$,
starting from a set of $\cglegroup$-orbits in the fiber $\fiberM{I}{R_{\mathrm{o}},\nu_{\mathrm{o}},\mu_{\mathrm{o}}}$ and
carrying them into another fiber $\fiberM{I}{R_{\mathrm{n}},\nu_{\mathrm{n}},\mu_{\mathrm{n}}}$ over a point 
$(R_{\mathrm{n}},\nu_{\mathrm{n}},\mu_{\mathrm{n}}) \ne (R_{\mathrm{o}},\nu_{\mathrm{o}},\mu_{\mathrm{o}})$ in the base 
$\baseSpace$ using a path following method \cite{bk:ortega}.

Indeed, structural changes associated with additional symmetry breaking or gaining (vanishing and appearance of 
symmetry classes in $\fiberM{I}{R,\nu,\mu}$)
were frequently uncovered in the parameter regions traversed.
This includes the identification of new symmetry classes (see, for example, the
metamorphosis of the moduli space along the path $\Ai{13}$ --
Section~\ref{sec:discussion_particular_solutions} and Figure~\ref{fig:sol13_abs}).
Thus a complex and interesting structure of the fibered space $(\moduliSpaceI,\baseSpace,\pi|\moduliSpaceI)$ was revealed.
Sections~\ref{sec:problem_statement} and \ref{sec:results} describe in detail the additional symmetries that
are being gained or broken at particular values of the {CGLE} parameters.
Such (abrupt) structural changes amount to a kind of ``phase transition'' in the moduli space,
which is an interesting aspect to be considered as a focus for a subsequent detailed investigation.

To put all of the above in context, in contrast to studies which are concerned with continuation (deformation) of a single critical,
singular or other non-generic point (with or without symmetry) of the space of solutions of \refEqn{eqn:cgle_pde} (for example,
solutions which are stable, steady-states, or of the aforementioned traveling waves class) and their bifurcations,
ours should be viewed as a study ``in the large'' of properties of the 3-dimensional family of  moduli spaces  
$\fiberM{I}{R,\nu,\mu}$ of invariant solutions within a certain general functional class, being approximated with a 
spectral-Galerkin discretization (the details of which appear in Section~\ref{sec:discretization_nleqns}).

As we focus on solutions of the {CGLE} having the invariance \refEqn{eq:Ainvariant},
henceforth when we refer to \textit{$\cglegroup$-orbits (of solutions of the {CGLE})} we mean $\cglegroup$-orbits of 
\textit{invariant} solutions of the {CGLE} which satisfy \refEqn{eq:Ainvariant}, unless otherwise indicated.
We also note that during the continuation, most often
the final parameter region of interest was sought by moving in the direction of varying  values of $R$.
However, at times we had to venture into a subdomain of the
{CGLE} parameter space $\baseSpace$ by moving in a direction of varying $\nu$ and $\mu$ as well.
Newton's method, which is commonly used in path following methods \cite{bk:allgower,bk:kuznetsov}, is employed to 
solve an underdetermined system of nonlinear algebraic equations resulting from the discretization of the {CGLE}.
To the best of our knowledge, the way in which the Newton step is computed here is new.
The approach is conceptually simple, yet that is where its value lies:  it led to the efficient computation of an accurate 
Newton step, making the solution of a computationally challenging problem with a large number
(up to 32,260) of unknowns practical without the need of a cluster or supercomputer.  
These and other aspects of the numerical methodology
are discussed in Section~\ref{sec:methodology}.  We note here that
the Newton step used is defined from the Moore-Penrose inverse \cite{bk:benisrael}.
This is one technique used in numerical continuation \cite{bk:allgower,wulff06}, without the
need to define phase, or gauge, conditions \cite{bk:kuznetsov} to augment the underdetermined system.
This offers an advantage since the best (or a suitable) choice of phase conditions
may be problem dependent.  The question of whether to impose phase conditions, or to simply work directly with the
underdetermined system, thus arises.  Here we chose to explore the latter approach.  
However, understanding the advantages that working 
with phase conditions may offer over the chosen approach is important and should be considered as a follow-up investigation.

As a bi-product of our study we note that,
taking the presence of positive Lyapunov exponents for typical (that is, non-invariant) solutions
as an indication of chaotic dynamics \cite{bk:ott}, both the initial and final parameter regions in our study
exhibit chaotic behavior.  Specifically, non-invariant solutions 
in the initial and final parameter region have, respectively, $5$ and $16$ positive 
Lyapunov exponents.\footnote{Lyapunov exponents for typical (non-invariant) solutions
were computed using the technique from \cite{benettin80a}.}
This provides  another motivation for conducting this study, which is to evaluate 
the potential benefits of using numerical continuation (on problems with a large number of unknowns) 
to continue multiple, distinct, unstable 
invariant solutions from one chaotic regime into another with the aim of ending, again, with multiple, distinct,
unstable invariant solutions in the final region.
One question that arises (see also \cite{chandler13}) is whether a significant number
of the distinct solutions used as initial points to continue on will actually lead to solutions in distinct 
$\cglegroup$-orbits in the final parameter region.  
While the possibility of this not happening cannot be ruled out, 
we found that options like alternating the choice
of continuation parameter, or particular settings for tuning parameters in the numerical solvers, can increase the possibility of
reaching a multitude of distinct $\cglegroup$-orbits in the final desired parameter region.

A detailed account of the results obtained is provided in Section~\ref{sec:results}.  We note here that 
the set of $\cglegroup$-orbits in the fiber $\fiberM{I}{R_{\mathrm{o}},\nu_{\mathrm{o}},\mu_{\mathrm{o}}}$ used as starting
points in the path following method correspond to the first 15 $\cglegroup$-orbits listed in the 
Appendix from \cite{lopez05}; these were selected simply to follow the order listed in said Appendix.  
The $\cglegroup$-orbits from this initial set were carried from the point $(R_{\mathrm{o}},\nu_{\mathrm{o}},\mu_{\mathrm{o}}) = (16,-7,5)$ in the {CGLE}
parameter space to the point $(R_{\mathrm{n}},\nu_{\mathrm{n}},\mu_{\mathrm{n}}) = (100,-7,5)$.
The number of unknowns to solve for ranged between 4,000 and 32,260.
Both the number of 15 $\cglegroup$-orbits from \cite{lopez05} and the final point  $(R_{\mathrm{n}},\nu_{\mathrm{n}},\mu_{\mathrm{n}})$ in the {CGLE} parameter 
space were chosen because we deemed them to be sufficient to help us gain insight into the symmetry changes
occurring in the spaces $\fiberM{I}{R,\nu,\mu}$ of the fibered space $(\moduliSpaceI,\baseSpace,\pi|\moduliSpaceI)$,
as well as to allow us to evaluate the potential for success of the proposed 
approach for computing multiple unstable invariant solutions in fixed parameter regions of a dynamical system which exhibits 
chaotic behavior.  

The initial set of $\cglegroup$-orbits led to distinct, new $\cglegroup$-orbits of invariant 
solutions of the {CGLE} along the continuation paths and in the final parameter region.  The solutions
in the resulting $\cglegroup$-orbits are unstable,
and have multiple modes and frequencies active in their spatial and temporal spectra, respectively.
The fact that the computed solutions are unstable suggests that they may belong to the set of
(infinitely many) unstable periodic orbits embedded in chaotic attractors \cite{bk:chaosbook,lan10,chandler13}.
This direction, by itself, is certainly very interesting to pursue in a future study of the dynamics
of the {CGLE}.

To conclude the introduction we note that
previous numerical continuation studies of the {CGLE} include \cite{takac98}, 
where bifurcations from a stable rotating wave to two-tori (of the generalized traveling wave class) were identified. 
Values of $\Lx = 1$ and $R \le 180$ in the formulation \refEqns{eqn:cgle_pde}{eqn:cgle_bcs} were
considered,
giving $L \approx 13.42$ for the maximum length of the spatial period in the formulation \refEqn{eqn:cgle_pde2}.
In comparison, the values $\Lx = 2\pi$ and $R \le 100$ in our study yield a maximum value of 
$L \approx 62.83$ in \refEqn{eqn:cgle_pde2}.  
The values of $\nu$ and $\mu$ used in \cite{takac98} are different from those in the current study,
but in both cases they belong to the Benjamin-Feir unstable region $1 + \mu \nu < 0$ \cite{vanSaarloos94}.
A different study \cite{brusch01} considers traveling waves solutions, where the
{CGLE} reduces to a system of three coupled ordinary differential equations ({ODEs}).  
Continuation was performed on the system of three {ODEs} for different values of 
$L$ up to $512$ and various chaotic regions were classified.

Other studies can be found in \cite{Moon}, where transition to chaos from a limit cycle of the {CGLE} is investigated, 
\cite{keefe85}, in which the bifurcation structure and dynamics of even solutions of the {CGLE} are analyzed,
and \cite{lloyd05}, which studies the dynamics of the {CGLE} in heteroclinic cycles,
focused on invariant $\integers_2$-subspaces.
A numerical study on solutions fixed by the $\integers_2$-symmetry of the {CGLE} and
their stability with respect to symmetry-breaking perturbations appears in \cite{aston99},
where values of $R = 1.05, 16, 36$ are considered, and a spatial period of $\Lx = 2\pi$ was used
(the latter being the same as in the current study).
The subsequent study \cite{aston00} considers symmetry-breaking perturbations for solutions fixed
by the spatial translation symmetry, for (discrete) values of the parameter $R$ in the range $[4.2, 80]$.

\section{Invariant Solutions of the {CGLE} and their Properties}
\label{sec:problem_statement}

The {CGLE} has a number of well known symmetries that are central to its behavior \cite{aranson02}.
In particular, equations \refEqns{eqn:cgle_pde}{eqn:cgle_bcs} have a three-parameter group  
\begin{equation}    \label{eqn:Ggroup}
    \cglegroup = \torus^2 \times \reals
\end{equation}
of continuous symmetries generated by space-time translations $x \rightarrow x + \spacetrans$, $t \rightarrow t + \timetrans$
and a rotation $A \rightarrow \e^{\ii\Arot}A$ of the complex field $A(x, t)$, 
in addition to being invariant under the action of the discrete group of transformations $A(x,t) \rightarrow A(-x,t)$ 
of spatial reflections.  In other words, if $A(x, t)$ is a solution of equations \refEqns{eqn:cgle_pde}{eqn:cgle_bcs}, 
then so are 
\begin{align}   
    \e^{\ii\Arot} A(x, t)&,            \label{cglesym1} \\
    A(x + \spacetrans, t)&,       \label{cglesym2} \\
    A(x, t + \timetrans)&,          \label{cglesym3} \\
    A(-x, t)&,                              \label{cglesym4}
\end{align}
for any $(\Arot, \spacetrans, \timetrans) \in \cglegroup$.  
In the present study it is the group $\cglegroup$ generated by the continuous symmetries \refEqns{cglesym1}{cglesym3} 
which (explicitly) enters the problem formulation. Namely, for a given solution $A(x,t)$ of the {CGLE}, let us consider the 
isotropy subgroup $\cgleisotropy{A}$ of $\cglegroup$ at $A$,
\begin{equation}     \label{eqn:isotropy_subgroup}
    \cgleisotropy{A}  =  \{ (\varphi, S, T) \in \cglegroup \ \, | \ \,  A(x,t) =  \e^{\ii \varphi} A(x+S, t+T) \},
\end{equation}
which consists of elements of the symmetry group $\cglegroup = \torus^2 \times \reals$  leaving $A$ invariant.
With that in mind, we pose the problem:
seek solutions $A(x,t)$ of the {CGLE} satisfying
\begin{equation}   \label{eqn:cgle_invariant_solution}
    A(x,t)  =  \e^{\ii \varphi} A(x+S, t+T),
\end{equation}
for $(\varphi, S, T) \in \cglegroup$ also unknown and to be determined. 
In other words, find orbits $\cglegroup \cdot A$ of $\cglegroup$ generated by solutions $A$ of the {CGLE} which are
invariant under the action of some subgroup $\cgleisotropy{A} \subset \cglegroup$,
that is, $\cgleisotropy{A} \cdot A = A$.  Here, at least one subgroup of $\cgleisotropy{A}$ generated by an element 
$( \varphi, S, T ) \in \cglegroup$ is also to be determined.

As is clear from \refEqn{eqn:cgle_invariant_solution}, the case $\varphi = S = 0$, $T > 0$, would result in a time-periodic
solution.  Within the more general context of the problem of seeking solutions of a dynamical system fixed by the action of 
a subgroup of the system's symmetry group (which also contains time translation), as is the case resulting from $T > 0$ and nonzero $\varphi$ or $S$ 
in \refEqn{eqn:cgle_invariant_solution}, such invariant solutions are also referred to as \textit{relative} time-periodic solutions.  
Since the solutions sought must satisfy the boundary (space-periodicity) condition \refEqn{eqn:cgle_bcs},
it is easy to see that if $S = \Lx / q$ for some integer $q > 1$, then $|A(x,t)| = |A(x,t+qT)|$, 
whereas if both $S = \Lx / q$ and $\varphi = 2\pi / q$ for some integer $q > 1$, then $A(x,t) = A(x,t+qT)$ and, 
therefore, $(0,0,qT) \in \cgleisotropy{A}$ (i.e., $A$ is time-periodic, with time period $qT$).

Notice that if $(\varphi, S, T) \in \cgleisotropy{A}$, the triples $(j\varphi, jS, jT), \ j \in \integers,$ are also elements
of the isotropy subgroup $\cgleisotropy{A}$.  Hence, $(\varphi, S, T)$ generates a subgroup of $\cgleisotropy{A}$.
Thus, the problem that we aim to solve numerically can be described succinctly as follows:
\begin{enumerate}
    \item Given a point $\pt{p}{0} = (\pt{R}{0}, \pt{\nu}{0}, \pt{\mu}{0})$ in the parameter space of the {CGLE}, find a 
              solution $\Ap{\pt{p}{0}}(x,t)$ of the {CGLE} and a generator $(\varphi(\pt{p}{0}), S(\pt{p}{0}), T(\pt{p}{0}))$ of a 
              subgroup of the isotropy subgroup $\cgleisotropy{A_{p_{\sss 0}}}$, such that condition
              \refEqn{eqn:cgle_invariant_solution} holds.  That is, $\Ap{\pt{p}{0}}$ is an invariant solution of the {CGLE} under 
              the action of the subgroup of $\cgleisotropy{A_{p_{\sss 0}}}$ generated by $(\varphi(\pt{p}{0}), S(\pt{p}{0}), T(\pt{p}{0}))$.
    \item Then, starting from $\pt{p}{0} = (\pt{R}{0}, \pt{\nu}{0}, \pt{\mu}{0})$, vary the point $p = (R, \nu, \mu)$ along a 
              subspace in the parameter space of the {CGLE},
              ending at a point $\pt{p}{\mathrm{N}} = (\pt{R}{\mathrm{N}}, \pt{\nu}{\mathrm{N}}, \pt{\mu}{\mathrm{N}})$, 
              to find a sequence of new invariant solutions
              $\Ap{p}(x,t)$  and generators $(\varphi(p), S(p), T(p))$ of subgroups of their corresponding isotropy
              subgroups $\cgleisotropy{A_{p}}$.
\end{enumerate}
In reference \cite{lopez05} we found 77 distinct unstable invariant solutions (that is, 77 $\cglegroup$-orbits generated by
distinct invariant solutions) of the {CGLE} at the point 
$\pt{p}{0} = (\pt{R}{0}, \pt{\nu}{0}, \pt{\mu}{0}) = (16,-7,5)$ of the parameter space of the {CGLE}, thus addressing the
first part of the problem.  Here, we take the first 15 of these solutions, per the listing
from the Appendix in \cite{lopez05}, and address the second part of the problem.  Specifically, using numerical
continuation (as described in Section~\ref{sec:methodology}) we found 15 sequences (or discrete continuation paths)
\begin{equation}    \label{eqn:Asequence}
    \begin{split}
        \Ai{i}  =  \{ \Ap{\pts{p}{k}{(i)}}&(x,t) \, ; \, (\varphi(\pts{p}{k}{(i)}), S(\pts{p}{k}{(i)}), T(\pts{p}{k}{(i)}))  \ \, | \ \,  \\
                                          \pts{p}{k}{(i)} &= (\pts{R}{k}{(i)}, \pts{\nu}{k}{(i)}, \pts{\mu}{k}{(i)}) \in 
                                          [9,100] \times [-7,-2.7] \times [-0.05,5.98],  \ \ 0 \le k \le \Ni \},
    \end{split}
\end{equation}
$i = 1,\ldots,15$, of new invariant solutions \refEqn{eqn:cgle_invariant_solution} of the {CGLE} and corresponding
generators of subgroups of their isotropy subgroups $\cgleisotropy{ \Ap{\pts{p}{k}{(i)}}}$.  In \refEqn{eqn:Asequence},
the number $\Ni$ of invariant solutions in a sequence is at least 100 and,
for each $i = 1,\ldots,15$, the final point $\pts{p}{\Ni}{(i)}$ in the {CGLE} parameter space was fixed at 
$\pts{p}{\Ni}{(i)} = (\pts{R}{\Ni}{(i)}, \pts{\nu}{\Ni}{(i)}, \pts{\mu}{\Ni}{(i)}) = (100,-7,5)$.  
Thus, the sequences \refEqn{eqn:Asequence} can be thought of as a deformation of an initial set of distinct
$\cglegroup$-orbits at $\pt{p}{0} = (16,-7,5)$ into a final set of $\cglegroup$-orbits at $\pt{p}{\mathrm{N}} = (100,-7,5)$,
which in this study are also distinct with the only exception being that the final orbits in the sequences $\Ai{2}$ and $\Ai{4}$
at $\pt{p}{\mathrm{N}}$ happened to coincide (details are provided in Section~\ref{sec:results}).
In other words, if we think of the space of $\cglegroup$-orbits as fibered over the parameter space $\baseSpace$ of the {CGLE} 
(the base of the fibered space), then the sequences (or continuation paths) $\Ai{i}$ can be thought of as (discrete) sections of the
fibered space $(\moduliSpaceI,\baseSpace,\pi|\moduliSpaceI)$.  
Interestingly, several of the sequences $\Ai{i}$ that we have computed contain solutions with additional symmetries
(which we describe in detail later in this section), thus revealing an intricate structure of the fibered space $(\moduliSpaceI,\baseSpace,\pi|\moduliSpaceI)$.

Note that the meaning of the space-periodicity boundary condition \refEqn{eqn:cgle_bcs} is that any solution $A(x,t)$ in the
class of solutions of the {CGLE} that we seek has a subgroup in its isotropy subgroup $\cgleisotropy{A}$ which is generated
by $(0, \Lx, 0)$.  In other words we restrict, a priori, the class of solutions of the {CGLE} that we look for to the one that
contains, at a minimum, solutions with symmetry \refEqn{eqn:cgle_bcs}.  This allows us to represent $A(x,t)$ as a Fourier series
\begin{equation}  \label{eqn:xFseries}
     A(x,t)   =  \sum_{m \in \integers} \am{m}(t) \e^{\ii \km x},
\end{equation}
where $\km = 2\pi m / \Lx$ denotes the $m$-th wavenumber in the expansion. 
From the group-invariance condition \refEqn{eqn:cgle_invariant_solution} it then follows that the complex-valued 
Fourier coefficient functions $\am{m}(t)$ in \refEqn{eqn:xFseries} satisfy
\begin{equation}   \label{eqn:am_invariant}
    \am{m}(t)  =  \e^{\ii\varphi} \e^{\ii \km S} \am{m}(t +T)
\end{equation}
for all $m \in \integers$.  Because of the presence of symmetry \refEqn{eqn:cgle_bcs}, the solutions sought
can be restricted to those with elements $(\varphi,S,T) \in \cglegroup$ having $S \in [0,\Lx)$.

Moreover, since the {CGLE} is invariant under the action of the group $\integers_2$ of spatial reflections
$A(x,t) \rightarrow A(-x,t)$, to any solution $A(x,t)$ of the {CGLE} having $(0,\Lx,0)$ and $(\varphi,S,T)$ as
generators of subgroups of the isotropy subgroup $\cgleisotropy{A}$ (defined in \refEqn{eqn:isotropy_subgroup}) there 
corresponds a solution $\tilde{A}(x,t) := A(-x,t)$ having $(0,\Lx,0)$ and $(\varphi,\Lx\!-\!S,T)$ as generators of 
subgroups of the isotropy subgroup $\cgleisotropy{\tilde{A}}$.  This can be seen from the chain of equalities 
\begin{eqnarray}    \label{eqn:Smap}
        \tilde{A}(x,t) \ := \  A(-x,t)   &=&  \e^{\ii\varphi}A(-x + S, t + T)  \hspace*{3.25em} \mbox{by \refEqn{eqn:cgle_invariant_solution}}   \nonumber \\
                                               &=&  \e^{\ii\varphi}\tilde{A}(x - S, t + T)  \hspace*{4em} \mbox{by definition of $\tilde{A}$}   \nonumber \\
                                               &=&  \e^{\ii\varphi}\tilde{A}(x + (\Lx - S), t + T)  \quad \mbox{by \refEqn{eqn:cgle_bcs}.}
\end{eqnarray}
To express the above in a more symmetric form, let us introduce $\delta = |\Lx/2 - S|$.  Then, if $A(x,t)$ is a solution 
of the {CGLE} having $(0,\Lx,0)$ and $(\varphi, \Lx/2 \pm \delta, T)$ as generators of subgroups of the isotropy 
subgroup $\cgleisotropy{A}$, the solution $\tilde{A}(x,t) := A(-x,t)$ has $(0,\Lx,0)$ and $(\varphi, \Lx/2 \mp \delta, T)$ 
as generators of subgroups of the isotropy subgroup $\cgleisotropy{\tilde{A}}$.
We shall call the invariant solutions $(A; \varphi, \Lx/2 \pm \delta, T)$ and $(\tilde{A}; \varphi, \Lx/2 \mp \delta, T)$,
as well as their corresponding orbits $\cglegroup \cdot A$ and $\cglegroup \cdot \tilde{A}$,
\textit{conjugate} to each other under the (involutive) action of the group $\integers_2$ of spatial reflection
symmetry of the {CGLE}. 

Now, while invariance of solutions of the {CGLE} other than that defined by  \refEqn{eqn:cgle_invariant_solution} 
and \refEqn{eqn:cgle_bcs} is not part of the problem formulation 
\refEqns{eqn:isotropy_subgroup}{eqn:cgle_invariant_solution}, it is clearly not excluded from it.
The {CGLE} may admit solutions having symmetries other than (or in addition to) that defined by
\refEqn{eqn:cgle_invariant_solution} and several of the solutions resulting from our study do have additional symmetries.  
In what follows we discuss some such symmetries and their properties.
We emphasize that our treatment on additional types of symmetries exhibited by solutions of the {CGLE} is not
exhaustive, but rather inclusive of material relevant to the discussion on our results in Section~\ref{sec:results}.

For instance, there may exist solutions of the {CGLE} satisfying
\begin{eqnarray}
     A(x,t)  & = &   \e^{\ii 2\pi / \lsym} A(x + \Lx / \lsym, t),     \hspace{2.25em} \mbox{for some $\lsym \in \naturals$, $\lsym > 1$,}
                       \label{additionalsym1}  \\
     A(x, t)  & = &   \hspace*{1em}A(-x + 2\ce, t)          \hspace{3.75em} \mbox{for some $\ce \in \reals$,}            
                       \label{additionalsym2}  \\
     A(x, t)  & = &   -A(-x + 2\co, t)                                  \hspace{4.00em} \mbox{for some $\co \in \reals$.}  
                       \label{additionalsym3}
\end{eqnarray}
Note that \refEqn{additionalsym1} describes solutions fixed by a composition of the actions \refEqn{cglesym2} and 
\refEqn{cglesym1}, and gives $(2\pi / \lsym, \Lx / \lsym, 0)$ as one generator of a subgroup of $\cgleisotropy{A}$.  
From \refEqn{additionalsym1} it is also clear that the absolute value of such a solution has spatial period of $\Lx / \lsym$. 
Furthermore, by substituting condition \refEqn{additionalsym1} into the Fourier series expansion \refEqn{eqn:xFseries}
it follows that the Fourier coefficient functions $\am{m}(t)$ of a solution with symmetry \refEqn{additionalsym1} satisfy
\begin{equation}   \label{eqn:Fcoefs_lsym}
    \am{m}(t)   =  \left \{  \begin{array}{ll}
                                     \mbox{nonzero}  &  \mbox{if $m \in \{\lsym \tilde{m} - 1, \ \tilde{m} \in \integers \}$}  \\
                                     0                           &   \mbox{otherwise}.
                                     \end{array}
                         \right . 
\end{equation}

Symmetries \refEqn{additionalsym2} and \refEqn{additionalsym3} describe solutions that are, respectively,
even about $x = \ce$ or odd about $x = \co$ for some real numbers $\ce, \co$.
(These solutions are fixed by a composition of the actions \refEqn{cglesym4}, \refEqn{cglesym2}, and \refEqn{cglesym1}.)
From \refEqn{additionalsym2} and the periodic boundary condition \refEqn{eqn:cgle_bcs} it follows that a solution
even about $x=\ce$ is also even about $x = \ce + \Lx/2$; similarly a solution
odd about $x=\co$ is also odd about $x = \co + \Lx/2$.
The Fourier coefficient functions $\am{m}(t)$ in \refEqn{eqn:xFseries} of a solution even about $x = \ce$ satisfy 
\begin{equation}   \label{eqn:Fcoefs_even}
    \am{-m}(t)   =  \am{m}(t) \, \e^{\ii \km  2\ce}   \ , \ \ m = 0, 1, 2, \ldots \, ,
\end{equation}
whereas for a solution odd about $x = \co$ one has
\begin{equation}   \label{eqn:Fcoefs_odd}
    \am{-m}(t)   =  - \am{m}(t) \, \e^{\ii \km 2\co}   \\ , \ \ m = 0, 1, 2, \ldots \, .
\end{equation}
From \refEqn{eqn:Fcoefs_lsym}, \refEqn{eqn:Fcoefs_even}, and \refEqn{eqn:Fcoefs_odd} it follows that 
restricting the search for solutions to those possessing symmetries \refEqn{additionalsym1}, \refEqn{additionalsym2}, or
\refEqn{additionalsym3} would lead to a reduction in the number of unknowns.  However, as already mentioned,  
we did not make a priori such a restriction in order to allow for a more general set of solutions.
Finally, we point out that a solution having both symmetries \refEqn{additionalsym1} and \refEqn{additionalsym2} also satisfies
\begin{equation}    \label{eqn:lsym_even}
    A(-x + 2(\ce + \Lx/(2\lsym)), t)  = \e^{\ii 2\pi / \lsym} A(x, t) \, .
\end{equation}
In particular, note that a solution satisfying \refEqn{additionalsym1} for $\lsym = 2$ and which is even about
$x = \ce$ is also odd about $x = \co = \ce + \Lx / 4$.

For a solution satisfying \refEqn{eqn:cgle_invariant_solution} and \refEqn{additionalsym1} it follows that
$(\varphi - 2\pi / \lsym, S - \Lx / \lsym, T)$ is another generator of a subgroup of $\cgleisotropy{A}$.
In particular, if it happens that for such a solution one has $S = \Lx / \lsym$, then $|A(x,t)| = |A(x,t+T)|$ is satisfied, 
whereas if both $S = \Lx / \lsym$ and $\varphi = 2\pi / \lsym$, then $A(x,t) = A(x,t+T)$ also holds.  
As for invariant solutions \refEqn{eqn:cgle_invariant_solution} that also possess symmetry \refEqn{additionalsym2},
note that, for each $m = 0, 1, 2, \ldots$,
\begin{eqnarray}    \label{eqn:am_even1}
    \am{-m}(t)    & = &  \e^{\ii\varphi} \e^{-\ii \km S} \am{-m}(t + T)
                                             \hspace*{4.0em} \mbox{by \refEqn{eqn:am_invariant}}    \nonumber \\
                         & = &  \e^{\ii\varphi} \e^{-\ii \km S}  \am{m}(t + T) \, \e^{\ii  \km 2\ce}   
                                            \hspace*{1.5em} \mbox{by \refEqn{eqn:Fcoefs_even}.}
\end{eqnarray}
On the other hand, for each $m = 0, 1, 2, \ldots$,
\begin{eqnarray}    \label{eqn:am_even2}
    \am{-m}(t)   & = &  \am{m}(t) \, \e^{\ii  \km  2\ce} 
                                            \hspace*{7.0em}  \mbox{by \refEqn{eqn:Fcoefs_even}}    \nonumber \\
                        & = &  \e^{\ii\varphi} \e^{\ii \km S}  \am{m}(t + T) \, \e^{\ii  \km  2\ce} 
                                            \hspace*{1.5em} \mbox{by \refEqn{eqn:am_invariant}.}
\end{eqnarray} 
From \refEqn{eqn:am_even1} and \refEqn{eqn:am_even2} it follows that we must have $\e^{-\ii \km S}  \ = \  \e^{\ii \km S}$
for all $m \in \integers$, which holds whenever $S$ is an integer multiple of  $\Lx/2$ (recall that $\km = 2\pi m / \Lx$).
The case for invariant solutions \refEqn{eqn:cgle_invariant_solution} with the additional symmetry \refEqn{additionalsym3}
is analogous.
Therefore, solutions satisfying \refEqn{eqn:cgle_invariant_solution} which also posses symmetries 
\refEqn{additionalsym2} or \refEqn{additionalsym3} exist in subspaces of the solution space
$( A; \varphi, S, T )$ for which either $S=0$ or $S = \Lx/2$
(since, by the periodic boundary conditions \refEqn{eqn:cgle_bcs}, $S$ can be restricted to be in the interval $[0, \Lx)$).

\begin{figure}[!t]  
        \centering
        \subfloat[$t=0$]{
                {\includegraphics[trim = 32mm 0mm 20mm 10mm, clip,width=0.175\textwidth]{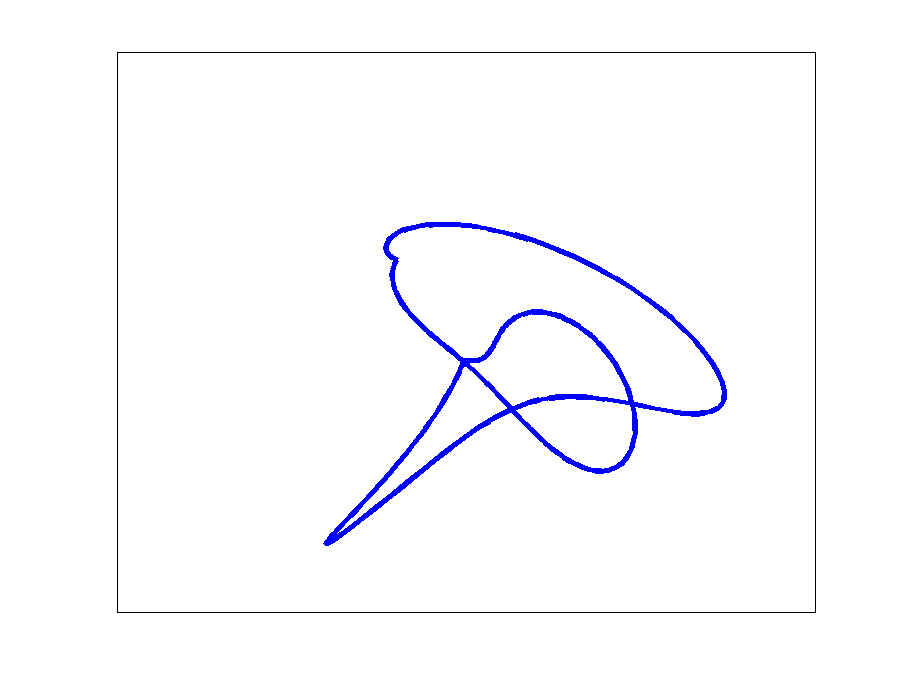}}
                \put(-50,0){\tiny{$\Re(A)$}}
                \put(-85,25){\rotatebox{90}{\tiny{$\Im(A)$}}}
                \label{fig:invariant_sol_t0}}
        \hspace*{1.0ex}
        \subfloat[$t = T/3$]{
                {\includegraphics[trim = 32mm 0mm 20mm 10mm, clip,width=0.175\textwidth]{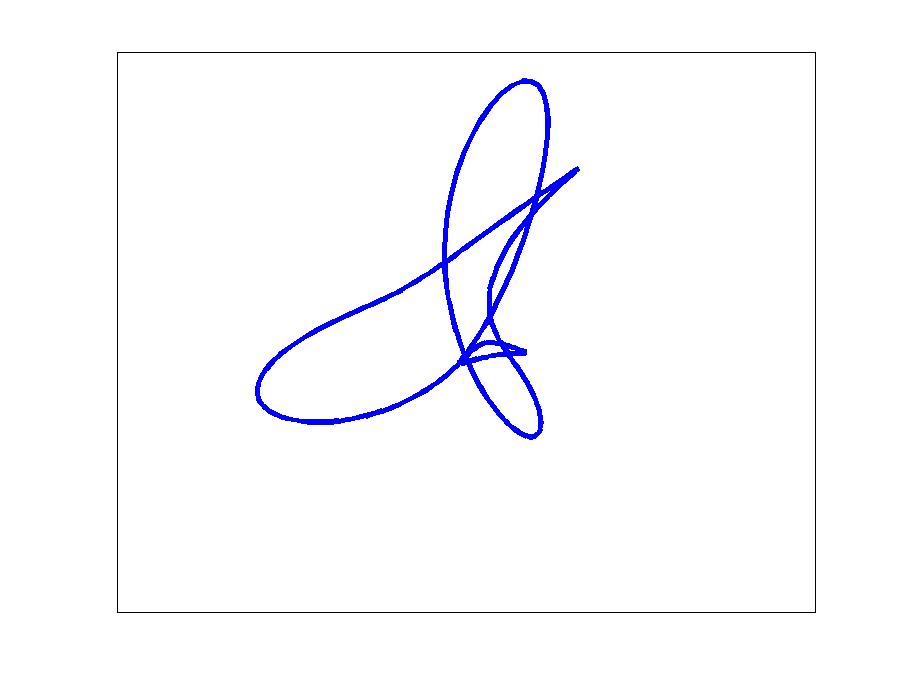}}
                \put(-50,0){\tiny{$\Re(A)$}}
                \put(-83,25){\rotatebox{90}{\tiny{$\Im(A)$}}}
                \label{fig:invariant_sol_t1}}
        \hspace*{1.0ex}
        \subfloat[$t=2T/3$]{
                {\includegraphics[trim = 32mm 0mm 20mm 10mm, clip,width=0.175\textwidth]{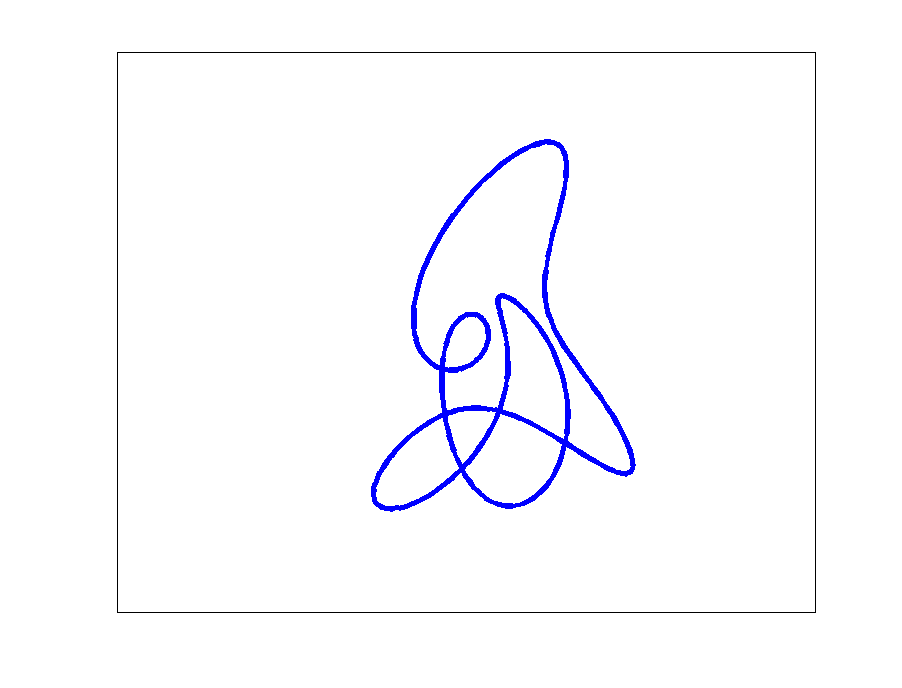}}
                \put(-50,0){\tiny{$\Re(A)$}}
                \put(-83,25){\rotatebox{90}{\tiny{$\Im(A)$}}}
                \label{fig:invariant_sol_t2}}
        \hspace*{1.0ex}
        \subfloat[$t=T$]{
                {\includegraphics[trim = 32mm 0mm 20mm 10mm, clip,width=0.175\textwidth]{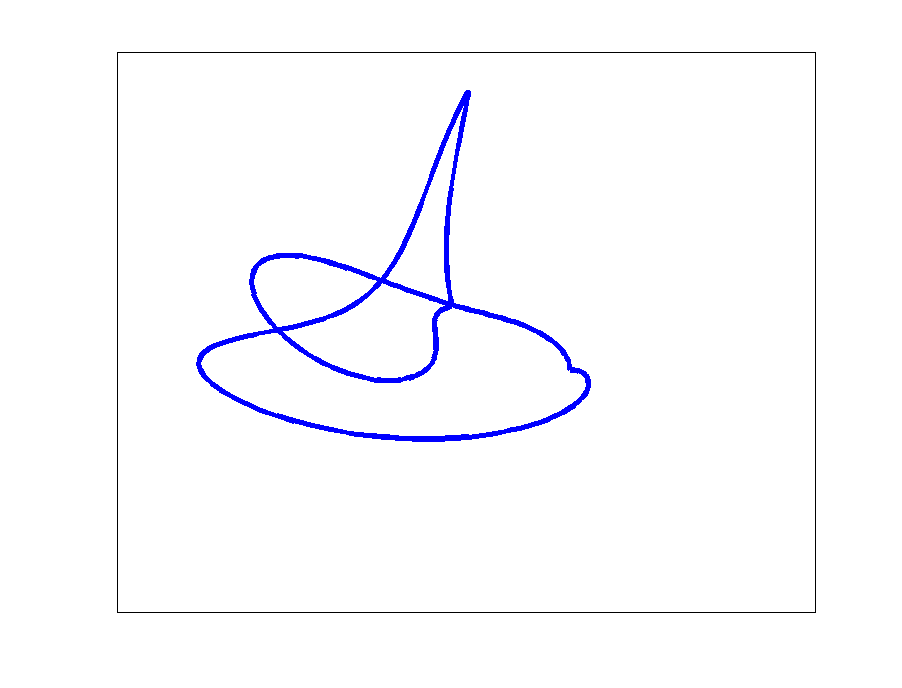}}
                \put(-50,0){\tiny{$\Re(A)$}}
                \put(-83,25){\rotatebox{90}{\tiny{$\Im(A)$}}}
                \label{fig:invariant_sol_t3}}
	  
        \vspace*{1.0ex}
        \subfloat[$|A|$]{
                {\includegraphics[trim = 22mm 12mm 15mm 10mm, clip,width=0.27\textwidth]{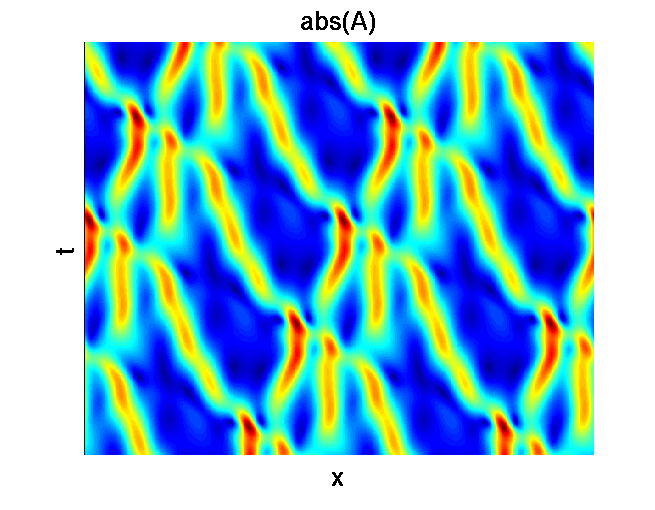}}
                \put(-60,-6){\small{$x$}}
                \put(-125,45){\rotatebox{90}{\small{$t$}}}
                \label{fig:invariant_sol_abs}}
        \caption{{Solution having symmetry  \refEqn{eqn:cgle_invariant_solution}, with $T> 0$ and nonzero
                 $\varphi$ and $S$.
                }}
        \label{fig:illustrate_invariant_sol}
\end{figure}
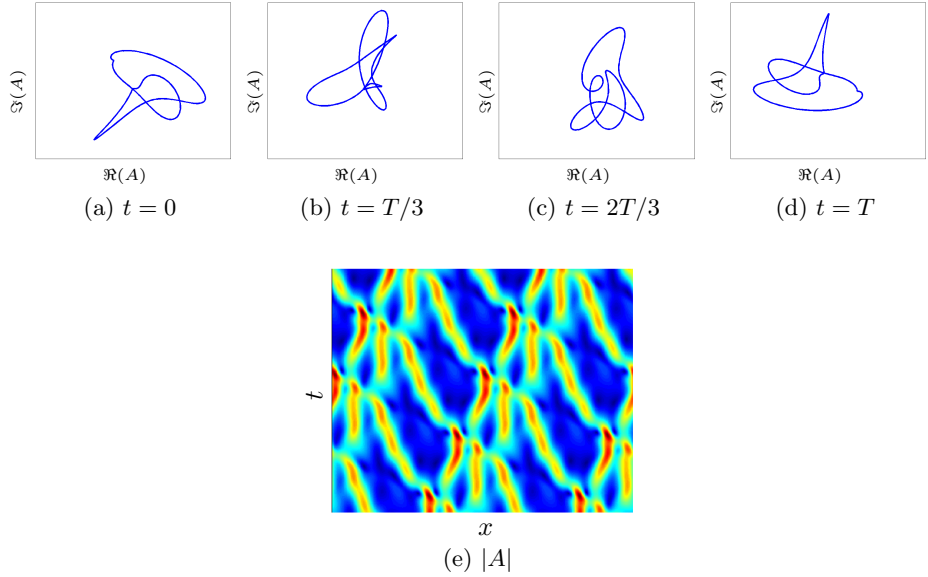
To illustrate some of the aforementioned symmetries of solutions of the {CGLE}, 
Figures~\refFigsSecs{fig:illustrate_invariant_sol}{fig:illustrate_additional_symm} display several plots that aid
in visualizing the invariant properties.  Figure~\ref{fig:illustrate_invariant_sol} shows a solution of the 
{CGLE} having symmetry \refEqn{eqn:cgle_invariant_solution}, but none of \refEqns{additionalsym1}{additionalsym3}.  
This solution belongs to the sequence $\Ai{5}$ (see \refEqn{eqn:Asequence}) resulting from the 
numerical continuation procedure to be described in Section~\ref{sec:methodology}, that is, from the continuation path 
for the sequence listed with id~5 in Tables~\refFigsSecs{tbl:properties_solutions}{tbl:properties_solutions_2} (refer to Section~\ref{sec:results}).  
The time evolution, represented as curves on the plane with coordinates defined by the real part $\Re(A)$ and 
imaginary part $\Im(A)$ of the solution $A(x,t)$ at different times within the interval $[0,T]$, is depicted in
Figures~\refFigsSecs{fig:invariant_sol_t0}{fig:invariant_sol_t3}.
The excitation of multiple temporal frequencies is apparent from these curves.  For single-frequency solutions
$A(x,t) = B(x) \e^{\ii\omega t}$ or generalized traveling waves $A(x,t) = \rho(x-vt) \e^{\ii\phi(x-vt)} \e^{\ii\omega t}$
(where $\omega$ is some single frequency), plots of this kind would show, except for a rotation, the same curve
at each point in time.  Therefore it is clear that the solution depicted in
Figures~\refFigsSecs{fig:invariant_sol_t0}{fig:invariant_sol_t3} is not of either of these single-frequency types.
Note also that the curve at time $t=T$ differs only by a rotation from that at time $t=0$ due to the rotation of the complex field
$A(x,0) \rightarrow \e^{\ii \varphi}A(x+S,T)$.
Repeated patterns resulting from invariance due to time periodicity and the nonzero space translation $S$ is 
better observed from surface plots of the absolute value $|A|$ of $A(x,t)$ over several space and time periods, 
as in Figure~\ref{fig:invariant_sol_abs}.

\begin{figure}[!t]  
        \centering
        \hspace*{-0.5ex}
        \subfloat[$\Re(A)$]{
                {\includegraphics[trim = 5mm 5mm 0mm 0mm, clip,width=0.38\textwidth]{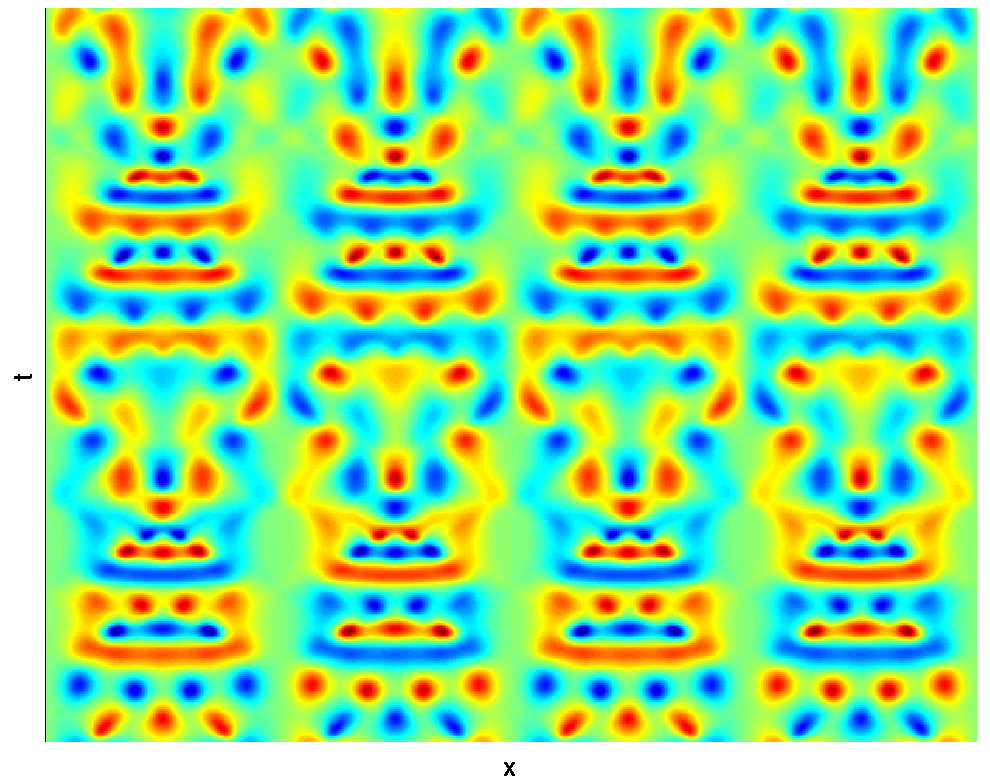}}
                \put(-80,-7){\small{$x$}}
                \put(-170,65){\rotatebox{90}{\small{$t$}}}
                \label{fig:additional_symm_real}}
        \hspace*{1.0ex}
        \subfloat[$\Im(A)$]{
                {\includegraphics[trim = 5mm 5mm 0mm 0mm, clip,width=0.38\textwidth]{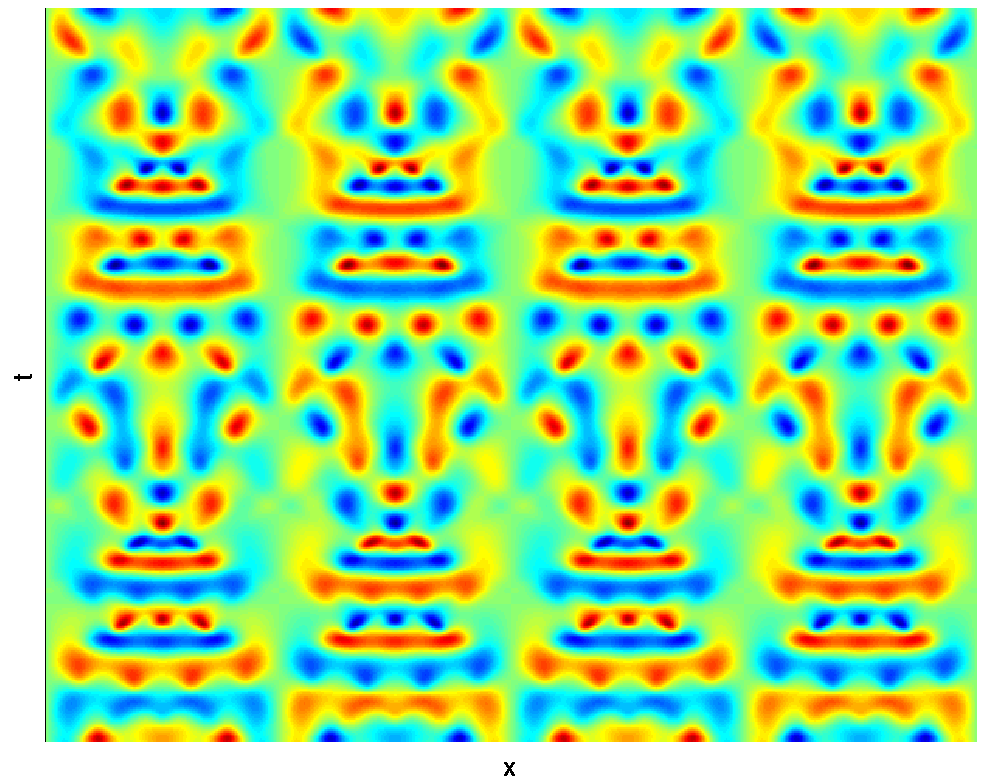}}
                \put(-80,-7){\small{$x$}}
                \put(-170,65){\rotatebox{90}{\small{$t$}}}
                \label{fig:additional_symm_imag}}
	    
        \subfloat[$|A|$]{
                {\includegraphics[trim = 5mm 5mm 0mm 0mm, clip,width=0.38\textwidth]{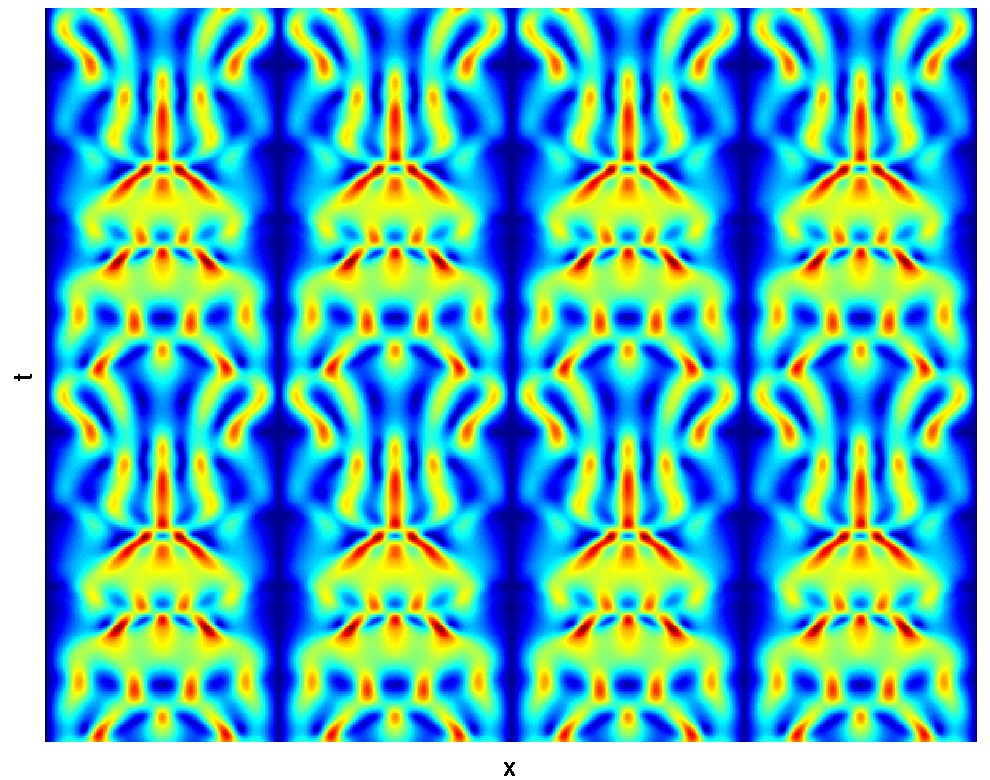}}
                \put(-80,-7){\small{$x$}}
                \put(-170,65){\rotatebox{90}{\small{$t$}}}
                \label{fig:additional_symm_abs}}
        \hspace*{1.0ex}
        \subfloat[$\Im(A) \mbox{ vs. } \Re(A),  \ t = 0, T/3, 2T/3, T$]{
                {\includegraphics[trim = 30mm 15mm 20mm 10mm, clip,width=0.38\textwidth]{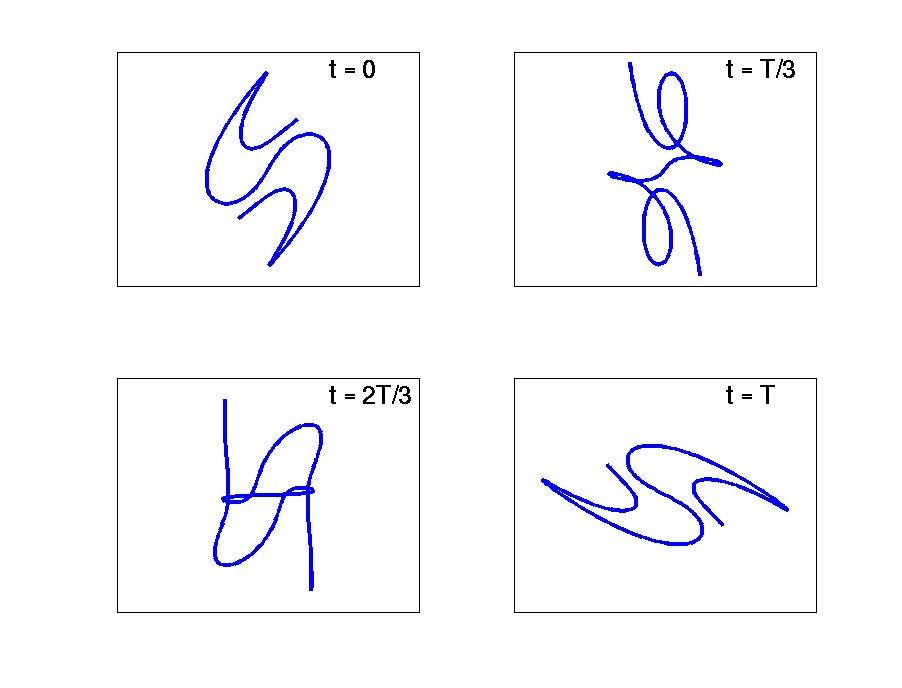}}
                \put(-80,-7){\small{$\ $}}
                \label{fig:additional_symm_time_evolution}}
        \caption{{Solution having symmetries \refEqn{additionalsym1}, for $\lsym=2$,
                 \refEqn{additionalsym2}, and \refEqn{additionalsym3}, in addition
                 to \refEqn{eqn:cgle_invariant_solution}. 
                }}
        \label{fig:illustrate_additional_symm}
\end{figure}

A solution possessing all of the symmetries \refEqns{additionalsym1}{additionalsym3}, in addition to the
symmetry \refEqn{eqn:cgle_invariant_solution}, is shown in Figure~\ref{fig:illustrate_additional_symm}.  
The plots represent a solution which belongs to the sequence of solutions under id~15 
in Tables~\refFigsSecs{tbl:properties_solutions}{tbl:properties_solutions_3} (refer to Section~\ref{sec:results}), computed at the point
$(R,\nu,\mu) = (100,-7,5)$ of the {CGLE} parameter space.  Surface plots of the real part $\Re(A)$, 
imaginary part $\Im(A)$, and absolute value $|A|$ of the solution $A(x,t)$ are shown in
Figures~\refFigsSecs{fig:additional_symm_real}{fig:additional_symm_abs}, where $(x,t) \in [0,2\Lx] \times [0,2T]$,
that is, the surfaces are plotted over two space and two time periods.  
For this solution, symmetry \refEqn{additionalsym1} holds with $\lsym=2$.
Since, in addition, the solution is even about $x = \tilde{m}\Lx/4$, for $\tilde{m} \in \integers$ odd, it follows from
\refEqn{eqn:lsym_even} that the solution is also odd about $x = (\tilde{m} + 1) \Lx/4$. 
Finally, the absolute value of the solution has spatial period $\Lx/2$ and is time-periodic, with period $T$.
As seen in Figures~\refFigsSecs{fig:additional_symm_real}{fig:additional_symm_time_evolution}, pattern similarities 
in both space and time are easily observed in the presence of the additional symmetries
\refEqns{additionalsym1}{additionalsym3}.

We conclude this section by noting the following fact.
Suppose that for some $(\tilde{\varphi}, c, \tilde{T}) \in \cglegroup$, where $\cglegroup$
is the group of continuous symmetries of the {CGLE} (refer to \refEqn{eqn:Ggroup}),
a solution $A(x,t)$ of the {CGLE} has the symmetry
\begin{equation}    \label{additionalsym4}
    A(x,t)  =   \e^{\ii \tilde{\varphi}} A(-x + c, t + \tilde{T}).
\end{equation}
The Fourier coefficient functions $\am{m}(t)$ in \refEqn{eqn:xFseries} of such a solution satisfy
 \begin{equation}   \label{eqn:Fcoefs_addsymm4}
    \am{-m}(t)   =  \am{m}(t + \tilde{T}) \,  \e^{\ii \tilde{\varphi}}\e^{\ii \km  c}   \ , \ \ m = 0, 1, 2, \ldots \, .
\end{equation}
The right-hand side of \refEqn{additionalsym4} is the result of the (left) action on $A(x,t)$ of the composition
$(\tilde{\varphi}, c, \tilde{T}) \circ (x \rightarrow -x)$, and after a (left) action of said composition on both sides
of \refEqn{additionalsym4} one obtains that
\begin{eqnarray} 
    A(x,t)  & = &  \e^{\ii \tilde{\varphi}} A(-x + c, t + \tilde{T})  \nonumber  \\
               & = &  \e^{\ii 2 \tilde{\varphi}} A(x, t + 2\tilde{T}).   \label{eqn:element_isotropy_G2}
\end{eqnarray}
Hence $(2 \tilde{\varphi}, 0, 2\tilde{T}) \in \cgleisotropy{A}$,
where $\cgleisotropy{A}$ is the isotropy subgroup defined in \refEqn{eqn:isotropy_subgroup}.
Also, note that we have
\begin{equation*}
    [  (\tilde{\varphi}, c, \tilde{T}) \circ (x \rightarrow -x)  ]^{2}  \ = \  (2 \tilde{\varphi}, 0, 2\tilde{T})
        \ \in \   \cglegroup,
\end{equation*} 
for any element $(\tilde{\varphi}, c, \tilde{T})$ in the group $\cglegroup$ of continuous symmetries of the {CGLE}.
Conversely, let $(\tilde{\varphi}, 0, \tilde{T}) \in \cgleisotropy{A}$ for some solution $A(x,t)$ of the {CGLE}.
Then we have
\begin{equation*} 
    [  (\tilde{\varphi}/2 + k \pi, \, c \, , \tilde{T}/2) \circ (x \rightarrow -x)  ]^{2}  \ = \  (\tilde{\varphi}, 0, \tilde{T})
\end{equation*} 
for every $k \in \integers$ and $c \in \reals$.
Therefore, solutions of the {CGLE} with symmetry \refEqn{additionalsym4} do possess symmetry 
\refEqn{eqn:cgle_invariant_solution} of the type we seek, and, conversely, a solution with symmetry
\refEqn{eqn:cgle_invariant_solution} may also possess the additional symmetry \refEqn{additionalsym4}.
An example of solutions having both symmetries \refEqn{eqn:cgle_invariant_solution} and \refEqn{additionalsym4}
is described in Section~\ref{sec:results} (cf.~Figure~\ref{fig:sol13_abs}).  Such solutions appeared in the continuation 
path for the sequence listed with id~13 in Tables~\refFigsSecs{tbl:properties_solutions}{tbl:properties_solutions_3} (refer to Section~\ref{sec:results}),
that is, in the sequence $\Ai{13}$ (see \refEqn{eqn:Asequence}).

Finally, note that if a solution $A(x,t)$ of the {CGLE} having symmetry \refEqn{eqn:cgle_invariant_solution}
for $S = 0$ or $S = \Lx/2$ also satisfies
\begin{equation}    \label{additionalsym5}
    A(x,t)   =    \e^{\ii \hat{\varphi}} A(-x + \hat{c}, t)
\end{equation}
for some real numbers $\hat{\varphi}$ and $\hat{c}$
(an example being solutions with the additional symmetry \refEqn{additionalsym2}, where $\hat{\varphi}$ = 0,
or with the additional symmetry \refEqn{additionalsym3}, for which $\hat{\varphi} = \pi$), then
\[   
    A(x,t)   \ = \    \e^{\ii (\hat{\varphi} + 2 \varphi)} A(-x + \hat{c}, t + 2T).   
\]
That is, such solution $A$ also has symmetry \refEqn{additionalsym4}.
Therefore, one should expect to find solutions having both symmetries \refEqn{eqn:cgle_invariant_solution} 
and \refEqn{additionalsym4} in subspaces of the space of solutions $( A; \varphi, S, T )$ for which
either $S = 0$, by \refEqn{eqn:element_isotropy_G2}, or $S \in \{0,\Lx/2\}$, if symmetry \refEqn{additionalsym5} 
is also present.

\section{Numerical Method}
\label{sec:methodology}

As noted in Section~\ref{sec:introduction}, having computed previously in \cite{lopez05} a set of unstable invariant 
solutions of the {CGLE} for fixed values of the parameters $(R,\nu,\mu)$, our first goal is to 
employ numerical continuation to carry solutions of this initial set into solutions in a regime with a different 
set of  parameter values $(R,\nu,\mu)$.  
To achieve this, we discretize using Fourier series expansions in both space and time to derive an underdetermined 
system of nonlinear algebraic equations from which invariant solutions of the {CGLE} are sought.
This discretization was used in the previous study \cite{lopez05}.  The associated material which is directly relevant to the current
study is summarized in Sections~\ref{sec:discretization_nleqns} and \ref{sec:discretization_jacobian} below in order to make
the present account self-contained.  Details concerning the numerical continuation, which was not a component of the 
previous study \cite{lopez05}, are provided in Section~\ref{sec:continuation_procedure}.

\subsection{Derivation of Nonlinear Algebraic Equations}
\label{sec:discretization_nleqns}

Since the boundary conditions \refEqn{eqn:cgle_bcs}  are periodic in $x$, 
we use the spatial Fourier series \refEqn{eqn:xFseries} and substitute into the {CGLE}~\refEqn{eqn:cgle_pde}
to obtain an infinite system of ordinary differential equations ({ODEs}),
\begin{equation}   \label{eqn:odes}
    \frac {\dd \am{m}}{\dd t}  =  R \am{m} - \km^2 (1 + \ii\nu) \am{m} - (1 + \ii\mu) \sum_{m_1+m_2-m_3=m}
	\am{m_1} \am{m_2} \am{m_3}^{*},
\end{equation}
for the complex-valued functions $\am{m}(t)$.  
Under this transformation the symmetries \refEqns{cglesym1}{cglesym4} of equations \refEqns{eqn:cgle_pde}{eqn:cgle_bcs}
become symmetries of~\refEqn{eqn:odes}.  Thus, if $\mmbf{a}(t) = (\am{m}(t))$ is a solution of the system of 
{ODEs} \refEqn{eqn:odes}, then so are 
\begin{align}   
    (\e^{\ii\Arot} \am{m}(t))&,                    \label{odesym1} \\
    (\e^{\ii m\spacetrans} \am{m}(t))&,   \label{odesym2} \\
    (\am{m}( t + \timetrans))&,                 \label{odesym3} \\
    (\am{-m}( t))&,                                     \label{odesym4} 
\end{align}
for any $(\Arot, \spacetrans, \timetrans) \in  \torus^2 \times \reals$.
In particular, \refEqn{odesym1} and \refEqn{odesym2} say that the {ODEs}~\refEqn{eqn:odes} are invariant
under the $\torus^2$-action
\begin{equation*}  \label{eqn:torusact}
    (\Arot, \spacetrans)\cdot (\am{m}(t))  = (\e^{\ii\Arot} \e^{\ii m \spacetrans} \am{m}(t)).
\end{equation*}

We employ a spectral-Galerkin projection obtained by fixing an even
number $N_x$ and truncating the expansion \refEqn{eqn:xFseries} to include only the terms with indices $m$ satisfying 
$-N_x/2+1 \leq m \leq N_x/2-1$.  We then work with the corresponding finite system of {ODEs} which results from \refEqn{eqn:odes}
after the Galerkin projection.  Much
accumulated theory and computation shows that for sufficiently large $N_x$ the behavior of this truncation
captures the essential features of the dynamics of~\refEqns{eqn:cgle_pde}{eqn:cgle_bcs} \cite{Doelman, Jolly}.

From the condition \refEqn{eqn:cgle_invariant_solution} defining an invariant solution of the {CGLE}, 
it follows that the corresponding solution $\mmbf{a}(t)$ of the system of {ODEs} \refEqn{eqn:odes} satisfies 
\begin{equation}   \label{eqn:am_rpo}
    \am{m}(t)  =  \e^{\ii\varphi} \e^{\ii \km S} \am{m}(t +T)
\end{equation}
for all $m$ and $t$ (and where $\varphi, S, T$ are to be determined).  
It is easy to see that the set of functions 
\begin{equation}   \label{eqn:am_ansatz}
    \am{m}(t)  =  \e^{-\ii\frac{\varphi}{T}t} \e^{-\ii \km \frac{S}{T}t}
            \sum_{n \in \integers} \amn{m}{n} \e^{\ii \omegan t} \, ,
\end{equation}
where $\omegan = 2 \pi n / T$ denotes the $n$-th frequency in the expansion, are a solution of the system of 
functional equations \refEqn{eqn:am_rpo}. Hence, they provide an appropriate
representation for invariant solutions of the system of {ODEs} \refEqn{eqn:odes}. 
Substituting \refEqn{eqn:am_ansatz} into the truncated system of {ODEs} \refEqn{eqn:odes} and using again
a Galerkin projection obtained by fixing an even number $N_t$, so that the summation
index in~\refEqn{eqn:am_ansatz} runs over the range $-N_t/2+1 \le n \le N_t/2-1$, 
results in a system of nonlinear algebraic equations,
\begin{equation}  \label{eqn:cgle_nleqns}
    \Fargs  \: = \:  \FLargs + \FNLargs  \: = \:  \mmbf{0},
\end{equation}
for the complex Fourier coefficients $\{\amn{m}{n}\}$ and elements $(\varphi, S, T)$ of the isotropy subgroup
\refEqn{eqn:isotropy_subgroup}.  In \refEqn{eqn:cgle_nleqns}, 
$\mmbf{\hat{a}}$ denotes a vector with components given by the coefficients $\{\amn{m}{n}\}$
and the vectors $\FLargs$ and $\FNLargs$ are defined as
\begin{equation}  \label{eqn:Flin}
   \FLargs  \ \equiv \  \left \{ \ii \left (\frac{2\pi n}{T} - \frac{\varphi}{T} - \km \frac{S}{T} \right )	\amn{m}{n}
        - R \amn{m}{n} + \km^2 (1 + \ii\nu) \amn{m}{n} \right \}  \, ,
\end{equation}
and
\begin{equation}  \label{eqn:Fnonlin}
    \FNLargs  \ \equiv \   \left \{ (1 + \ii\mu) \sum_{m_1+m_2-m_3=m} \left ( \sum_{n_1+n_2-n_3=n}
    \amn{m_1}{n_1} \amn{m_2}{n_2} \amn{m_3}{n_3}^{*} \right ) \right \} \,. 
\end{equation}
Note that the components of the vector $\FL$  in \refEqn{eqn:Flin} correspond to the discretization of the linear terms in 
the {CGLE} and those of $\FNL$  in \refEqn{eqn:Fnonlin} to that of the nonlinear term $(1 + \ii\mu) A|A|^2$.
Furthermore, in defining the vector
$\mmbf{\hat{a}}$ (and similarly for $\FL$ and $\FNL$) we are implicitly assigning
an ordering on the coefficients $\{\amn{m}{n}\}$ that uniquely determines an indexing for the components
of $\mmbf{\hat{a}}$.  Henceforth, such a convention should be understood whenever applicable.
Finally, we will use the notation in~\refEqn{eqn:cgle_nleqns} to denote both the system of complex
equations and the system obtained by splitting \refEqn{eqn:cgle_nleqns} into its
real and imaginary parts, as it should be clear from the context which case applies.

Splitting the equations into their real and imaginary parts, one has that  \refEqn{eqn:cgle_nleqns} is an underdetermined 
system of $2(N_x-1)(N_t-1)$ real equations in $2(N_x-1)(N_t-1) + 3$ real unknowns.  Solutions of this system of equations 
will give the desired invariant solutions of the truncated system of {ODEs} via the expansion \refEqn{eqn:am_ansatz}.  
We note here that with the introduction of the representation \refEqn{eqn:am_ansatz},
the symmetry group $\cglegroup = \torus^2 \times \reals$ of \refEqn{eqn:cgle_pde} and \refEqn{eqn:odes} descends
to the symmetry group $\torus^3 = \torus^2 \times \mathrm{S}^1$ of \refEqn{eqn:cgle_nleqns}, acting on the space
$(\{\amn{m}{n}\}, \varphi,S,T)$ of solutions of \refEqn{eqn:cgle_nleqns}.  Henceforth, by a slight abuse of notation, 
we refer to both symmetry groups $\torus^2 \times \reals$ and $\torus^3$ as $\cglegroup$.

The symmetries~\refEqns{odesym1}{odesym4} of the {ODEs} \refEqn{eqn:odes}
induce symmetries of the system of algebraic equations \refEqn{eqn:cgle_nleqns}.  Note that if 
$(\{\amn{m}{n}\}, \varphi,S,T)$ solves \Feqzero, then for any $(\Arot, \spacetrans, \timetrans) \in \torus^3$
\begin{align}
    (\{ \e^{\ii\Arot} \amn{m}{n} \}, \varphi, S, T)&,                    \label{Fsym1}    \\
    (\{ \e^{\ii m\spacetrans} \amn{m}{n} \}, \varphi, S, T)&,    \label{Fsym2}     \\
    (\{ \e^{\ii n \timetrans} \amn{m}{n} \}, \varphi, S, T)&,       \label{Fsym3} 
\end{align}
\begin{align}
    (\{ \amn{-m}{n} \}, \varphi, -S, T)&,                                     \label{Fsym4} 
\end{align}
are also solutions.  From the continuous symmetries \refEqns{Fsym1}{Fsym3}, it follows that the set of
solutions of \Feqzero\ splits into orbits $\mathcal{O}_{(\nleargs)}$ of the symmetry group $\torus^3$,
\begin{equation}    \label{eqn:torus3orbits}
    \mathcal{O}_{(\nleargs)}  :=  \left \{ (\Arot,\spacetrans,\timetrans) \cdot 
    (\nleargs) \; | \: (\Arot,\spacetrans,\timetrans) \in \torus^3 \right \},
\end{equation}
where the action of $\torus^3$ on a point $(\nleargs)$ is defined by
\begin{equation}  \label{eqn:torus3act}
     (\Arot,\spacetrans,\timetrans) \cdot (\nleargs)  \ = \ (\{ \e^{\ii\Arot} \e^{\ii m \spacetrans}
     \e^{\ii n \timetrans}  \amn{m}{n} \}, \varphi, S, T). 
\end{equation}
That is, $\torus^3$ acts on $\mmbf{\hat{a}}$ via multiplication by the matrix $\diag(\e^{\ii\Arot}
\e^{\ii m \spacetrans} \e^{\ii n \timetrans})$, and it acts trivially on $(\varphi, S, T)$. 
Finally we note that, for the system of nonlinear algebraic equations \refEqn{eqn:cgle_nleqns},
the transformation \refEqn{Fsym4} induced by \refEqn{cglesym4} maps a solution
\begin{equation}    \label{Fsym5a}
    (\{ \amn{m}{n} \}, \varphi, \Lx/2 \pm \delta, T)
\end{equation}
of \refEqn{eqn:cgle_nleqns} to another (conjugate) solution
\begin{equation}    \label{Fsym5b}
    (\{ \amn{-m}{n+m} \}, \varphi, \Lx/2 \mp \delta, T) \, 
\end{equation}
of \refEqn{eqn:cgle_nleqns}, where, again, $\delta = |\Lx/2 - S|$. 
(Refer to the paragraph containing \refEqn{eqn:Smap}.)

\subsection{Jacobian Matrix}
\label{sec:discretization_jacobian}

The \jacobian\ matrix of the system $\Feqzero$ of nonlinear algebraic equations \refEqn{eqn:cgle_nleqns} is dense so,
as the number of unknowns (and equations) increases,
it becomes unfeasible to solve linear systems with the \jacobian\ as coefficient matrix using direct methods.
However, matrix-vector products with the \jacobian\ matrix of $\F$ can be computed efficiently for the problem at hand,
making the use of iterative methods for solving linear systems a viable option.
We proceed to review the calculation of this matrix-vector product since it is an essential feature of the Newton step 
computation employed in the numerical continuation.

Let $\Ja$ denote the matrix whose columns correspond to derivatives of $\F$ with respect to the real and imaginary 
parts of the unknowns $\{\amn{m}{n}\}$, and let $\hat{\mmbf{v}}$ be a vector with components given by the coefficients 
in the truncated Fourier series expansion of a function $V(x,t)$.  Assume that $\Ja$ is evaluated at a given
point $(\nleargs)$.  The product $\productJa{v}$ can then be computed as\footnote{Note that in the right-hand side
of \refEqn{eqn:Jv} we are actually using $\hat{\mmbf{v}}$ to denote a vector with the \emph{complex} numbers 
$\{\hat{v}_{m,n}\}$ as components, whereas in the left-hand side of \refEqn{eqn:Jv} $\hat{\mmbf{v}}$ denotes a vector 
with \emph{real} components that are the real and imaginary parts of the coefficients $\{\hat{v}_{m,n}\}$.  We make
use of this slight abuse of notation in this paper since the intended meaning should be clear from the context.}
\begin{equation}   \label{eqn:Jv}
    \productJa{v}  = 
    \JL(\hat{\mmbf{v}}, \varphi, S, T) + 
    \JNL(\hat{\mmbf{a}}, \hat{\mmbf{v}}),
\end{equation}
where $\JL(\hat{\mmbf{v}}, \varphi, S, T) \equiv \FL(\hat{\mmbf{v}}, \varphi, S, T)$ (as defined in \refEqn{eqn:Flin})
and $\JNL(\hat{\mmbf{a}}, \hat{\mmbf{v}})$ is a vector with components given by the coefficients in the truncated 
Fourier series expansion (in both space and time) of \mbox{$(1+i\mu) (A^2V^{*} + 2|A|^2V)$}.
This matrix-vector product operation follows from the discretization (analogous to that used for the {CGLE}) of the 
first variational derivative of equation \refEqn{eqn:cgle_pde},
\begin{equation*}
    \frac {\partial V}{\partial t}   =   R V + (1 + i \nu)
    \frac{\partial^2 V}{\partial x^2} - (1 + i \mu) (A^2V^{*} + 2|A|^2V).
\end{equation*}
Furthermore, as can be seen from system \refEqns{eqn:cgle_nleqns}{eqn:Fnonlin}, the operation of computing a matrix-vector 
product with the columns of the \jacobian\ matrix of $\F$ corresponding to the derivatives with respect to 
$\varphi$, $S$, and $T$ poses no difficulty. 

It follows then that matrix-vector products with the \jacobian\ matrix of $\F$ can be easily computed without the need of 
explicitly calculating the (full) \jacobian.  Note also from \refEqn{eqn:Flin} that the portion of the \jacobian\ matrix coming from 
the discretized linear terms \FLargs\ in the {CGLE} is a block diagonal matrix, with $2 \times 2$ blocks, whose components 
are easily computed.  Hence, solving linear systems with this block diagonal matrix poses no complications.  This is 
advantageous since this block diagonal matrix provides an effective preconditioner for some matrix-free iterative 
methods when solving linear systems having the \jacobian\ as coefficient matrix for the problem at hand.  
(Refer to Section \ref{sec:continuation_procedure}.)

Finally, we note that the matrix $\Ja$ (refer to \refEqn{eqn:Jv}),
whose columns correspond to derivatives of $\F$ with respect to the real and imaginary parts of the 
unknowns $\{\amn{m}{n}\}$, is singular at a solution $(\nleargs)$ of \Feqzero.  
This is relevant for the computation of the Newton step, discussed in Appendix~\ref{app:newton_step}.
The vectors in the null space of $\Ja$ result
from a basis for the space of infinitesimal generators of the action~\refEqn{eqn:torus3act}
of $\torus^3$ on the point $(\nleargs)$.  The reader may consult \cite{lopez05} for further details.

\subsection{Numerical Continuation of Solutions}
\label{sec:continuation_procedure}

The numerical continuation was done using the Library of Continuation Algorithms ({LOCA}) software
package \cite{loca}, specifically with the aid of the algorithms provided to track steady state solutions 
of discretized {PDEs} as a function of a single parameter.  
The option of pseudo arc length continuation was used in order to allow for turning points \cite{bk:allgower} to be followed.  
Although we are not computing steady state solutions in this study, 
it is clear that the feature of tracking steady state solutions in the {LOCA} package 
provides the capability of solving a system of nonlinear algebraic equations using numerical continuation
(which is what we need). 
We thus take advantage of this feature, particularly to handle the step size control, 
that is, changes in the value of the continuation parameter, including that in the vicinity of turning points, at
each continuation step.
A general description of the continuation procedure appears next, along with details concerning the input required to be
supplied by the user to the {LOCA} routines.  For specific information on the implementation of capabilities used as provided 
by the {LOCA} package (that is, without us making any modifications to the {LOCA} software), like that of the computation of 
changes in the value of the continuation parameter, the user is referred to the {LOCA} documentation \cite{loca}.

Let $\contParam{}$ denote the continuation parameter.  Since we perform single-parameter continuation,
$\contParam{}$ will correspond to one of the {CGLE} parameters $R$, $\nu$, or $\mu$.
Set $\uAlgo{} = (\hat{\mmbf{a}}, \varphi, S, T)$ and let $\fAlgo(\uAlgo{}; \contParam{})  = \mmbf{0}$ denote the system
$\Fargs  =  \mmbf{0}$ of nonlinear algebraic equations \refEqn{eqn:cgle_nleqns}, for a particular point
$(R,\nu,\mu)$ in the {CGLE} parameter space.  (Note that the point $(R,\nu,\mu)$ is associated with the continuation
parameter $\contParam{}$.)
The continuation process can then be described generally as follows:
\begin{myenumerate}
     \item Set the initial values of the {CGLE} parameters $(R,\nu,\mu)$ and solution $\uAlgo{0}$.
     \item Select one of the {CGLE} parameters $R$, $\nu$, or $\mu$ to be used as the continuation parameter  
              and set the initial value $\contParam{0}$ of the continuation parameter.
     \item Set the desired final value $\contParam{final}$ of the continuation parameter.
     \item Set the remaining inputs to the {LOCA} software \cite{loca}.  
              These include values for the 
              minimum $\stepChange{\mathrm{min}}$ and maximum $\stepChange{\mathrm{max}}$ 
              changes allowed in the continuation parameter at each continuation step.   
              Specific details appear in Section~\ref{sec:comments_numerics}, where
              comments on aspects related to the numerical simulations are provided.
     \item Set the maximum number $\contStep{max}$ of continuation steps.
     \item Begin loop:  For $\contStep{} = 1, 2, \ldots, \contStep{max}$
              \begin{myenumerate}
                   \item Determine the change $\stepChange{\contStep{}}$  in the value of the 
                            continuation parameter \cite{loca}.
                   \item Update the value of the continuation parameter:
                            $\contParam{\contStep{}} = \contParam{\contStep{}-1} + \stepChange{\contStep{}}$.
                   \item Solve the system $\fAlgo(\uAlgo{\contStep{}}; \contParam{\contStep{}})  = \mmbf{0}$  
                            for $\uAlgo{\contStep{}}$, providing $\uAlgo{\contStep{}-1}$ as the initial
                            guess for the nonlinear equations solver.  Details on the numerical solution of 
                            $\fAlgo  = \mmbf{0}$
                            are given in the next paragraph below.
                   \item If the nonlinear equations solver fails to converge, decrease the magnitude of $\stepChange{\contStep{}}$.
                            If $|\stepChange{\contStep{}}| < |\stepChange{\mathrm{min}}|$, exit the loop, indicating failure.
                            Otherwise, reset $\contParam{\contStep{}} = \contParam{\contStep{}-1} + \stepChange{\contStep{}}$ 
                            and go to step (c) above.
                   \item If $\contParam{\contStep{}} = \contParam{final}$, exit the loop, indicating convergence to a solution
                            at the final value $\contParam{final}$ of the continuation parameter.
              \end{myenumerate}
              End loop
\end{myenumerate}

Newton's method is used to solve the system of nonlinear algebraic equations in the numerical continuation, 
and the user must supply the {LOCA} package with a routine for computing the Newton step.  That is, the user must 
provide a routine that solves a linear system having as coefficient matrix the \jacobian\ of the system of nonlinear 
algebraic equations.   For this purpose, we employed iterative methods for solving linear systems,
specifically the {GMRES} solver from the Meschach software package \cite{meschach}.
The computation of the nonlinear terms 
$\FNL$ in \refEqn{eqn:cgle_nleqns} and $\JNL$ in \refEqn{eqn:Jv}, needed, respectively, for the evaluation of $\F$ and 
that of the product of the \jacobian\ matrix and a vector, was done using the {FFTW} software package \cite{fftw}.
For further efficiency in the calculations we used {POSIX} threads (pthreads) programming \cite{pthreads} in our
routines, taking thus advantage of the multiple cores available nowadays in personal workstations.

Rather than augmenting the system \refEqn{eqn:cgle_nleqns} with an additional set of equations in order to work with 
an equal number of equations and unknowns \cite{lopez05}, we work with the 
underdetermined system \refEqn{eqn:cgle_nleqns} and consider here a Newton step defined from the Moore-Penrose 
inverse \cite{bk:benisrael}.  This yields a minimum norm solution of the system of linear equations with the 
underdetermined \jacobian\ as coefficient matrix, and is one approach used in numerical continuation 
methods \cite{bk:allgower,wulff06}.  
A detailed description of the computation of the Newton step for the present study appears in Appendix~\ref{app:newton_step}.
We note that the use of conceptually simple techniques led to an accurate and efficient computation of the Newton step.
The techniques employed made it practical to solve a computationally challenging problem without the need of
a cluster or supercomputer.

\section{Numerical Study and Results}
\label{sec:results}

The procedure described in Section~\ref{sec:methodology} was applied to a subset of the unstable invariant solutions 
of the {CGLE} computed at the point $(R, \nu, \mu) = (16, -7, 5)$ of the {CGLE} parameter space (without employing
continuation) in the preceding  study \cite{lopez05}, in order to carry them into solutions in a regime with an increased 
value of the parameter $R$, namely to the region at the point $(R, \nu, \mu) = (100, -7, 5)$ of the {CGLE} parameter space.  
As indicated in Section~\ref{sec:introduction}, chaotic behavior is exhibited both at the initial and final parameter regions.

A summary of the obtained results is gathered in Tables~\refFigsSecs{tbl:properties_solutions}{tbl:properties_solutions_3} 
and Figure~\ref{fig:unstable_dim}.  Already from them, we see that our probe into the moduli space  
of $\cglegroup$-orbits reveals a complicated and interesting structure. To start with, along each of the continuation
paths $\Ai{i}$, $i = 1,\ldots,15$, we were able to compute a number $\Ni$ 
(listed in Table~\ref{tbl:properties_solutions_2}) of new distinct $\cglegroup$-orbits of solutions of the {CGLE} 
(each one of which corresponds to a distinct invariant solution of the {CGLE}).
Thus, the continuation paths $\Ai{i}$, $i = 1,\ldots,15$, can be thought of as (discrete) sections
of the fibered space of $\cglegroup$-orbits over the space of parameters of the {CGLE}.

Before describing the content of  Tables~\refFigsSecs{tbl:properties_solutions}{tbl:properties_solutions_3}, let us list the
possibilities that may occur when numerically continuing a set of distinct $\cglegroup$-orbits.
(These are analogous to the cases listed later in Section~\ref{sec:discussion_particular_solutions} where we examine  
continuation paths which revisit a fiber over a point in the {CGLE} parameter space 
after a series of steps while performing continuation for a single $\cglegroup$-orbit.)
Suppose that $\Ap{\pts{p}{0}{(i)}} \in \Ai{i}$ and $\Ap{\pts{p}{0}{(j)}} \in \Ai{j}$, $i\ne j$, are two invariant solutions representing 
distinct $\cglegroup$-orbits at the initial point $\pts{p}{0}{(i)} = \pts{p}{0}{(j)} = (\pt{R}{0}, \pt{\nu}{0}, \pt{\mu}{0})$ in the
{CGLE} parameter space, where $\Ai{i}$ and $\Ai{j}$ are, respectively, the continuation paths \refEqn{eqn:Asequence}
emanating from each one of the two initial invariant solutions.  Given two points $\pts{p}{k}{(i)}$ and $\pts{p}{l}{(j)}$ in the {CGLE}
parameter space and two invariant solutions $\Ap{\pts{p}{k}{(i)}} \in \Ai{i}$ and $\Ap{\pts{p}{l}{(j)}} \in \Ai{j}$,
it may happen that $\pts{p}{k}{(i)} = \pts{p}{l}{(j)}$, and we have to consider several possibilities.  Namely, whether
the invariant solutions $\Ap{\pts{p}{k}{(i)}}$ and $\Ap{\pts{p}{l}{(j)}}$ represent 
(i) the same $\cglegroup$-orbit,
that is, $(\varphi(\pts{p}{k}{(i)}), S(\pts{p}{k}{(i)}), T(\pts{p}{k}{(i)})) =  (\varphi(\pts{p}{l}{(j)}), S(\pts{p}{l}{(j)}), T(\pts{p}{l}{(j)}))$
and there exists some $(\Arot,\spacetrans,\timetrans) \in \cglegroup$ such that 
$\Ap{\pts{p}{k}{(i)}} = (\Arot,\spacetrans,\timetrans) \cdot \Ap{\pts{p}{l}{(j)}}$;
(ii) different, but conjugate, $\cglegroup$-orbits, as defined in Section~\ref{sec:problem_statement}
(see also \refEqns{Fsym5a}{Fsym5b}); or
(iii) different, non-conjugate, $\cglegroup$-orbits.

After a careful analysis of all computed solutions we found that there were solutions $\Ap{\pts{p}{k}{(2)}} \in \Ai{2}$ 
and $\Ap{\pts{p}{l}{(4)}} \in \Ai{4}$
which represent the same $\cglegroup$-orbit, at points $\pts{p}{k}{(2)} = \pts{p}{l}{(4)}$ in the range
$(R,\nu,\mu) \in ([16.5,54.4] \cup [54.7,100]) \times [-7] \times [5]$ of {CGLE} parameter values.
Furthermore, there were 
(a) solutions $\Ap{\pts{p}{k}{(11)}} \in \Ai{11}$ and $\Ap{\pts{p}{l}{(14)}} \in \Ai{14}$
which represent the same $\cglegroup$-orbit, at points $\pts{p}{k}{(11)} = \pts{p}{l}{(14)}$ in the range
$(R,\nu,\mu) \in [11.5,60]  \times [-7] \times [5]$;
(b) solutions $\Ap{\pts{p}{k}{(11)}} \in \Ai{11}$, $\Ap{\pts{p}{l}{(14)}} \in \Ai{14}$, and $\Ap{\pts{p}{j}{(15)}} \in \Ai{15}$
which represent the same $\cglegroup$-orbit, at points $\pts{p}{k}{(11)} = \pts{p}{l}{(14)} = \pts{p}{j}{(15)}$ in the range
$(R,\nu,\mu) \in [16,60]  \times [-7] \times [5]$;
(c) solutions $\Ap{\pts{p}{l}{(14)}} \in \Ai{14}$ and $\Ap{\pts{p}{j}{(15)}} \in \Ai{15}$
which represent the same $\cglegroup$-orbit, at points $\pts{p}{l}{(14)} = \pts{p}{j}{(15)}$ in the range
$(R,\nu,\mu) \in [16,85]  \times [-7] \times [5]$; and
(d) solutions $\Ap{\pts{p}{k}{(5)}} \in \Ai{5}$ and $\Ap{\pts{p}{l}{(6)}} \in \Ai{6}$
which represent conjugate $\cglegroup$-orbits, at points $\pts{p}{k}{(5)} = \pts{p}{l}{(6)}$ in the range
$(R,\nu,\mu) \in [17,59.5]  \times [-7] \times [5]$. 
The latter case is illustrated in Figure~\ref{fig:several_RT_paths}, where the graphs in $(R,T)$-space
for the sequences $\Ai{5}$ and  $\Ai{6}$ overlap for the aforementioned range of $R$. 
Therefore,  multiple representatives of same $\cglegroup$-orbits were
carefully accounted for and only one of them was taken as representative of the corresponding distinct $\cglegroup$-orbit.

Figure~\ref{fig:several_RT_paths} illustrates as well that, while for the sequences $\Ai{1}$ and $\Ai{9}$, for example, the
numerical continuation progressed in a relatively smooth manner, such was not the case in general.
Turning points and overlapping paths, exemplified by the depiction of the graphs $T=T(R)$ for sequences $\Ai{5}$, $\Ai{6}$, 
and $\Ai{8}$ in Figure~\ref{fig:several_RT_paths}, were frequently encountered. 
These features revealed intricate and challenging parameter regions for traversal. 
(Details appear in Section~\ref{sec:discussion_particular_solutions} below.)
\begin{figure}[!t]  
        \centering
        \subfloat[continuation paths in $(R,T)$ space]{
                {\includegraphics[trim = 10mm 0mm 0mm 0mm, clip, width=0.5\textwidth]{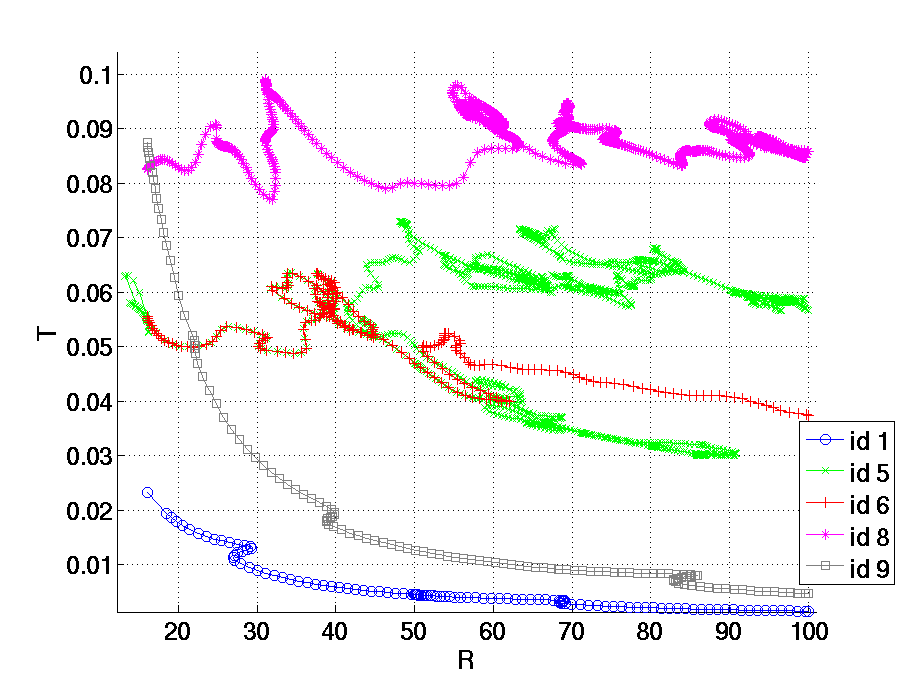}}
                \label{fig:several_RT_paths}}
        \subfloat[unstable dimension]{
                {\includegraphics[trim = 10mm 0mm 5mm 0mm, clip, width=0.5\textwidth]{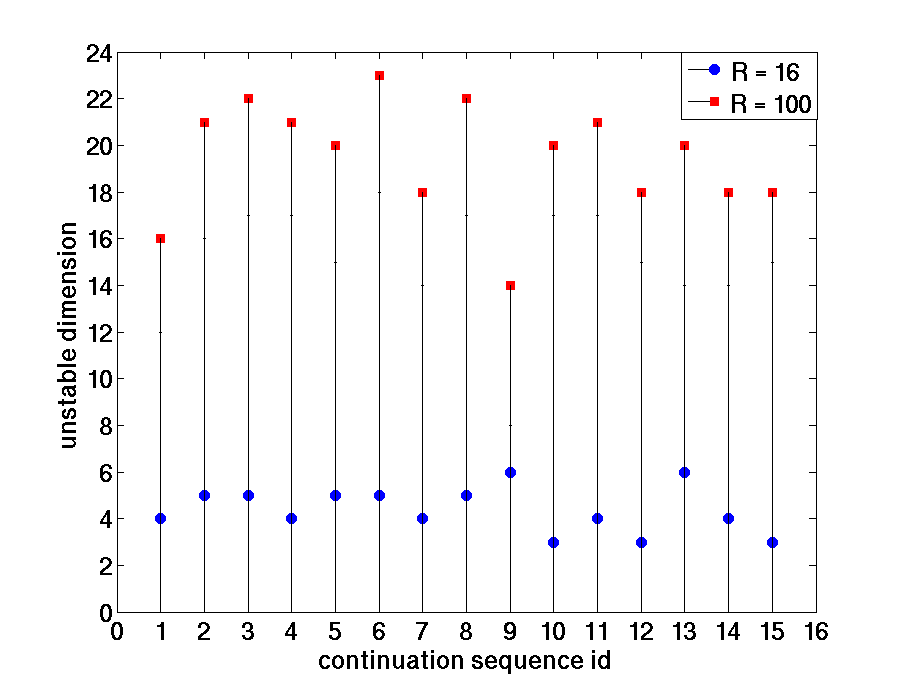}}
                \label{fig:unstable_dim}}
        \caption{{
                       (a) Representative continuation paths, depicted by plotting the time period $T$ as a function
                             of the parameter $R$.
                       (b) Unstable dimension of solutions at the initial and final parameter regions
                             for each sequence $\Ai{i}$, $i = 1,\ldots,15$.
                      }}
        \label{fig:summary_solns}
\end{figure}

Table~\ref{tbl:properties_solutions} lists the values  of $(\varphi,S,T)$ for the solutions at the starting {CGLE} 
parameter values of $(R, \nu, \mu) = (16, -7, 5)$ and at the final values of $(R, \nu, \mu) = (100, -7, 5)$, as well the 
unstable dimension\footnote{The unstable dimension of an invariant (or relative time-periodic) solution is the number of
eigenvalues of the associated relative monodromy matrix having magnitude greater than one; see \cite{lopez05}.} 
and spatial period of the invariant solutions at the aforementioned parameter values.  
Per the third column in Table~\ref{tbl:properties_solutions}, the listed solutions are unstable.  
As seen from Table~\ref{tbl:properties_solutions} and the depiction in Figure~\ref{fig:unstable_dim}, the solutions 
used as initial points for the numerical continuation have unstable dimension  ranging between 3 and 6, whereas the 
new solutions in the final parameter region have unstable dimension between 14 and 23.
The time period for the initial solutions is in the range $T \in (0.02, 0.12)$ (or $T \in (0.32,1.92)$ for the 
formulation \refEqn{eqn:cgle_pde2} of the {CGLE}); for the new solutions in the final parameter region 
we have $T \in (0.001, 0.11)$ (or $T \in (0.1,11)$ for the formulation \refEqn{eqn:cgle_pde2} of the {CGLE}).
No truly time-periodic solutions were identified (although their existence in the regions traversed is not ruled out), 
as all solutions have a nonzero value for the rotation angle $\varphi$.

Except for the sequence $\Ai{1}$, listed with id~1 in Tables~\refFigsSecs{tbl:properties_solutions}{tbl:properties_solutions_2},
for which the solution at the final
parameter values has only a few temporal frequencies active and appears to be close to a single-frequency solution, 
all of the solutions have broad spatial and temporal spectra.
Also, aside from the sequences $\Ai{2}$ and $\Ai{4}$, listed, respectively, with ids~2 and 4 in
Tables~\refFigsSecs{tbl:properties_solutions}{tbl:properties_solutions_3}, all of the resulting solutions retained the 
same spatial period of length $\Lx$ as that of the starting solutions.  The spatial period $\Lx/3$ of solutions in the
sequences $\Ai{2}$ and $\Ai{4}$ was acquired (for both sequences) at parameter values $(R,\nu,\mu) \approx (20.2,-7,5)$. 
The ending solutions in these two sequences are different elements of the same orbit \refEqn{eqn:torus3orbits}  
of the symmetry group $\cglegroup$ at the final point in parameter space, although the corresponding starting solutions belong 
to different orbits.  As for the other sequences, the solutions at the final point in parameter space belong to different 
orbits of the symmetry group $\cglegroup$.  
\begin{table}[!t]  
\tablesize
\centering
\begin{tabular}{|r|c|c|l|}
\hline
        &        &    unstable   &  \multicolumn{1}{c|}{spatial}   \\
id      & $(\varphi,\,S,\,T)_{\sss {(R,\nu,\mu) = (16,-7,5)}}  
              \rightarrow (\varphi,\,S,\,T)_{\sss {(R,\nu,\mu) = (100,-7,5)}}$
	         &    dimension   &  \multicolumn{1}{c|}{period}  \\
\hline  \hline
1  &  (5.3622, 3.8544, 0.0233) \ra (5.9158, 3.8856, 0.0015)  &  4 \ra 16  &  $\Lx$ \ra $\Lx$  \\
2  &  (2.8849, 3.0956, 0.0539) \ra (0.1088, 3.1416, 0.0130)  &  5 \ra 21  &  $\Lx$ \ra $\Lx / 3$  \\
3  &  (0.0011, 3.9709, 0.0539) \ra (5.5905, 2.2876, 0.0193)  &  5 \ra 22  &  $\Lx$ \ra $\Lx$  \\
4  &  (2.9343, 3.1416, 0.0540) \ra (0.1088, 3.1416, 0.0130)  &  4 \ra 21  &  $\Lx$ \ra $\Lx / 3$  \\
5  &  (4.6093, 1.4537, 0.0547) \ra (0.9483, 1.0333, 0.0567)  &  5 \ra 20  &  $\Lx$ \ra $\Lx$  \\
6  &  (4.5165, 4.7061, 0.0556) \ra (5.2620, 5.1417, 0.0374)  &  5 \ra 23  &  $\Lx$ \ra $\Lx$  \\
7  &  (0.2436, 2.3887, 0.0608) \ra (4.0066, 3.1416, 0.0319)  &  4 \ra 18  &  $\Lx$ \ra $\Lx$  \\
8  &  (4.7959, 3.0824, 0.0825) \ra (3.8358, 3.4537, 0.0859)  &  5 \ra 22  &  $\Lx$ \ra $\Lx$  \\
9  &  (0.2876, 2.4431, 0.0875) \ra (0.3410, 1.3964, 0.0047)  &  6 \ra 14  &  $\Lx$ \ra $\Lx$  \\
10 &  (5.0251, 3.1416, 0.0895) \ra (5.3410, 3.1728, 0.0491)  &  3 \ra 20  &  $\Lx$ \ra $\Lx$  \\
11 &  (2.6023, 3.1719, 0.1046) \ra (1.5060, 3.2037, 0.0762)  &  4 \ra 21  &  $\Lx$ \ra $\Lx$  \\
12 &  (2.6575, 3.1209, 0.1078) \ra (4.5024, 2.4768, 0.0754)  &  3 \ra 18  &  $\Lx$ \ra $\Lx$  \\
13 &  (6.0553, 0.0032, 0.1106) \ra (2.5186, 0.0000, 0.0803)  &  6 \ra 20  &  $\Lx$ \ra $\Lx$  \\
14 &  (2.6063, 3.1057, 0.1128) \ra (4.0182, 3.2164, 0.0948)  &  4 \ra 18  &  $\Lx$ \ra $\Lx$  \\
15 &  (2.2500, 3.1416, 0.1146) \ra (1.7332, 3.1416, 0.1020)  &  3 \ra 18  &  $\Lx$ \ra $\Lx$  \\
\hline
\end{tabular}
\caption{{Properties of solutions at initial and final points of continuation.}}
\label{tbl:properties_solutions}
\end{table}

Breaking or gaining of the additional symmetries \refEqn{additionalsym2} or \refEqn{additionalsym3}
was often detected, and gain of the additional symmetries \refEqn{additionalsym1} and \refEqn{additionalsym4}
was also uncovered. (More details appear in Tables~\refFigsSecs{tbl:properties_solutions_2}{tbl:properties_solutions_3}
and Section~\ref{sec:discussion_particular_solutions}.) We did not observe a change in stability of the solutions
at the points where additional symmetries were gained or broken, but the unstable dimension 
would usually change at said points (with an increase or decrease of $1$ or $2$). 
Table~\ref{tbl:properties_solutions_2} indicates which additional symmetries, if any, the invariant solutions posses, whether
continuation was done only on the parameter $R$ or not 
(as will be discussed in Section~\ref{sec:discussion_particular_solutions}), as well as the number $\Ni$ of 
distinct {CGLE} parameter points for which solutions were found in each sequence $\Ai{i}$, $i = 1,\ldots,15$
(see \refEqn{eqn:Asequence}).
To determine the number $\Ni$, we counted two points in the resulting numerical continuation path of the 
sequence $\Ai{i}$, say $\pts{p}{j}{(i)} = (\pts{R}{j}{(i)}, \pts{\nu}{j}{(i)}, \pts{\mu}{j}{(i)})$ and 
$\pts{p}{k}{(i)} = (\pts{R}{k}{(i)}, \pts{\nu}{k}{(i)}, \pts{\mu}{k}{(i)})$,  where $j \ne k$, as distinct if 
$|| \pts{p}{j}{(i)} - \pts{p}{k}{(i)} ||_{2} \ge 0.05$.
Approximate values of the {CGLE} parameters $(R,\nu,\mu)$ at which any additional symmetry was gained or broken during
the numerical continuation are listed in Table~\ref{tbl:properties_solutions_3}, only for those sequences where
symmetry gaining or breaking behavior occurred.
\begin{table}[!t]  
\tablesize
\centering
\begin{tabular}{|r|l|l|l|c|c|}
\hline
    &  \multicolumn{3}{c|}{additional symmetries}                   & continuation   &        \\    \cline{2-4}
id  &  start of continuation &  in between  &  end of continuation  & on $R$ only    & $\Ni$  \\  
\hline  \hline
1  &  none & none & none & yes   & 112 \\
2  &  none & \refEqn{additionalsym2} & \refEqn{additionalsym2} & yes & 188  \\ 
3  &  \refEqn{additionalsym1}, $\lsym=2$ & \refEqn{additionalsym1}, $\lsym=2$ & \refEqn{additionalsym1}, $\lsym=2$ & no & 179 \\
4  &  \refEqn{additionalsym2} & \refEqn{additionalsym2} & \refEqn{additionalsym2} & yes & 233 \\
5  &  none & none & none & yes  & 634 \\
6  &  none & none & none & yes  & 191  \\
7  &  none & \refEqn{additionalsym1}, $\lsym=2$,               & \refEqn{additionalsym1}, $\lsym=2$,               & no & 179 \\
   &       & \refEqn{additionalsym2}, \refEqn{additionalsym3}  & \refEqn{additionalsym2},  \refEqn{additionalsym3} &    &     \\
8  &  none & \refEqn{additionalsym2}  & none & yes & 489  \\
9  &  \refEqn{additionalsym1}, $\lsym=3$  &  \refEqn{additionalsym1}, $\lsym=3$ &  \refEqn{additionalsym1}, $\lsym=3$ & yes  & 116 \\
10 &  \refEqn{additionalsym2}  & \refEqn{additionalsym2} & none & no & 385  \\ 
11 &  none & \refEqn{additionalsym1}, $\lsym=2$,               & \refEqn{additionalsym1}, $\lsym=2$  & yes  & 415 \\
   &       & \refEqn{additionalsym2}, \refEqn{additionalsym3}  &                                     &      &     \\
12 &  none & none & none & no & 472 \\
13 &  none &  \refEqn{additionalsym4} &  \refEqn{additionalsym4} & no & 526  \\ 
14 &  none & \refEqn{additionalsym1}, $\lsym=2$,               & \refEqn{additionalsym1}, $\lsym=2$  & yes & 615 \\
   &       & \refEqn{additionalsym2}, \refEqn{additionalsym3}  &                                     &     &       \\
15 &  \refEqn{additionalsym1}, $\lsym=2$,  & \refEqn{additionalsym1}, $\lsym=2$,   &  \refEqn{additionalsym1}, $\lsym=2$,    & yes & 438  \\ 
   &  \refEqn{additionalsym2}, \refEqn{additionalsym3} & \refEqn{additionalsym2}, \refEqn{additionalsym3}  & \refEqn{additionalsym2}, \refEqn{additionalsym3}     &  &  \\ 
\hline
\end{tabular}
\caption{{Additional symmetries associated to $\cglegroup$-orbits along the continuation paths.}}
\label{tbl:properties_solutions_2}
\end{table}
\begin{table}[!t]   
\tablesize
\centering
\begin{tabular}{|r|l|}
\hline 
\vspace*{-0.7em}      &     \\
 id   & approximate $(R,\nu,\mu)$ values: type of symmetry gained/broken      \\
\hline  \hline
2  & $(16.5,-7,5)$: \refEqn{additionalsym2} gained  \ra
        $(20.6,-7,5)$: \refEqn{additionalsym2} broken  \\
     & \hspace*{1em}\ra $(53.4,-7,5)$: \refEqn{additionalsym2} gained      \\    
\vspace*{-0.7em}    &     \\  
4  &   $(20.6,-7,5)$: \refEqn{additionalsym2} broken  \ra
          $(54.7,-7,5)$: \refEqn{additionalsym2} gained      \\
\vspace*{-0.7em}    &     \\  
7  &  $(100,-7,0.3)$: \refEqn{additionalsym1}, $\lsym=2$, gained \ra
      $(100,-7,-0.02)$: \refEqn{additionalsym2}, \refEqn{additionalsym3} gained   \\
\vspace*{-0.7em}    &     \\  
8  &  $(32.6,-7,5)$: \refEqn{additionalsym2} gained \ra
      $(85.5,-7,5)$: \refEqn{additionalsym2} broken     \\
\vspace*{-0.7em}    &     \\  
10 &  $(82.8,-7,5.97)$: \refEqn{additionalsym2} broken \ra
          $(82.9,-7,5.97)$: \refEqn{additionalsym2} gained \\
     & \hspace*{1em}\ra  $(85.5,-7,5.97)$: \refEqn{additionalsym2} broken     \\ 
\vspace*{-0.7em}    &     \\  
11 &  $(12.17,-7,5)$: \refEqn{additionalsym1}, $\lsym=2$,
      \refEqn{additionalsym2},
      \refEqn{additionalsym3} gained \ra
      $(60.1,-7,5)$: \refEqn{additionalsym2},
      \refEqn{additionalsym3} broken    \\
\vspace*{-0.7em}    &     \\  
13 &  $(16,-5.6,3.4)$: \refEqn{additionalsym4} gained  \\ 
\vspace*{-0.7em}    &     \\  
14 &  $(16.2,-7,5)$: \refEqn{additionalsym2} gained \ra
          $(15.8,-7,5)$: \refEqn{additionalsym1}, $\lsym=2$,
          \refEqn{additionalsym3} gained  \\
     & \hspace*{1em}\ra  $(80.6,-7,5)$: \refEqn{additionalsym2},
         \refEqn{additionalsym3} broken   \ra
         $(80.8,-7,5)$: \refEqn{additionalsym2},
          \refEqn{additionalsym3} gained \\
     & \hspace*{1em}\ra  $(78.7,-7,5)$: \refEqn{additionalsym2},
          \refEqn{additionalsym3} broken   \\
\vspace*{-0.7em}    &     \\  
15 &  $(72.5,-7,5)$: \refEqn{additionalsym2}, \refEqn{additionalsym3} broken  \ra
          $(71.2,-7,5)$: \refEqn{additionalsym2}, \refEqn{additionalsym3} gained  \\
     & \hspace*{1em}\ra   $(84.8,-7,5)$: \refEqn{additionalsym2}, \refEqn{additionalsym3} broken    \ra
          $(67.7,-7,5)$: \refEqn{additionalsym2}, \refEqn{additionalsym3} gained  \\
     & \hspace*{1em}\ra  $(80.9,-7,5)$: \refEqn{additionalsym2}, \refEqn{additionalsym3} broken  \ra 
          $(84.9,-7,5)$: \refEqn{additionalsym2}, \refEqn{additionalsym3} gained    \\
\hline
\end{tabular}
\caption{{Summary of symmetries gained/broken along the continuation paths.} }
\label{tbl:properties_solutions_3}
\end{table}

\subsection{Features from the Solution Process}
\label{sec:discussion_particular_solutions}

Recall that we start the continuation from a point ($\cglegroup$-orbit)
in the fiber over the initial point $(R,\nu,\mu) = (16,-7,5)$ in the base
(space of parameters) tracing a path of $\cglegroup$-orbits (invariant solutions) 
which belong to fibers of the moduli space over the moving point in the base. 
(This was done 15 times starting from 15 different points ($\cglegroup$-orbits) in the moduli space belonging
to the fiber over the initial point $(R,\nu,\mu) = (16,-7,5)$ in the base.)  Given the nature of the continuation
method used and its implementation,
one may revisit a fiber over a particular point $(R,\nu,\mu)$ in the base several times during the continuation process. 
In other words, a continuation path in the moduli space of $\cglegroup$-orbits may turn around. 

To give an idea of the performance of the methodology employed,
Figures~\refFigsSecs{fig:sol8_continuation_paths}{fig:sol8_surfplot} show several plots corresponding to application
of the procedure for the sequence $\Ai{8}$ (see \refEqn{eqn:Asequence}), listed with id~8 in 
Tables~\refFigsSecs{tbl:properties_solutions}{tbl:properties_solutions_3}.  
Continuation in this case was done on the {CGLE} parameter $R$ only.
Paths resulting from the continuation appear in Figure~\ref{fig:sol8_continuation_paths}.
Specifically, Figures~\ref{fig:sol8_RvsT}, \ref{fig:sol8_RvsS}, and \ref{fig:sol8_Rvsphi} depict the resulting continuation paths 
by displaying, respectively, the values of the time period $T$, space translation $S$, 
and rotation\footnote{Since it was not strictly necessary, the constraint $\varphi \in [0, 2\pi)$ was not 
explicitly enforced when solving the system of nonlinear algebraic equations \refEqn{eqn:cgle_nleqns}.
Furthermore, after performing a series of preliminary test runs, we found no advantage
(from a computational point of view) in enforcing it.
The values of $\varphi$ in Figure~\ref{fig:sol8_Rvsphi}
are displayed as they resulted from the solution of the system \refEqn{eqn:cgle_nleqns}, and should be 
taken modulo an integer multiple of $2\pi$, mapping them back to the interval $[0,2\pi)$.}
$\varphi$ as functions of the continuation parameter $R$.
Note from Figure~\ref{fig:sol8_RvsS} that within the range of $R \approx 32.6$ through $R \approx 85.5$
the value of $S$ remained constant.  The start of this interval of constant $S$ corresponds to a step in the continuation 
process at which the resulting solution gained the additional symmetry \refEqn{additionalsym2}; this symmetry was broken 
at the point in the path where $S$ ceases to be constant.
(Recall that solutions with symmetries \refEqn{eqn:cgle_invariant_solution} and \refEqn{additionalsym2} exist in subspaces
of the solution space $(A; \varphi,S,T)$ for which either $S=0$ or $S=\Lx/2$; see Section~\ref{sec:problem_statement}.)
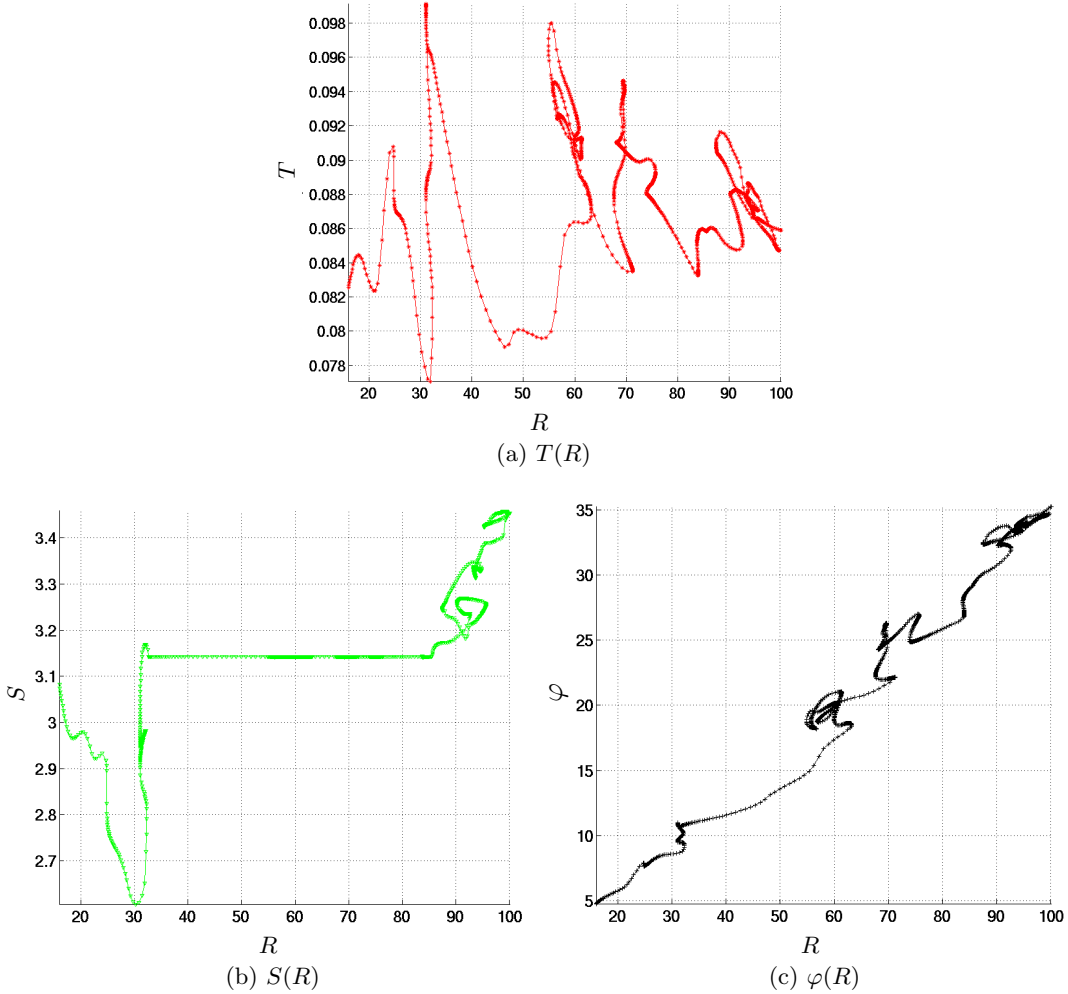
\begin{figure}[!t]  
        \centering
        \subfloat[$T(R)$]{
                {\includegraphics[trim = 12mm 13mm 18mm 12mm, clip, width=0.43\textwidth]{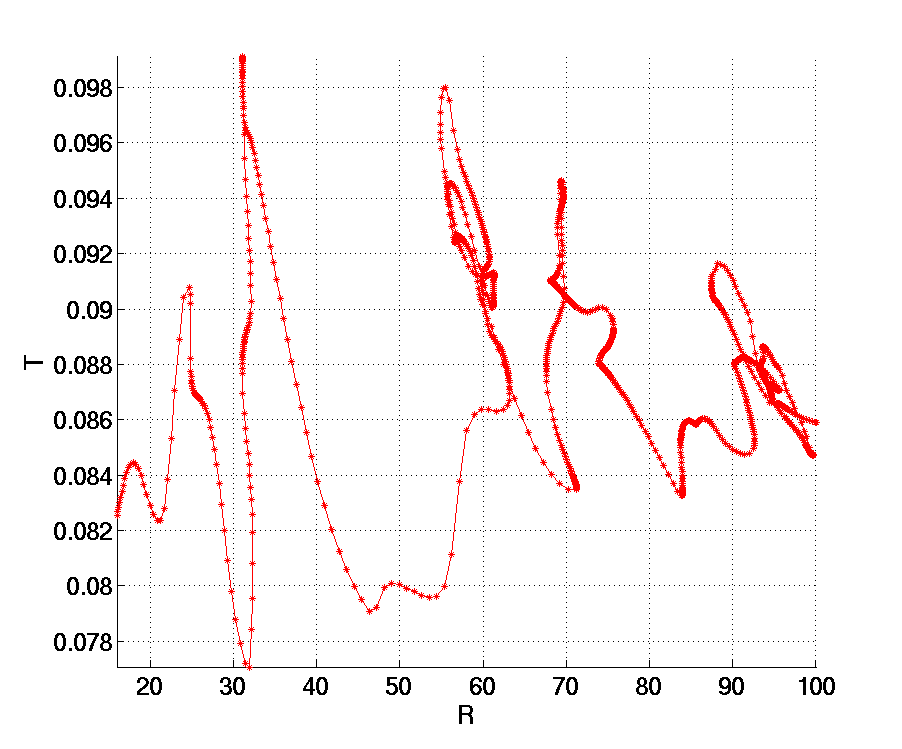}}
                \put(-195,85){\rotatebox{90}{\footnotesize $T$}}
                \put(-100,-10){\footnotesize $R$}
                \label{fig:sol8_RvsT}}
        \\
        \subfloat[$S(R)$]{
                {\includegraphics[trim = 21mm 13mm 18mm 12mm, clip, width=0.43\textwidth]{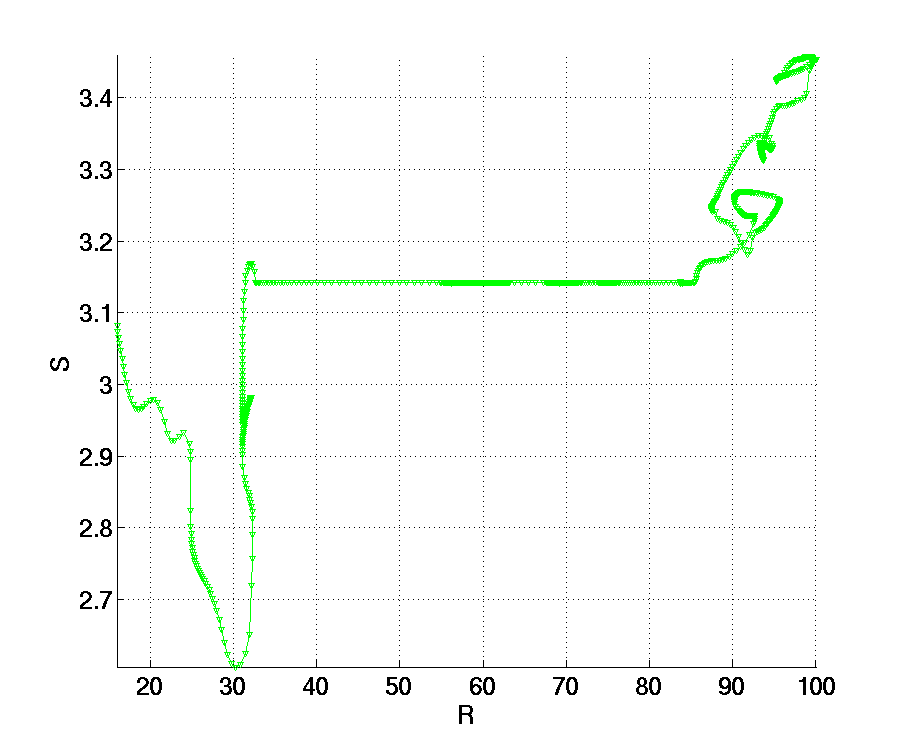}}
                \put(-195,85){\rotatebox{90}{\footnotesize $S$}}
                \put(-100,-10){\footnotesize $R$}
                \label{fig:sol8_RvsS}}
        \hspace*{2.5ex}
        \subfloat[$\varphi(R)$]{
                {\includegraphics[trim = 23mm 13mm 18mm 12mm, clip, width=0.43\textwidth]{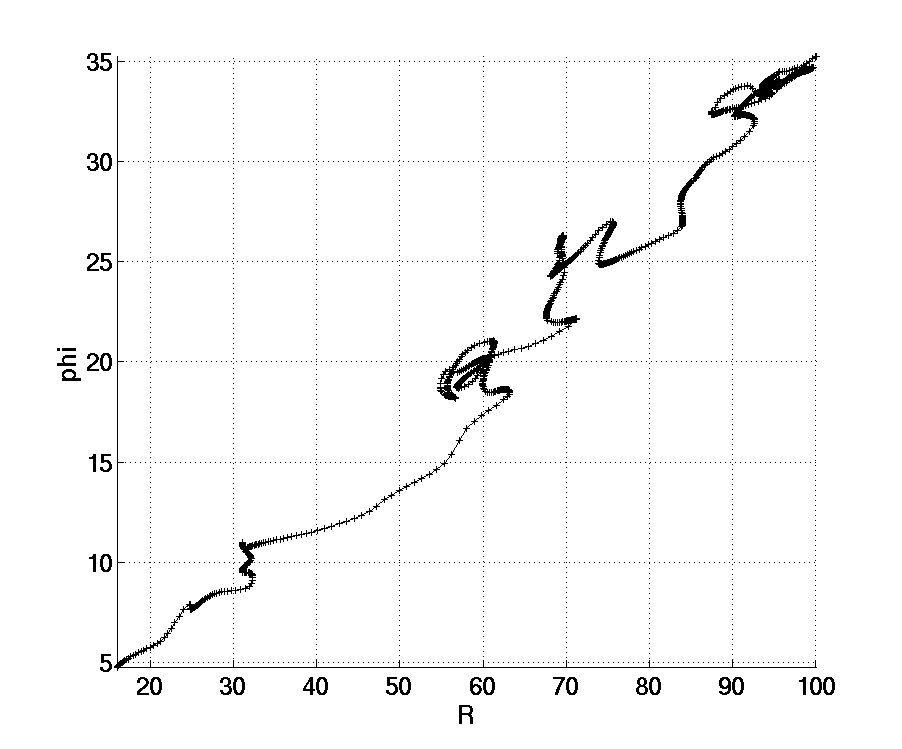}}
                \put(-195,85){\rotatebox{90}{$\varphi$}}
                \put(-100,-10){\footnotesize $R$}
                \label{fig:sol8_Rvsphi}}
        \caption{{Continuation sequence $\Ai{8}$:
                       Paths traversed by $T$, $S$, and $\varphi$, as functions of $R$.
                       }}
        \label{fig:sol8_continuation_paths}
\end{figure}

The depictions in Figure~\ref{fig:sol8_continuation_paths} make it convenient to identify turning points in the continuation 
path and, for a given (fixed) value of the continuation parameter, whether there may exist multiple solutions of $\Feqzero$ in 
the path.  For example, in Figure~\ref{fig:sol8_RvsT} one can identify four points where the line $R = 90$ intersects
the curve $T(R)$.  These four points correspond to four invariant solutions computed at the same 
particular point $(R, \nu, \mu)$ in parameter space. Then, the multitude of solutions associated with this point in parameter space
can be inspected to determine whether they are different elements of the same orbit \refEqn{eqn:torus3orbits} 
of the symmetry group $\cglegroup$, whether they belong to conjugate orbits of the symmetry group, or whether they 
belong to different (non-conjugate) orbits of the symmetry group.

\begin{figure}[!t]  
        \centering
        \subfloat[spatial spectrum]{
                {\includegraphics[trim = 20mm 10mm 20mm 12mm, clip, width=0.46\textwidth]{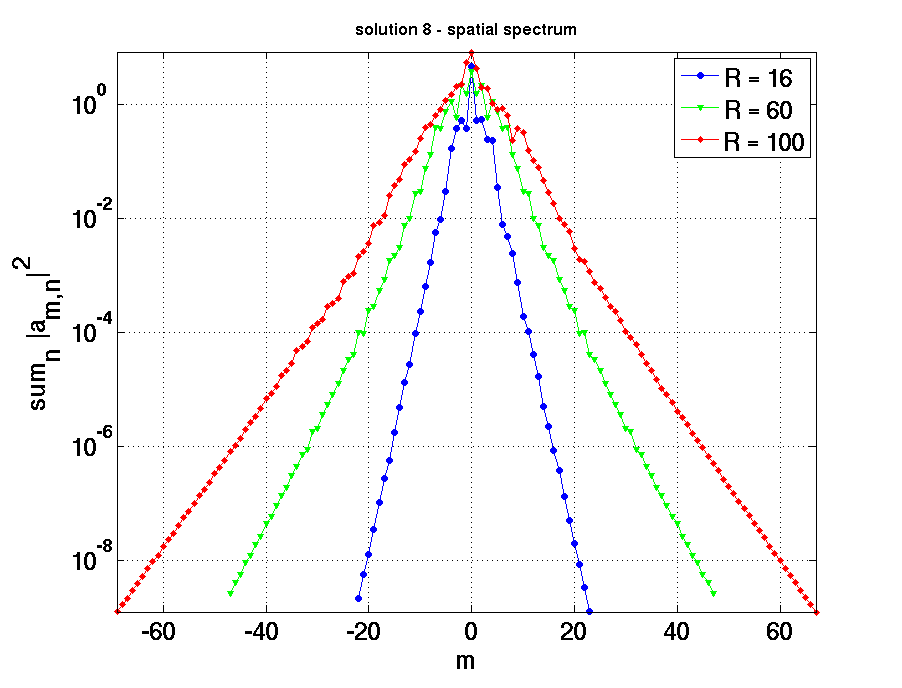}}
                \put(-210,55){\rotatebox{90}{\footnotesize $\sum_n |\amn{m}{n}|^2$}}
                \put(-100,-7){\small $m$}
                \label{fig:sspec}}
        \hspace*{3.0ex}
        \subfloat[temporal spectrum]{
                {\includegraphics[trim = 20mm 10mm 20mm 12mm, clip, width=0.46\textwidth]{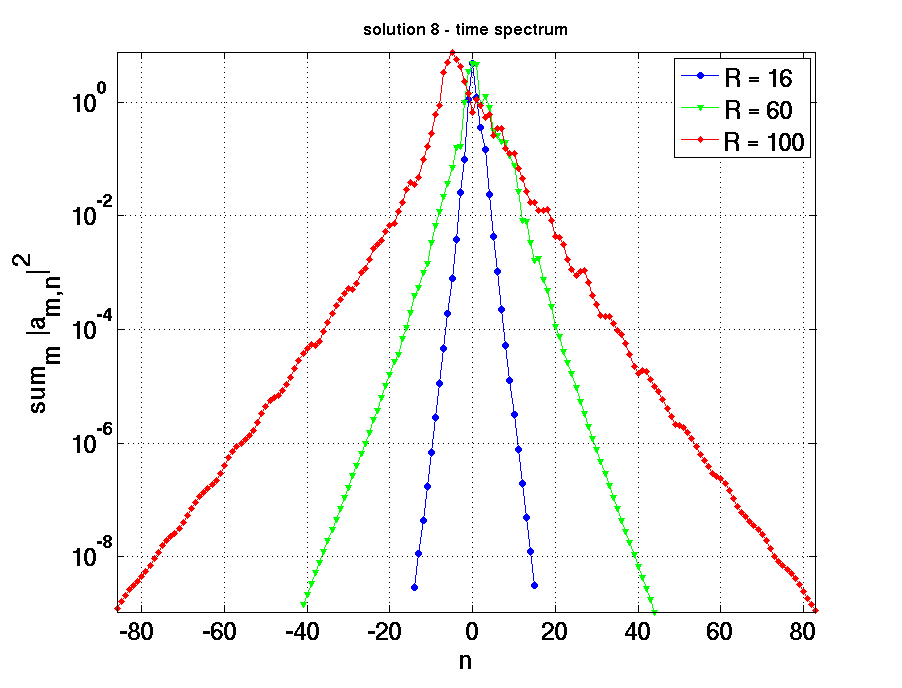}}
                \put(-210,55){\rotatebox{90}{\footnotesize $\sum_m |\amn{m}{n}|^2$}}
                \put(-100,-7){\small $n$}
                \label{fig:tspec}}
        \caption{{Spectra for three solutions in the sequence $\Ai{8}$,
                       at values of $R=16, 60, 100$ ($\nu =-7$, $\mu=5$).}}
        \label{fig:sol8_spectra}
\end{figure}

\begin{figure}[!t] 
        \centering
        \subfloat[$\Re(A)$, $R=16$]{
                {\includegraphics[trim = 7mm 5mm 10mm 10mm, clip, width=0.29\textwidth]{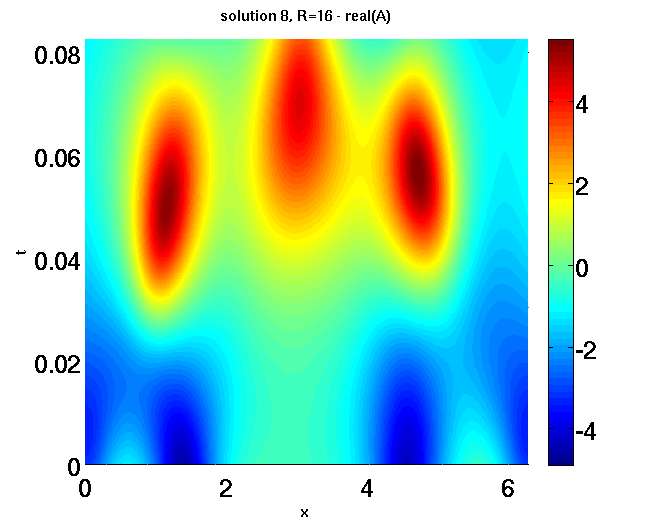}}
                \put(-65,-4){\footnotesize{\textit{x}}}
                \put(-130,52){\rotatebox{90}{\footnotesize{{\textit{t}}}}}
                \label{fig:R16_real}}
        \hspace*{0.5ex}
        \subfloat[$\Im(A)$, $R=16$]{
                {\includegraphics[trim = 7mm 5mm 10mm 10mm, clip, width=0.29\textwidth]{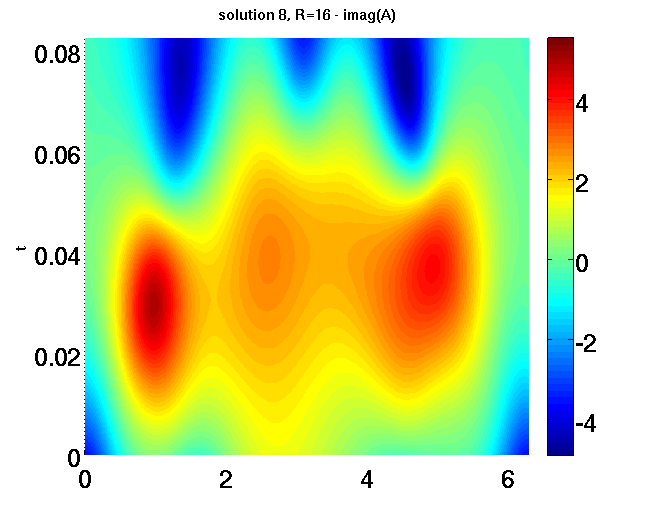}}
                \put(-65,-4){\footnotesize{\textit{x}}}
                \put(-130,52){\rotatebox{90}{\footnotesize{{\textit{t}}}}}
                \label{fig:R16_imag}}
        \hspace*{0.5ex}
        \subfloat[$|A|$, $R=16$]{
                {\includegraphics[trim = 7mm 5mm 10mm 10mm, clip, width=0.29\textwidth]{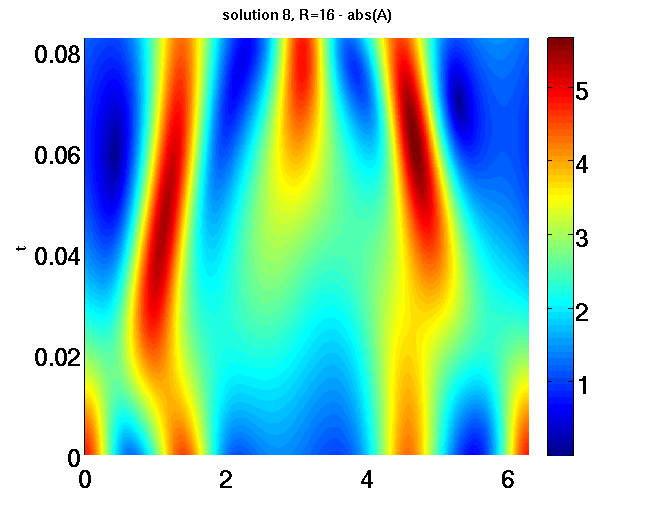}}
                \put(-65,-4){\footnotesize{\textit{x}}}
                \put(-130,52){\rotatebox{90}{\footnotesize{{\textit{t}}}}}
                \label{fig:R16_abs}}

        \vspace*{0.02ex}
        \subfloat[$\Re(A)$, $R=60$]{
                {\includegraphics[trim = 7mm 5mm 10mm 10mm, clip, width=0.29\textwidth]{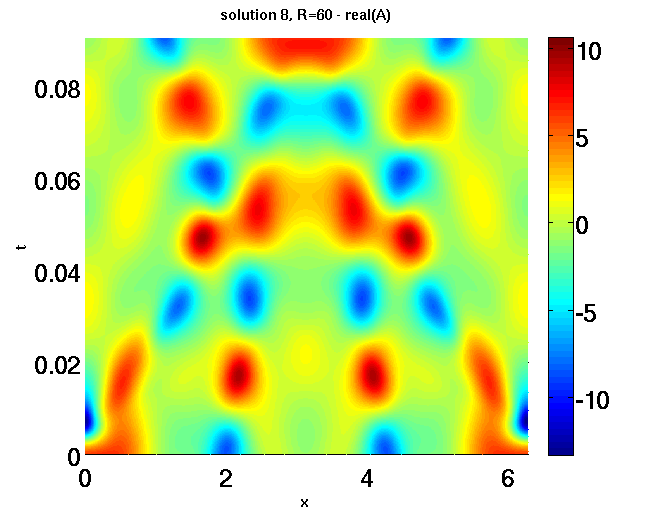}}
                \put(-65,-4){\footnotesize{\textit{x}}}
                \put(-130,52){\rotatebox{90}{\footnotesize{{\textit{t}}}}}
                \label{fig:R60_real}}
        \hspace*{0.5ex}
        \subfloat[$\Im(A)$, $R=60$]{
                {\includegraphics[trim = 7mm 5mm 10mm 10mm, clip, width=0.29\textwidth]{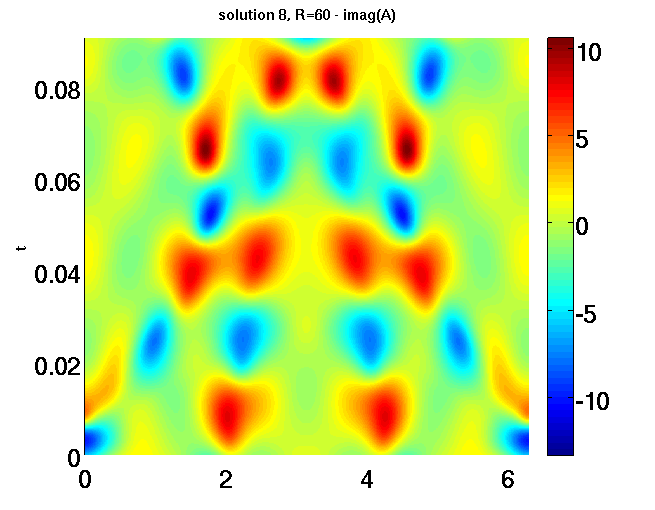}}
                \put(-65,-4){\footnotesize{\textit{x}}}
                \put(-130,52){\rotatebox{90}{\footnotesize{{\textit{t}}}}}
                \label{fig:R60_imag}}
        \hspace*{0.5ex}
        \subfloat[$|A|$, $R=60$]{
                {\includegraphics[trim = 7mm 5mm 10mm 10mm, clip, width=0.29\textwidth]{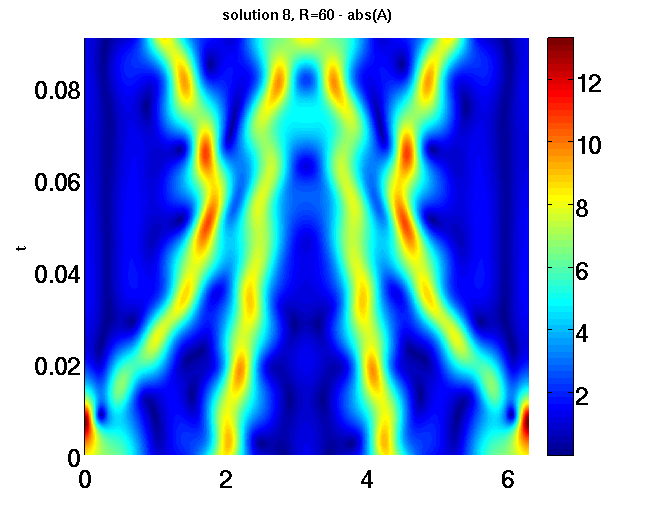}}
                \put(-65,-4){\footnotesize{\textit{x}}}
                \put(-130,52){\rotatebox{90}{\footnotesize{{\textit{t}}}}}
                \label{fig:R60_abs}}

        \vspace*{0.02ex}
        \subfloat[$\Re(A)$, $R=100$]{
                {\includegraphics[trim = 7mm 5mm 10mm 10mm, clip, width=0.29\textwidth]{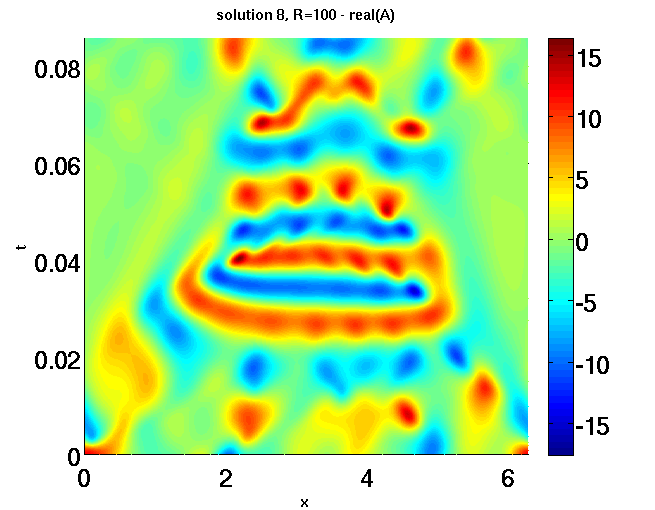}}
                \put(-65,-4){\footnotesize{\textit{x}}}
                \put(-130,52){\rotatebox{90}{\footnotesize{{\textit{t}}}}}
                \label{fig:R100_real}}
        \hspace*{0.5ex}
        \subfloat[$\Im(A)$, $R=100$]{
                {\includegraphics[trim = 7mm 5mm 10mm 10mm, clip, width=0.29\textwidth]{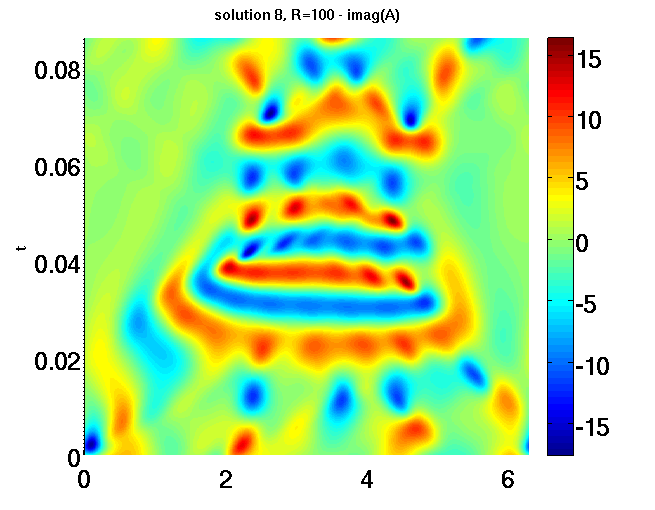}}
                \put(-65,-4){\footnotesize{\textit{x}}}
                \put(-130,52){\rotatebox{90}{\footnotesize{{\textit{t}}}}}
                \label{fig:R100_imag}}
        \hspace*{0.5ex}
        \subfloat[$|A|$, $R=100$]{
                {\includegraphics[trim = 7mm 5mm 10mm 10mm, clip, width=0.29\textwidth]{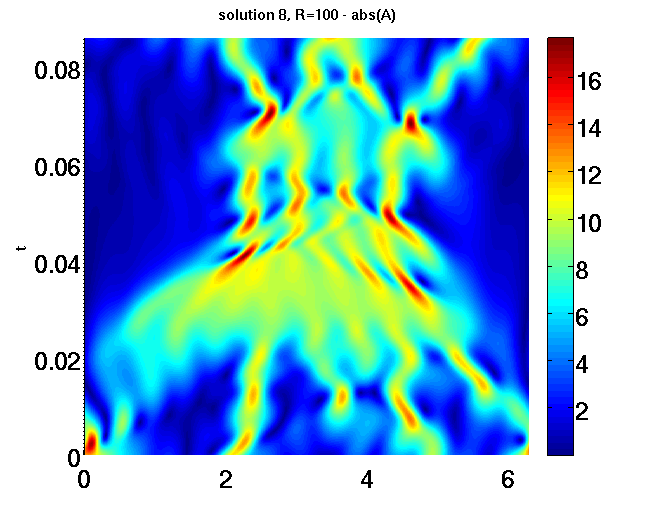}}
                \put(-65,-4){\footnotesize{\textit{x}}}
                \put(-130,52){\rotatebox{90}{\footnotesize{{\textit{t}}}}}
                \label{fig:R100_abs}}       
        \caption{{Surface plots of the real part $\Re(A)$, imaginary part $\Im(A)$, and absolute value $|A|$
                 for three solutions in sequence $\Ai{8}$,
                 at values of $\nu=-7, \mu=5$ and $R=16$ (top), $R=60$ (middle), and $R=100$ (bottom).
                 Symmetry \refEqn{additionalsym2} was gained at $(R,\nu,\mu) \approx (32.6,-7,5)$ 
                 and broken at $(R,\nu,\mu) \approx (85.5,-7,5)$.  
                 Hence $A$ is even (about $x = \pi$) for $R=60$, but not for $R=16, 100$.
                 A movie depicting the continuation path in $(\varphi,S,T)$
                 space, as well as solutions represented by plotting $\Im(A(x,0))$ vs.~$\Re(A(x,0))$ 
                 at a sequence of continuation steps, can be found among the supplementary material associated with this article. }}
        \label{fig:sol8_surfplot} 
\end{figure}

Spectra for several solutions in the path from $R=16$ to $R=100$ are shown in Figure~\ref{fig:sol8_spectra}.
As expected, an increase in the value of $R$ requires more terms in the expansions \refEqn{eqn:xFseries} and 
\refEqn{eqn:am_ansatz} in order to keep a suitable decay in both the spatial and temporal spectra for the solutions.
Finally, surface plots of the real part $\Re(A)$, imaginary part $\Im(A)$, and absolute value $|A|$ for the solutions
whose spectra are shown in Figure~\ref{fig:sol8_spectra} appear in Figure~\ref{fig:sol8_surfplot}, where the aforementioned 
gain and, thereafter, loss of symmetry \refEqn{additionalsym2} can be observed.

The example above illustrates the general situation that one faces. By this we mean that, due to the use of the arc-length
continuation option from the {LOCA} package \cite{loca}, which was the appropriate choice for us because it allows for 
turning points in the path following process, it is possible (i.e., inherent in the continuation algorithm) that a point
$\pt{p}{l} = (\pt{R}{l}, \pt{\nu}{l}, \pt{\mu}{l})$ in the {CGLE} parameter space may return to itself, that is, 
$\pt{p}{k+l} = \pt{p}{l}$, after $k$ continuation steps.  In such a situation, we must consider different cases
for the solutions
$( \Ap{\pt{p}{k+l}}  ;  \varphi(\pt{p}{k+l}), S(\pt{p}{k+l}), T(\pt{p}{k+l}) )$ and
$( \Ap{\pt{p}{l}}  ;  \varphi(\pt{p}{l}), S(\pt{p}{l}), T(\pt{p}{l}) )$.  
Namely, whether said solutions represent
(i) the same $\cglegroup$-orbit,
that is, $(\varphi(\pt{p}{k+l}), S(\pt{p}{k+l}), T(\pt{p}{k+l})) =  (\varphi(\pt{p}{l}), S(\pt{p}{l}), T(\pt{p}{l}))$
and there exists some $(\Arot,\spacetrans,\timetrans) \in \cglegroup$ such that 
$\Ap{\pt{p}{k+l}} = (\Arot,\spacetrans,\timetrans) \cdot \Ap{\pt{p}{l}}$;
(ii) different, but conjugate, $\cglegroup$-orbits, as defined in Section~\ref{sec:problem_statement}
(see also \refEqns{Fsym5a}{Fsym5b}); or
(iii) different, non-conjugate, $\cglegroup$-orbits.
Along a continuation path, say $\Ai{i}$, many returns to a same point $\pts{p}{l}{(i)}$ do occur.
However, we include only one of the invariant solutions computed at $\pts{p}{l}{(i)}$ in the
count $\Ni$ in Table~\ref{tbl:properties_solutions_2}, since the presentation of the complete analysis of the multitude 
of invariant solutions that were computed at such ``revisited'' points is out of the scope of this paper.

Challenging behavior that arose during the numerical continuation was often due to traversal of values of the 
continuation parameter in a cyclic manner, specifically related to the cases (i) and (ii) listed in the previous paragraph.  
As a result, the continuation path within these cycles would contain  (different) elements in the same $\cglegroup$-orbit, 
or solutions representing conjugate $\cglegroup$-orbits (cf.~\refEqn{eqn:Smap} and \refEqns{Fsym5a}{Fsym5b}).
Within the cycles, solutions at the turning points in the continuation path were in the vicinity of solutions with additional 
symmetries, or near solutions with a smaller spatial period $\Lx / q_{\sss 1}$ for some integer $q_{\sss 1} > 1$ or smaller 
time period $T/q_{\sss 2}$ for some integer $q_{\sss 2} > 1$.  
Often the {LOCA} continuation algorithm \cite{loca} would exit from the cycles automatically, so that the
procedure would again start yielding solutions in different, non-conjugate, $\cglegroup$-orbits, as well as continue to 
make progress towards the goal of reaching the (desired) final point in the {CGLE} parameter space.  
However, sometimes the continuation algorithm would get caught in said cycles.
We discuss instances of these scenarios in the following paragraphs.

\begin{figure}[!t]  
        \centering
        \subfloat[$R$, $T$, $S$, $\varphi$ vs.~continuation step number]{
                {\includegraphics[trim = 15mm 10mm 20mm 13mm, clip, width=0.46\textwidth]{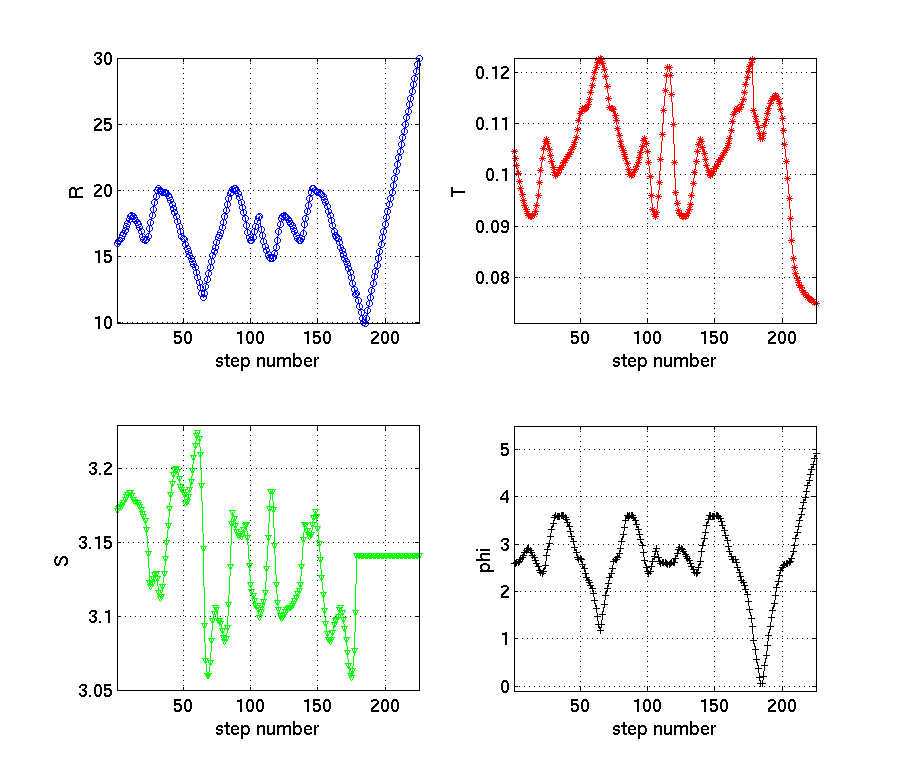}}
                \label{fig:sol11_RTSphi}}
        \subfloat[$S(R)$]{
                {\includegraphics[trim = 10mm 0mm 18mm 12mm, clip, width=0.46\textwidth]{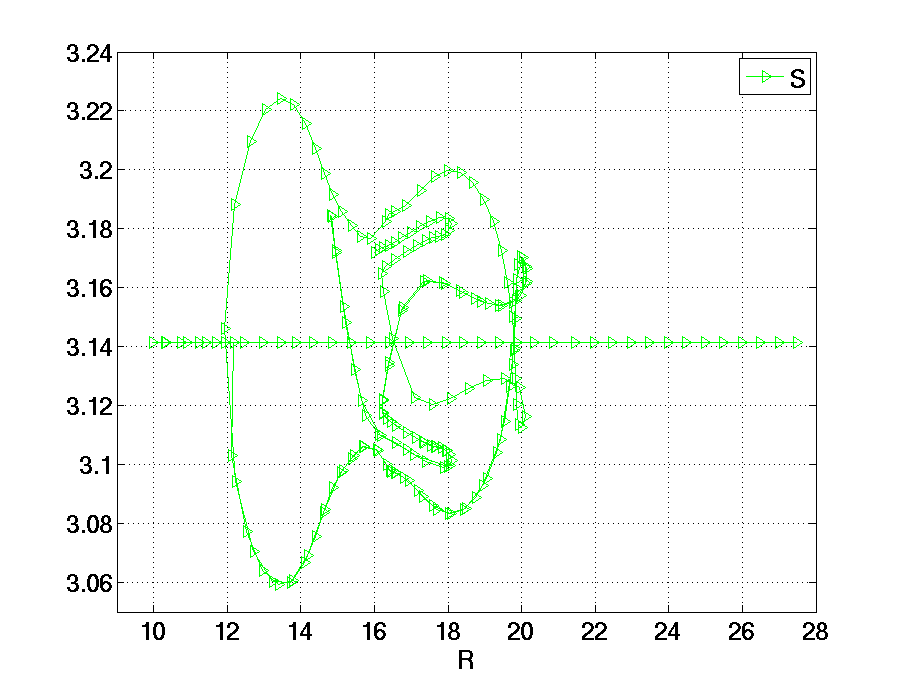}}
                \label{fig:sol11_RvsS}}
\\
        \vspace*{0.1ex}
        \subfloat[$\Re({\amn{0}{0}})$ vs.~$R$; \ $\Im({\amn{0}{0}})$ vs.~$R$]{
                {\includegraphics[trim = 15mm 0mm 15mm 10mm, clip, width=0.46\textwidth]{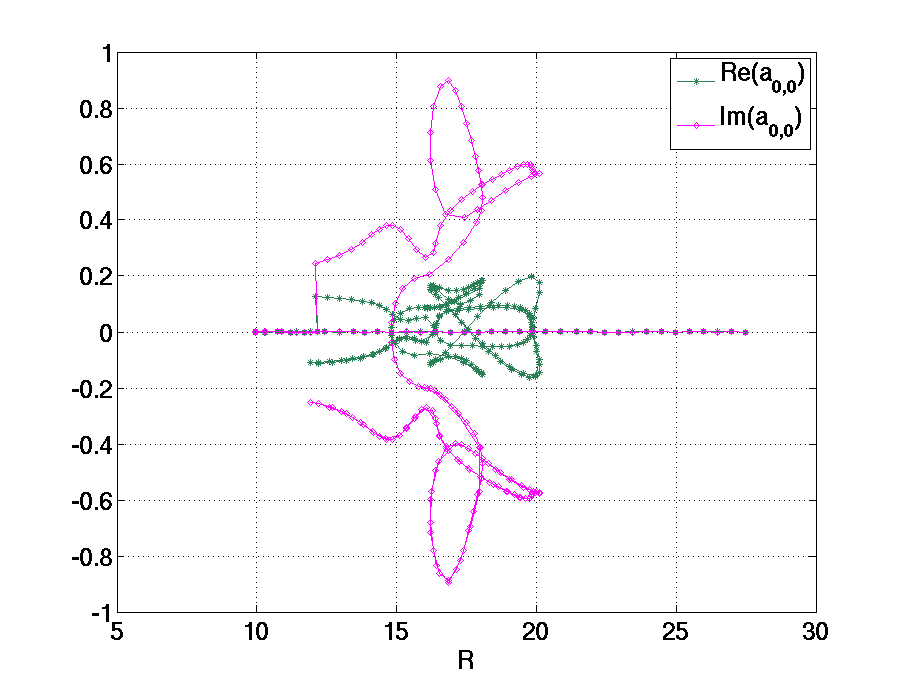}}
                \label{fig:sol11_Rvsa00}}
        \subfloat[$\Re({\amn{-1}{0}})$ vs.~$R$; \ $\Im({\amn{-1}{0}})$ vs.~$R$]{
                {\includegraphics[trim = 15mm 0mm 15mm 10mm, clip, width=0.46\textwidth]{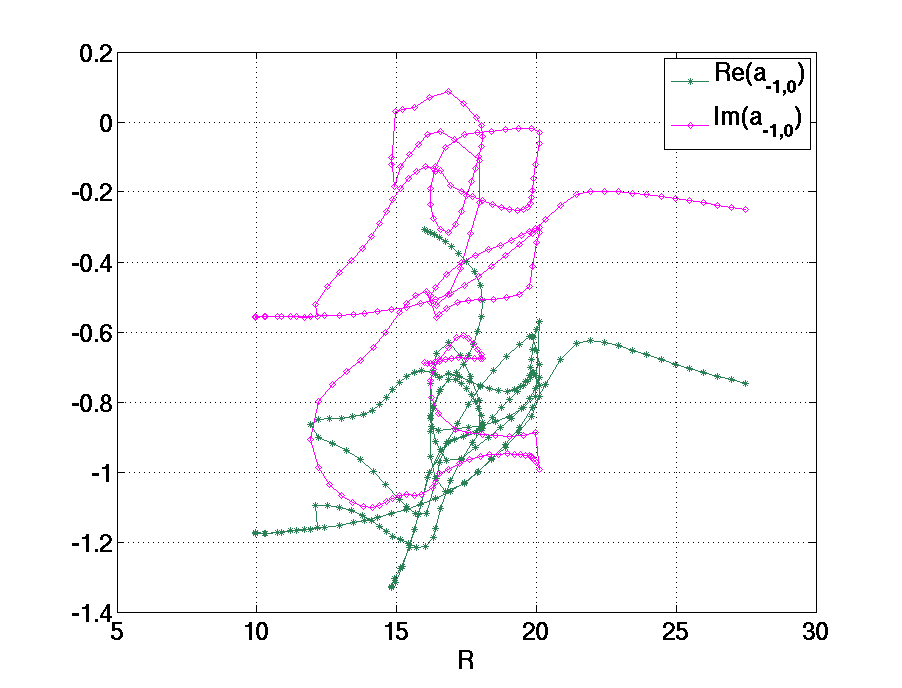}}
                \label{fig:sol11_Rvsam10}}
        \caption{{Cycling behavior during the numerical continuation for sequence $\Ai{11}$:
                       (a) Values of the continuation parameter $R$ and computed generator $(\varphi,S,T)$ as functions
                       of the continuation step number.
                       (b) $S$ as a function of $R$.
                       (c) $\Re({\amn{0}{0}})$ and $\Im({\amn{0}{0}})$ as functions of $R$.
                       (d) $\Re({\amn{-1}{0}})$ and $\Im({\amn{-1}{0}})$ as functions of $R$.
                       A movie depicting the continuation path 
                       in $(\varphi,S,T)$ space, as well as solutions represented by plotting 
                       $\Im(A(x,0))$ vs.~$\Re(A(x,0))$ at a sequence of continuation steps, can be found among the supplementary material associated with this article. 
                       The cycle in sub-figure~(b) above can be seen in the movie during the initial steps of the
                       path in $(\varphi,S,T)$ space. 
                       Symmetries  \refEqns{additionalsym1}{additionalsym3} were gained at $(R,\nu,\mu) \approx (12.17,-7,5)$;
                       symmetries  \refEqns{additionalsym2}{additionalsym3} were broken at $(R,\nu,\mu) \approx (60.1,-7,5)$.
                       }}
        \label{fig:sol11_cycles}
\end{figure}

An example of such cyclic behavior is depicted in Figure~\ref{fig:sol11_cycles}
for the sequence $\Ai{11}$, listed with id~11 in Tables~\refFigsSecs{tbl:properties_solutions}{tbl:properties_solutions_3}.
Continuation was done on the parameter $R$ only. 
Figure~\ref{fig:sol11_RTSphi} displays the values of the continuation parameter $R$, time period $T$, 
space translation $S$, and rotation $\varphi$ as functions of the continuation step number. 
Traversal of repeated values for $R$, $T$, $S$, and $\varphi$ is observed from the sub-figures in 
Figure~\ref{fig:sol11_RTSphi}, where it is also seen that the 
cycling behavior stops when $R \approx 12.17$, at around continuation step number 180 (where $S$ becomes constant), 
at which point the additional symmetries \refEqns{additionalsym1}{additionalsym3} are gained. 

Looking at Figure~\ref{fig:sol11_RvsS}, where the value of the space translation $S$ is plotted 
as a function of the continuation parameter $R$, one can see the cyclic behavior of $R$
resulting in a symmetric curve with respect to the horizontal line at the vertical axis value of $\Lx/2 = \pi$.  
The path depicted in Figure~\ref{fig:sol11_RvsS}, represented by the curve $S(R)$, contains conjugate solutions
(cf.~\refEqn{eqn:Smap} and \refEqns{Fsym5a}{Fsym5b}).
More precisely, the numerical continuation path for values of $R \in [12,20]$ that contains 
conjugate solutions is the one that yields the symmetric curve about the horizontal line at the
value of $\Lx/2 = \pi$ (seen in Figure~\ref{fig:sol11_RvsS}).  That is, points on the curve $S(R)$ that are
mirror images with respect to the line $\Lx/2 = \pi$ correspond to conjugate solutions under the spatial reflection
symmetry of the {CGLE}, which belong to conjugate orbits of the symmetry group $\cglegroup$.
The additional symmetry \refEqn{additionalsym2} 
was gained at a value of $R \approx 12.17$, and at this point the cycling behavior stops and the spatial 
translation $S$ takes on the value of $\Lx/2$, as solutions with symmetries \refEqn{eqn:cgle_invariant_solution} and
\refEqn{additionalsym2} exist in subspaces of the solution space $(A; \varphi,S,T)$ for which either $S=0$ or 
$S = \Lx/2$ (refer to Section~\ref{sec:problem_statement}).

Also, the additional symmetry \refEqn{additionalsym1}, for $\lsym = 2$, was gained along with the additional 
symmetry \refEqn{additionalsym2}.  Recall from \refEqn{eqn:Fcoefs_lsym} that the Fourier coefficients 
$\{\amn{m}{n}\}$ of solutions with symmetry \refEqn{additionalsym1}, for $\lsym = 2$, satisfy $\amn{m}{n} = 0$ 
if $m$ is even.  Thus, we can visualize gain of this additional symmetry by selecting a coefficient $\amn{m}{n}$, 
for some even $m$ and some $n$, and plotting its value as a function of the continuation parameter $R$, as done in
Figure~\ref{fig:sol11_Rvsa00} for the coefficient $\amn{0}{0}$.  
At the point when this additional symmetry is gained, 
for $R \approx 12.17$, we see that the real part of the coefficient $\amn{0}{0}$ goes from (around) 0.22 to 0, whereas 
its imaginary part goes from (around) 0.1 to 0.  Upon gaining the additional symmetry \refEqn{additionalsym1}, for $\lsym = 2$, 
the coefficient $\amn{0}{0}$ remains equal to zero, as seen in the path depicted in Figure~\ref{fig:sol11_Rvsa00}.
Finally, from \refEqn{eqn:Fcoefs_lsym}, a solution with symmetry \refEqn{additionalsym1}, for $\lsym = 2$, will have
nonzero Fourier coefficients $\{\amn{m}{n}\}$ for odd $m$.  This is depicted
for the coefficient $\amn{-1}{0}$, plotted as a function of the continuation parameter $R$, in Figure~\ref{fig:sol11_Rvsam10}.
As seen, the coefficient $\amn{-1}{0}$ remains nonzero after the additional symmetry \refEqn{additionalsym1}
is gained (at the same time when the cycling behavior stops) at a value of $R \approx 12.17$.

The aforementioned traversal of values of the continuation parameter in a cyclic manner was quite common, 
and often the {LOCA} continuation algorithm \cite{loca} would exit from the cycles automatically, that is, 
without us having to stop and restart the continuation with different values for the allowed increments on 
the continuation parameter.  Nevertheless, as an alternative for circumventing such cycling behavior,
we also experimented with taking the other parameters $\mu$ or $\nu$ in the {CGLE} \refEqn{eqn:cgle_pde} 
as continuation parameters. The continuation was always done on a single parameter at a time, while still 
all solutions were numerically continued from the regime with parameter values $(R, \nu, \mu) = (16, -7, 5)$ to the 
regime for $(R, \nu, \mu) = (100, -7, 5)$.  As a starting point for performing continuation on an alternate parameter, 
we would select a solution within the cycle for which the spectra (spatial or temporal, as appropriate) did not display
characteristics typical of that of solutions around the turning points in the cycle.   
As an example, with the solutions represented via \refEqn{eqn:am_ansatz}, given an integer $q_1 > 1$, the Fourier 
coefficients $\{\amn{m}{n}\}$ of a solution with spatial period $\Lx/q_1$  have a recognizable pattern
of zeros, namely, $\amn{m}{n} \!=\! 0$ if $m$ is not divisible by $q_1$.  So if the numerical continuation was caught
in a cycle where solutions around a turning point were close to a solution with spatial period $\Lx/q_1$,
as a starting point for performing continuation on an alternate parameter we could select a solution within 
the cycle for which $\sum_n | \amn{m}{n}|^2 > \varepsilon$, for $m = 0,\pm 1, \ldots, \pm q_1$, and some 
cutoff, say, $\varepsilon = 10^{-2}$.

\begin{figure}[!t]    
    \centering
    {\includegraphics[trim = 13mm 7mm 15mm 10mm, clip, width=0.55\textwidth]{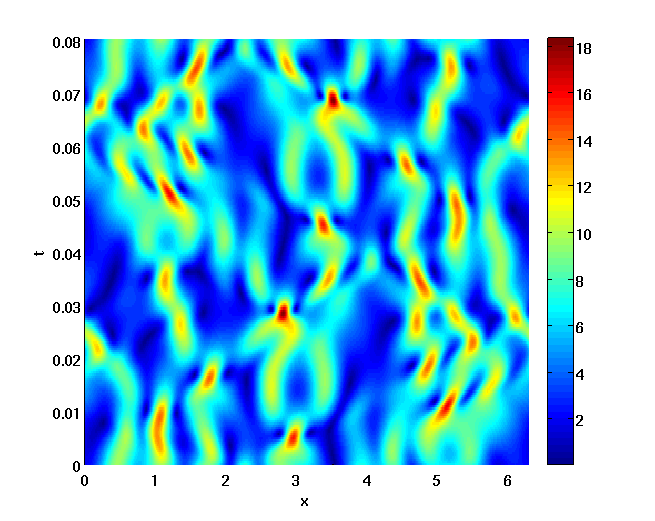}}
    \put(-122,-4){\textit{{x}}}
    \put(-250,110){\rotatebox{90}{\textit{{t}}}}
    \caption{{Surface plot of the absolute value $|A|$ for a solution in sequence $\Ai{13}$,
                   at parameter values $(R,\nu,\mu)=(100,-7,5)$.  The plot is over the 
                   space period $[0,\Lx]$ on the horizontal axis and time period $[0,T]$ on the vertical axis.
                   The solution has the additional symmetry \refEqn{additionalsym4},
                   with $c = \Lx, \tilde{\varphi} = \varphi/2 + \pi, \tilde{T} = T/2$, so $|A(x,t)| = |A(\Lx -x, t + T/2)|$,
                   as observed from the plot.
                   A movie depicting the continuation 
                   path in $(\varphi,S,T)$ space, 
                   as well as surface plots of $|A|$ at a sequence of continuation steps, can be found among the supplementary material associated with this article.
                   The additional symmetry \refEqn{additionalsym4} was gained at $(R,\nu,\mu) \approx (16,-5.6,3.4)$. }}
    \label{fig:sol13_abs}
\end{figure}

One particular case in which it was beneficial to alternate the continuation parameter was for the sequence $\Ai{13}$,
for which exiting automatically from cycling behavior in the vicinity 
of a solution that was even and had space period $\Lx/3$ and time period of $T/2$ was challenging. 
Hence we experimented with alternating the continuation parameter, as indicated in the previous paragraph. 
The continuation then led to a solution with the additional symmetry \refEqn{additionalsym4},
along with the invariance \refEqn{eqn:cgle_invariant_solution}.  This additional symmetry was gained
at parameter values of  $R = 16, \nu \approx -5.6, \mu \approx 3.4$.
Patterns resulting from the additional symmetry \refEqn{additionalsym4} can be visualized from the surface plot
shown in Figure~\ref{fig:sol13_abs}.

\subsection{Comments on Numerical Aspects}
\label{sec:comments_numerics}

Values of $N_x \in [48, 128]$ and  $N_t \in [48, 128]$ were used, respectively, in the truncation 
of the spatial Fourier series expansion \refEqn{eqn:xFseries} and the representation \refEqn{eqn:am_ansatz}.  
(For comparison, values of $N_x=32$ and $N_t = 48, 64$ were used in the preceding study \cite{lopez05}.) 
The number of terms $N_x$, $N_t$ in each expansion was chosen so that the 
solutions had a decay of at least $10^{-6}$ in their spatial and temporal spectra.  
The resulting number of unknowns for the system $\Feqzero$ of 
nonlinear algebraic equations \refEqn{eqn:cgle_nleqns} ranged between 4,000 and 32,260.

To solve the linear systems $\M_s\tilde{\x}_{{\sss k}} = \tilde{b}_{{\sss k}}$ using the {GMRES} iterative solver from the 
Meschach library \cite{meschach},
we set a tolerance of $10^{-9}$ for the residual $||\M_s\tilde{\x}_{{\sss k}} - \tilde{b}_{{\sss k}}||_2$
and a maximum of 3,000 {GMRES} iterations.  Recall that the solution of said linear systems is needed for the
computation of the Newton step, as described in Section~\ref{sec:continuation_procedure} and the associated
Appendix~\ref{app:newton_step}.
The number of iterations taken by the {GMRES} solver to meet the specified residual tolerance
ranged between 90 and 2,700.  In terms of actual computing time (on a ThinkPad W530 personal workstation with 
2.70 {GHz} processor speed), this translated to fractions of a second on the lower end to around 45-60 seconds on the 
higher end for the total time taken to compute the Newton step.  Occasionally the 
maximum number of {GMRES} iterations was reached, in which
case the computation of the Newton step was reported as failed to the main 
numerical continuation routine. However, in most cases convergence to the desired residual tolerance was reached with 
under 2,000 {GMRES} iterations.  Solution of the system of nonlinear algebraic equations typically took 2--6 iterations for 
Newton's method (a maximum of 10 Newton iterations was set).  Upon convergence of Newton's method, 
the residual $||\F||_2$ was on the order of $10^{-7}$ or less.

As for the {LOCA} numerical continuation library \cite{loca}, recall from Section~\ref{sec:continuation_procedure}
that we performed single-parameter continuation using the
option of arc-length continuation in order to allow for turning points in the path following process.  
Input information required by the {LOCA} library was set based on behavior observed for some initial runs
as well as on recommendations provided in the documentation \cite{loca}. 
In particular, we experimented with providing the {LOCA} library values in the range $[-0.1, 1.0]$ for the initial change 
in the continuation parameter and $[0.5, 2.0]$ for the maximum increment in the continuation parameter.  We found that 
it was best to set the initial change in the continuation parameter to be in the range $[0.01, 0.05]$, and to allow a maximum 
increment in the range $[0.5, 1.0]$.  Although larger values could also perform satisfactorily, in general we found that it was 
best for our problem to keep somewhat tight control on these increments in the sense that the number of failed attempts 
was then minimal (often zero).  In addition, allowing large increments led several of the solutions to a single-frequency solution 
in the range $R \in [9,10]$ of values of the continuation parameter.  With tighter bounds on the allowed increments, the numerical
continuation led to a larger variety of solutions, as discussed at the beginning of Section~\ref{sec:results} and in
Section~\ref{sec:discussion_particular_solutions}.

Upon reaching a solution of $\Feqzero$ at the final {CGLE} parameter values, the values of 
$N_x$ and $N_t$ in the truncated expansions \refEqn{eqn:xFseries} and \refEqn{eqn:am_ansatz} were increased to 
confirm that, with the increased number of terms in the expansions, Newton's method would converge to the same solution.
(That is, to confirm that the solution of $\Feqzero$ was numerically well defined.)  
In addition, the solution was validated against time integration of the truncated system of {ODEs} \refEqn{eqn:odes}.
Finally, we note that the computations were performed on a Thinkpad W530 personal workstation with 16 {GB} memory, 
four cores, with two threads per core, and 2.70 {GHz} processor speed.  Per the discussion in
Section~\ref{sec:continuation_procedure} and the corresponding Appendix~\ref{app:newton_step}, four threads were used 
concurrently when solving for the Newton step.

\section{Concluding Remarks}
\label{sec:conclusion}

Among aspects for further consideration we mention research on techniques that may help in
minimizing or circumventing excessive traversal of parameter values in a cyclic manner,  per the discussion in
Section~\ref{sec:discussion_particular_solutions}.
This could include alternative techniques for control of the step size in 
the continuation parameter or the use of multi-parameter continuation.
A comparison with alternatives to the use of the Moore-Penrose inverse for computing the Newton step,
specifically the use of phase or gauge conditions \cite{bk:kuznetsov},
should be also performed.
Such additional features will provide a more versatile setting in which to explore further larger parameter regions
(with an increasing number of unknowns and/or higher space dimension), and enhance the
understanding of the structure of the solution space of the {CGLE}, and in particular, the structure of the space of orbits
of its symmetry group.
In addition, the fact that the resulting solutions are unstable suggests that the solutions may belong to the set of
(infinitely many) unstable periodic orbits embedded in chaotic attractors \cite{bk:chaosbook,lan10,chandler13}.
This direction, by itself, is certainly very interesting to pursue in further studies of the dynamics of the {CGLE},
and on the potential use of such periodic orbits in the study of chaotic dynamical systems
\cite{bk:chaosbook,lan10,chandler13}.

\section*{Acknowledgements}

The author thanks Ognyan Stoyanov  for useful discussions on the topic of symmetry groups of differential equations and 
helpful feedback on a preliminary version of this paper, as well as for much help with installation of the Fedora operating
system prior to setting up and performing the computations described here.
The author also thanks the editor and anonymous reviewers for their time and feedback.

\appendix \section{Appendix: Newton Step Computation}
\label{app:newton_step}
\newcommand{\RefSecJacobian}{Section~\ref{sec:discretization_jacobian}}

We work with the underdetermined system of nonlinear algebraic equations \refEqn{eqn:cgle_nleqns}
and consider a Newton step defined from the Moore-Penrose inverse \cite{bk:benisrael}.  This yields a minimum norm 
solution of the system of linear equations with the underdetermined \jacobian\ of \refEqn{eqn:cgle_nleqns} as coefficient
matrix, and is one technique used in numerical continuation methods \cite{bk:allgower,wulff06}.  
However, instead of computing the desired Newton step directly from the linear system 
having the underdetermined \jacobian\ as coefficient matrix, we premultiply the linear system with a matrix 
composed of a subset of the columns of the \jacobian\ so that a numerical 
solution for the problem at hand may be obtained in a more efficient manner.
The approach is conceptually simple, yet that is where its value lies:  it allowed us to compute an accurate Newton step 
quickly and efficiently  and made the solution of a computationally challenging problem with a large number of unknowns 
practical without the need of a cluster or supercomputer. 
The details are explained next.

Let $\M$ be a $\m \times \n$ matrix, $\m < \n$, and assume $\M$ has rank $\m$. Let $b$ be a vector of size $\n \times 1$.
Recall that the minimum norm solution $\x$ of the system of linear equations
\begin{equation}    \label{eqn:underdet_sys1}
    \M  \x   =   b
\end{equation}
given by the Moore-Penrose inverse is \cite{bk:allgower}
\begin{equation}    \label{eqn:underdet_sys_solx}
    \x   =   \M^{\mathrm{T}} (\M \M^{\mathrm{T}})^{-1} b.
\end{equation}
Computing the solution $\x$ from \refEqn{eqn:underdet_sys_solx} thus requires the solution of a linear system having 
$\M \M^{\mathrm{T}}$ as coefficient matrix.\footnote{The components of the
\jacobian\ matrix in our computations are real, hence our use of the terms \emph{transpose} and \emph{symmetric} 
when referring to the matrix $\M$ in the discussion that follows, instead of the more general terms 
\emph{conjugate transpose} and \emph{Hermitian}.}
As is well known, the numerical solution of a system of linear algebraic equations may be obtained using a variety 
of methods, either direct \cite{bk:golub} or iterative ones \cite{bk:greenbaum}.  A main distinction between these 
two classes of methods is that the use of direct methods  requires explicitly the availability of the coefficient matrix, 
whereas for iterative methods what is needed is the ability to perform matrix-vector products with the coefficient 
matrix.  Hence, iterative methods are an attractive option when explicit computation 
of the coefficient matrix is not feasible or is inconvenient, and multiplication of the coefficient matrix and a vector can 
be performed (efficiently) without explicit computation of the coefficient matrix.  As discussed in \RefSecJacobian,
the latter applies to the problem at hand.  Therefore we considered the use 
of iterative methods, in particular the generalized minimal residual (GMRES) method \cite{bk:greenbaum,saad86} 
due to its robustness and suitability for non-symmetric systems.  Since the matrix $\M \M^{\mathrm{T}}$ in
\refEqn{eqn:underdet_sys_solx} is symmetric we also explored the possibility of using the conjugate gradient (CG) 
method for symmetric systems \cite{bk:greenbaum}, but, as will be discussed below, the {GMRES} method 
is a more suitable option for our problem.

The convergence behavior of iterative methods for solving linear systems is dependent on the method as well as on
various other factors, for example, certain properties of the coefficient matrix or the problem from which the linear 
system is derived.  In the case of the {GMRES} method, one desirable property for fast convergence is for the
eigenvalues of the coefficient matrix to be clustered around a few values, away from zero \cite{bk:greenbaum}.  
Another important issue is that the use of iterative methods for solving linear systems typically requires the use of a 
preconditioner in order to perform efficiently.  Generally speaking, preconditioning refers to multiplying the linear system 
(on the left or right) by a matrix such that the resulting system has the properties needed for optimal or enhanced
performance of the particular method under consideration (and yields a solution for the unpreconditioned (i.e., original)
linear system).  
For thorough treatments on iterative methods for solving linear systems, the interested reader is referred 
to \cite{bk:greenbaum} and references therein.
Here  we restrict ourselves to a brief discussion on the behavior resulting from the
use of the {GMRES} and {CG} methods for the problem at hand.

\begin{figure}[!t]  
    \begin{center}
         {\includegraphics[trim = 10mm 0mm 18mm 5mm, clip, width=0.60\textwidth]{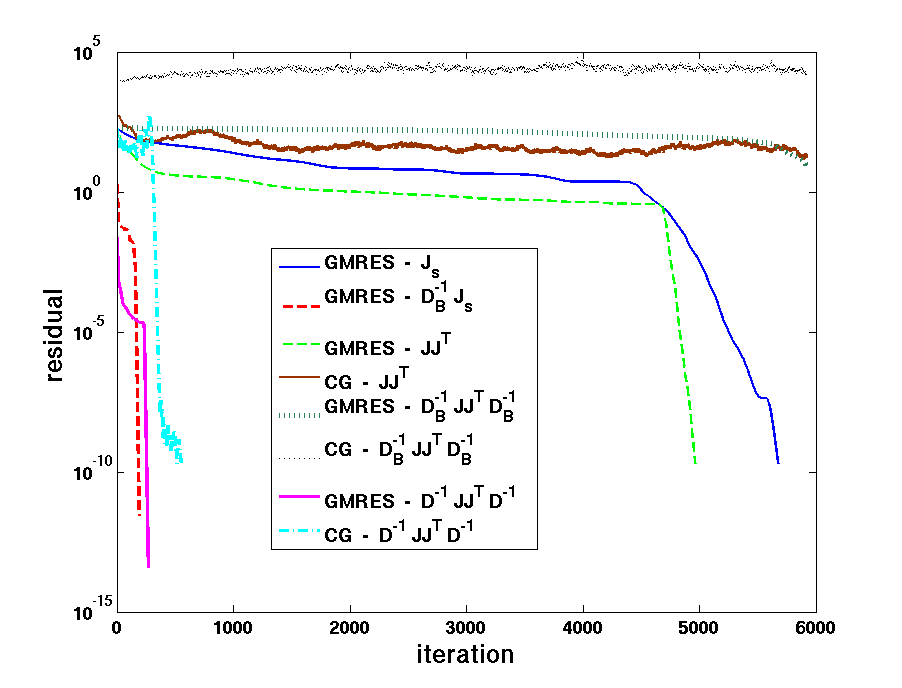}}
    \end{center}
    \caption{{Typical convergence behavior of {GMRES} and {CG} for the problem at hand.  
                   Suitable convergence behavior results for {GMRES} with $D^{-1}_{\sss B} \M_s$ as 
                   coefficient matrix, and for {GMRES} or {CG} with $D^{-1} \M \M^\mathrm{T} D^{-1}$ as coefficient 
                   matrix.  However, the only viable option for our problem is the {GMRES} method with $D^{-1}_{\sss B} \M_s$ as 
                   coefficient matrix. }}
    \label{fig:linsys_convergence}
\end{figure}
Figure~\ref{fig:linsys_convergence} depicts the typical convergence behavior exhibited by the {GMRES} and {CG} 
methods when used to compute the solution of systems of linear equations having as coefficient matrix the
\jacobian\ of the system \refEqn{eqn:cgle_nleqns} of nonlinear algebraic equations, after splitting the equations 
into their real and imaginary parts.
In the figure,  $\M$ denotes the underdetermined \jacobian\ matrix of the system \refEqn{eqn:cgle_nleqns},
$\M_s$ represents a square, non-singular matrix composed of a subset of columns of $\M$,
the preconditioner $D_{\sss B}$ is the \jacobian\ matrix of the terms in \refEqn{eqn:cgle_nleqns} linear in the unknowns 
$\{\amn{m}{n}\}$, with columns corresponding to the derivatives with respect to the real and imaginary parts of the coefficients
$\{\amn{m}{n}\}$ (it is a block diagonal matrix),
and the preconditioner $D$ is a diagonal matrix having the diagonal elements of $\M \M^\mathrm{T}$ on its diagonal.
As can be seen from Figure~\ref{fig:linsys_convergence}, suitable convergence behavior results only from the use 
of {GMRES} with $D^{-1}_{\sss B} \M_s$ as coefficient matrix (i.e., the use of {GMRES} to solve linear systems
with $\M_s$ as coefficient matrix and $D_{\sss B}$ as a preconditioner), 
as well as from using {GMRES} or {CG} with $D^{-1} \M \M^\mathrm{T} D^{-1}$ as coefficient matrix
(i.e., the use of {GMRES} or {CG} to solve linear systems
with $\M \M^\mathrm{T} $ as coefficient matrix and $D$ as a preconditioner).

The block diagonal preconditioner $D_{\sss B}$ is very effective when used with the {GMRES} method to solve linear 
systems with $\M_s$ as coefficient matrix, as seen in Figure~\ref{fig:linsys_convergence}.  For the example
depicted, it took 190 iterations to solve for a system having 5,925 unknowns.  The corresponding run with unpreconditioned 
{GMRES} required 5,676 iterations, making it impractical for our purposes.  Using the preconditioner $D$ and both the
{GMRES} and {CG} methods to solve for systems having $\M \M^\mathrm{T}$ as coefficient matrix also gave good 
results.  For the example in Figure~\ref{fig:linsys_convergence}, these two methods took, respectively, 269 and 552 
iterations to solve the preconditioned system.  However, per the discussion in \RefSecJacobian,
the block diagonal preconditioner $D_{\sss B}$ is readily available and easy to manipulate, whereas
assembling the preconditioner $D$ requires calculating the diagonal terms of the matrix $\M \M^\mathrm{T}$, 
and these terms are not readily available for our problem.  Therefore, the only viable option for us is the use of the 
{GMRES} method with preconditioner $D_{\sss B}$ to solve systems having $\M_s$ as the coefficient matrix.
As a result, the computation of the minimum norm solution $\x$ directly from 
\refEqns{eqn:underdet_sys1}{eqn:underdet_sys_solx} is unfeasible for the problem at hand.

Returning to the linear system in \refEqn{eqn:underdet_sys1}, we thus express $\M$ as composed by two matrices,
\begin{equation}    \label{eqn:A_decomposed}
    \M   =   \left [ \,\M_s\,  |  \,\M_r\, \right ],
\end{equation}
where $\M_s$ has dimension $\m \times \m$ (and is non-singular) and $\M_r$ has dimension $\m \times (\n \!-\! \m)$.
Now we consider the system $\M_s^{-1} \M \, \x   =   \M_s^{-1} b$ obtained by multiplying \refEqn{eqn:underdet_sys1}
with $\M_s^{-1}$, yielding
\begin{equation}    \label{eqn:underdet_sys3}
    \left [ \,I\,  |  \,\M_s^{-1} \M_r\, \right ]  \x   =   \M_s^{-1} b,
\end{equation}
where $I$ is the $\m \times \m$ identity matrix.  We work directly with the system \refEqn{eqn:underdet_sys3} and
compute the desired solution $\x$ by solving two sub-problems, namely:
\begin{enumerate}
    \item Compute the right-hand side $\M_s^{-1} b$, as well as the $\n \!-\! \m$ columns of the submatrix $\M_s^{-1} \M_r$ 
              of the coefficient matrix $\left [ \,I\,  |  \,\M_s^{-1} \M_r\, \right ]$ in \refEqn{eqn:underdet_sys3}.
              This sub-problem will therefore require the solution of $\n \!-\! \m \!+\! 1$ linear systems having $\M_s$ as 
              coefficient matrix.  It will be feasible if solving linear systems having $\M_s$ as coefficient matrix can be
              done efficiently and if $\n \!-\! \m \!+\! 1$ is small.  Both of these requirements are satisfied in our
              study since first, per the discussion from the previous paragraphs, we can use the {GMRES}
              method to solve the required linear systems efficiently, and second, for our problem, $\n \!-\! \m = 3$ due to the 
              3-tuple $( \varphi,S,T )$ of additional unknowns in the problem formulation.
    \item Upon completion of sub-problem 1, compute the minimum norm solution given by the Moore-Penrose inverse 
              for the underdetermined system of linear equations \refEqn{eqn:underdet_sys3}.
\end{enumerate}
Per sub-problem 1 above, one first needs to solve $k = 1 , \ldots, \n \!-\! \m \!+\! 1$ linear systems
\begin{equation}    \label{eqn:sub_linear_sys}
    \M_{{\sss s}} \, \tilde{\x}_{{\sss k}} = \tilde{b}_{{\sss k}}, 
\end{equation}
where, for  $k = 1, \ldots, \n \!-\! \m$, the $k$-th linear system \refEqn{eqn:sub_linear_sys} will have the right-hand 
side vector $\tilde{b}_{{\sss k}}$ equal to the $k$-th column of the matrix $\M_r$, so that the solutions $\tilde{\x}_{{\sss k}}$ 
of said $\n \!-\! \m$ linear systems yield the columns of the submatrix 
$\M_s^{-1} \M_r$ of the coefficient matrix $\left [ \,I\,  |  \,\M_s^{-1} \M_r\, \right ]$ in \refEqn{eqn:underdet_sys3}.
The solution of the remaining linear system \refEqn{eqn:sub_linear_sys}, with $\tilde{b}_{{\sss \n \!-\! \m \!+\! 1}} = b$, 
yields the right-hand side $\M_s^{-1} b$ in \refEqn{eqn:underdet_sys3}.  Solving these $\n \!-\! \m \!+\! 1$ 
linear systems \refEqn{eqn:sub_linear_sys} is not an obstacle since they can be solved either independently, 
in parallel, or with an implementation of a (direct or iterative) solver for linear systems that handles 
multiple right-hand sides.  Recall that the viable option for us is to use the {GMRES} method for linear systems.  
In our implementation, we combined it with the use of {POSIX} threads (pthreads) programming \cite{pthreads} in order 
to solve the required linear systems \refEqn{eqn:sub_linear_sys} in parallel.  Hence, in the current study, the solution 
of the $\n \!-\! \m \!+\! 1 = 4$  linear systems \refEqn{eqn:sub_linear_sys} having $\M_s$ as coefficient matrix was achieved
basically in the same amount of time as that required to solve a single such system.

Upon completion of  sub-problem 1, what remains to be done is to compute the minimum norm solution $\x$ from
the system in \refEqn{eqn:underdet_sys3}.  As previously noted, one desirable property for fast convergence of the 
{GMRES} method is for the eigenvalues of the coefficient matrix to be clustered around a few values, away from zero, 
since, typically, the number of iterations required for convergence when using the {GMRES} method  depends on the 
number of distinct eigenvalues of the coefficient matrix of the linear system \cite{bk:greenbaum}.
Thus, we also used the {GMRES} method to solve for the minimum norm solution $\x$ in 
\refEqn{eqn:underdet_sys3}, since this requires solving a linear system having
\begin{eqnarray}    \label{eqn:Bmatrix}
    I + (\U) (\U)^{\mathrm{T}}
\end{eqnarray}
as coefficient matrix, and the matrix \refEqn{eqn:Bmatrix}
has all but $\n \!-\! \m$ eigenvalues equal to one, with the remaining eigenvalues greater than or equal to one.
(Recall that  $\M_s$ is a non-singular $\m \times \m$ matrix and $\M_r$ has dimension $\m \times (\n\!-\!\m)$, where $\m < \n$.)
This follows from the fact that the eigenvalues $\lambda$ and eigenvectors $y$ of the 
matrix \refEqn{eqn:Bmatrix} satisfy
\begin{equation}    \label{eqn:eigs2}
    (\U) (\U)^{\mathrm{T}} \, y   \ = \    (\lambda - 1) \, y.    
\end{equation}
Noting that the null space of the matrix $(\U)^{\mathrm{T}}$ has dimension (at least) $\m - (\n \!-\! \m) = 2\m \!-\! \n$, 
it follows that the matrix $(\U) (\U)^{\mathrm{T}}$ has zero as an eigenvalue, that is, $\lambda = 1$,
with multiplicity (at least) $2\m \!-\! \n$.  Furthermore, $(\U) (\U)^{\mathrm{T}}$ is 
positive semi-definite and symmetric, so its eigenvalues $(\lambda \!-\! 1)$ in \refEqn{eqn:eigs2} are non-negative 
and, thus, the remaining (at most) $\n \!-\! \m$ eigenvalues of the matrix \refEqn{eqn:Bmatrix} satisfy $\lambda \ge 1$.  
The significance here is that solving a linear system with the matrix \refEqn{eqn:Bmatrix}  as coefficient matrix, 
which is required in order to compute the minimum norm solution $\x$ of system \refEqn{eqn:underdet_sys3}, 
should take $\n \!-\! \m \!+\! 1$ iterations if we use the {GMRES} method.  For our problem,  this means 
$\n \!-\! \m \!+\! 1 = 4$ iterations.  Note also that computing matrix-vector products with the matrix \refEqn{eqn:Bmatrix} 
can be easily done (since the $\n \!-\! \m = 3$ columns of the matrix $\M_s^{-1} \M_r$ have been previously computed 
and stored in memory).  Furthermore, no preconditioning is required to solve linear systems having \refEqn{eqn:Bmatrix} 
as coefficient matrix.  Hence, using the {GMRES} method to solve the aforementioned sub-problem 2 poses no difficulty
and results in a negligible amount of additional computing time when solving for the Newton step using the proposed
approach. 

Now, denoting the system \refEqn{eqn:cgle_nleqns} as $\Feqzero$, note that the vector $b$ in \refEqn{eqn:underdet_sys3}
corresponds to $-\F$ evaluated at the current solution estimate (from Newton's method).
Also, based on our presentation of the material, it may seem natural to consider the matrix $\Ja$ introduced in
\RefSecJacobian, whose columns correspond to derivatives of $\F$ with 
respect to the real and imaginary parts of the unknowns $\{\amn{m}{n}\}$, as that corresponding to the matrix 
$\M_s$ in \refEqns{eqn:A_decomposed}{eqn:underdet_sys3}.  
Recall, though, that $\Ja$ is singular at a solution of $\Feqzero$.
We therefore define the matrix $\M_s$ as that obtained by replacing three columns from $\Ja$ by the columns of 
the \jacobian\ matrix of $\F$ corresponding to derivatives with respect to the unknowns $\varphi$, $S$, and $T$.  
This proved effective in dealing with said singularity when computing the Newton step 
from the corresponding system \refEqn{eqn:underdet_sys3} during the numerical continuation.
(As discussed in \cite{lopez05}, the kernel of the \jacobian\ matrix at a solution of $\Feqzero$ is typically three-dimensional.)
The replaced columns from $\Ja$ define the columns of the matrix $\M_r$ in \refEqn{eqn:A_decomposed}, 
which are needed to construct the matrix $\M_s^{-1} \M_r$ in \refEqn{eqn:underdet_sys3}.  These three columns 
can be computed (efficiently) via matrix-vector products of the \jacobian\ matrix 
(see \RefSecJacobian) and standard basis vectors.

To summarize, we used the implementation of the {GMRES} method \cite{saad86} available from the 
Meschach software package \cite{meschach} to solve the required linear systems in
\refEqn{eqn:underdet_sys3}, as well as that having coefficient matrix \refEqn{eqn:Bmatrix}.  
For specific details about the {GMRES} method itself,
the reader is referred to \cite{saad86,meschach}.  Here we note that to solve, say, the linear system
$Ax=b$ using the {GMRES} solver \cite{meschach}, the user must provide a routine that computes the 
matrix-vector product $Ay$, given a vector $y$.  
Recall from the discussion in Section~\ref{sec:discretization_jacobian} that for our problem the required
matrix-vector products can be computed efficiently using fast Fourier transforms (FFTs).  In particular,
our implementation used the {FFTW} software package \cite{fftw} for this purpose.
If a preconditioner $D$ is to be used with the {GMRES} solver, the user must also provide a routine that
computes the solution $z$ of the linear system $Dz=d$, given a right-hand side vector $d$.
As noted in this appendix (see the discussion relating to Figure~\ref{fig:linsys_convergence}), 
this requirement does not pose a difficulty for us.  It required solving a linear system with a block diagonal
matrix, and such solution was easily implemented directly (that is, no iterative method was required).
In the calls to the {GMRES} solver \cite{meschach}, a tolerance of $10^{-9}$ was set for the residual, along with
a maximum of 3,000 {GMRES} iterations.  
An outline of the computation of the Newton step is as follows:
\begin{myenumerate}
    \item Compute the vector $b$ and the three columns of the matrix $\M_r$, as defined in the preceding paragraph.
             These compose the right-hand sides for the four linear systems 
             required to be solved in \refEqn{eqn:underdet_sys3}, which have $\M_s$ as coefficient matrix.
    \item Solve the linear systems in  \refEqn{eqn:underdet_sys3}, with the four right-hand sides computed
             in step 1 above.  The implementation was done using {POSIX} threads (pthreads) programming \cite{pthreads}, 
             so that the four linear systems were solved
             independently, in parallel.  In other words, four threads were used, each thread making a call to the 
             {GMRES} solver \cite{meschach}, with its corresponding right-hand side.  Since the four linear systems were
             solved concurrently, the computational time was essentially that required to solve a single such
             linear system.
    \item Compute the minimum norm solution given by the Moore-Penrose inverse 
             for the underdetermined system of linear equations \refEqn{eqn:underdet_sys3}.
             This was implemented with a call to the {GMRES} method \cite{meschach} to solve a linear system having
             coefficient matrix as in \refEqn{eqn:Bmatrix}.  Per the discussion following \refEqn{eqn:Bmatrix},
             the required matrix-vector products with the matrix \refEqn{eqn:Bmatrix} are easily implemented
             using the results from step 2 above.  Recall also from the discussion that computing this minimum norm solution
             required four iterations of the {GMRES} method and therefore the amount of computational time required
             in addition to that from step 2 above was negligible. 
\end{myenumerate}



\begin{thebibliography}{10}

\bibitem{bk:allgower}
Eugene~L. Allgower and Kurt Georg.
\newblock {\em Introduction to Numerical Continuation Methods}.
\newblock Classics in Applied Mathematics. SIAM, 2003.

\bibitem{aranson02}
Igor~S. Aranson and Lorenz Kramer.
\newblock The world of the complex {Ginzburg-Landau} equation.
\newblock {\em Rev. Mod. Phys.}, 74:99--143, 2002.

\bibitem{aston99}
Philip~J. Aston and Carlo~R. Laing.
\newblock Symmetry and chaos in the complex {Ginzburg-Landau} equation. {I}.
  {Reflectional} symmetries.
\newblock {\em Dynam. Stability Systems}, 14:233--253, 1999.

\bibitem{aston00}
Philip~J. Aston and Carlo~R. Laing.
\newblock Symmetry and chaos in the complex {Ginzburg-Landau} equation. {II}.
  {Translational} symmetries.
\newblock {\em Physica~D}, 135:79--97, 2000.

\bibitem{bk:benisrael}
Adi Ben-Israel and Thomas~N.E. Greville.
\newblock {\em Generalized Inverses : Theory and Applications}.
\newblock Springer-Verlag, 2003.

\bibitem{benettin80a}
G.~Benettin, L.~Galgani, A.~Giorgilli, and J.-M. Strelcyn.
\newblock {Lyapunov} characteristic exponents for smooth dynamical systems and
  for {Hamiltonian} systems; {A} method for computing all of them. {Part I:
  Theory}.
\newblock {\em Meccanica}, 15:9--20, 1980.

\bibitem{brusch01}
Lutz Brusch, Alessandro Torcini, Martin van Hecke, Mart{\'{\i}}n~G. Zimmermann,
  and Markus B{\"{a}}r.
\newblock Modulated amplitude waves and defect formation in the one-dimensional
  complex {Ginzburg-Landau} equation.
\newblock {\em Physica~D}, 160:127--148, 2001.

\bibitem{chandler13}
Gary~J. Chandler and Rich~R. Kerswell.
\newblock Invariant recurrent solutions embedded in a turbulent two-dimensional
  {Kolmogorov} flow.
\newblock {\em J. Fluid Mech.}, 722:554--595, 2013.

\bibitem{bk:chaosbook}
P.~Cvitanovi{\'{c}}, R.~Artuso, P.~Dahlqvist, R.~Mainieri, G.~Tanner,
  G.~Vattay, N.~Whelan, and A.~Wirzba.
\newblock {\em Chaos: Classical and Quantum}.
\newblock 2014.
\newblock Webbook available at \texttt{chaosbook.org}.

\bibitem{doelman89}
Arjen Doelman.
\newblock Slow time-periodic solutions of the {Ginzburg-Landau} equation.
\newblock {\em Physica~D}, 40:156--172, 1989.

\bibitem{Doelman}
Arjen Doelman and Edriss~S. Titi.
\newblock Regularity of solutions and the convergence of the {G}alerkin method
  in the {G}inzburg-{L}andau equation.
\newblock {\em Numer. Funct. Anal. Optim.}, 14(3-4):299--321, 1993.

\bibitem{fftw}
M.~Frigo and S.~G. Johnson.
\newblock {FFTW}: An adaptive software architecture for the {FFT}.
\newblock In {\em {ICASSP} Conference Proceedings}, volume~3, pages 1381--1384,
  1998.
\newblock \texttt{http://www.fftw.org/}.

\bibitem{bk:golub}
Gene~H. Golub and Charles F.~Van Loan.
\newblock {\em Matrix Computations}.
\newblock The Johns Hopkins University Press, 1996.

\bibitem{bk:greenbaum}
Anne Greenbaum.
\newblock {\em Iterative Methods for Solving Linear Systems}.
\newblock Frontiers in Applied Mathematics. SIAM, 1997.

\bibitem{holmes86}
Philip Holmes.
\newblock Spatial structure of time-periodic solutions of the {Ginzburg-Landau}
  equation.
\newblock {\em Physica~D}, 23:84--90, 1986.

\bibitem{Jolly}
M.~S. Jolly, R.~Temam, and C.~Xiong.
\newblock Convergence of a chaotic attractor with increased spatial resolution
  of the {G}inzburg-{L}andau equation.
\newblock {\em Chaos Solitons Fractals}, 5(10):1833--1845, 1995.

\bibitem{kapitula96}
Todd Kapitula and Stanislaus Maier-Paape.
\newblock Spatial dynamics of time-periodic solutions for the {Ginzburg-Landau}
  equation.
\newblock {\em Zeitschrift f{\"{u}}r Angewandte Mathematik und Physik},
  47:265--305, 1996.

\bibitem{keefe85}
Laurence~R. Keefe.
\newblock Dynamics of perturbed wavetrain solutions to the {G}inzburg-{L}andau
  equation.
\newblock {\em Stud. Appl. Math.}, 73(2):91--153, 1985.

\bibitem{bk:kuznetsov}
Y.~A. Kuznetsov.
\newblock {\em Elements of Applied Bifurcation Theory}.
\newblock Springer-Verlag, 1998.

\bibitem{lan10}
Y.~Lan.
\newblock Cycle expansions: From maps to turbulence.
\newblock {\em Commun Nonlinear Sci Numer Simulat}, 15:502--526, 2010.

\bibitem{levermore96}
C.~David Levermore and Marcel Oliver.
\newblock The complex {Ginzburg-Landau} equation as a model problem.
\newblock In {\em Dynamical Systems and Probabilistic Methods in Partial
  Differential Equations}, volume~31 of {\em Lectures in Appl. Math.}, pages
  141--190. Amer. Math. Soc., Providence, RI, 1996.

\bibitem{lloyd05}
D.~J.~B. Lloyd, A.~R. Champneys, and R.~E. Wilson.
\newblock Robust heteroclinic cycles in the one-dimensional complex
  {Ginzburg-Landau} equation.
\newblock {\em Physica~D}, 204:240--268, 2005.

\bibitem{lopez05}
Vanessa L\'{o}pez, Philip Boyland, Michael~T. Heath, and Robert~D. Moser.
\newblock Relative periodic solutions of the complex {Ginzburg-Landau}
  equation.
\newblock {\em {SIAM} J. Appl. Dyn. Syst.}, 4(4):1042--1075, 2005.

\bibitem{pthreads}
Linux~Programmer's Manual.
\newblock {POSIX} threads programming.
\newblock \texttt{http://man7.org/linux/man-pages/man7/pthreads.7.html}.

\bibitem{mielke02}
A.~Mielke.
\newblock The {Ginzburg-Landau} equation in its role as a modulation equation.
\newblock In {\em Handbook of Dynamical Systems, Vol. 2}, B. Fiedler, ed.,
  pages 759--834. Elsevier Science, 2002.

\bibitem{Moon}
H.~T. Moon, P.~Huerre, and L.~G. Redekopp.
\newblock Transitions to chaos in the {G}inzburg-{L}andau equation.
\newblock {\em Phys. D}, 7(1-3):135--150, 1983.

\bibitem{bk:ortega}
J.~M. Ortega and W.~C. Rheinboldt.
\newblock {\em Iterative Solution of Nonlinear Equations in Several Variables}.
\newblock Classics in Applied Mathematics. SIAM, 2000.

\bibitem{bk:ott}
Edward Ott.
\newblock {\em Chaos in Dynamical Systems}.
\newblock Cambridge University Press, 2002.

\bibitem{saad86}
Y.~Saad and M.~H. Schultz.
\newblock {GMRES}: A generalized minimal residual algorithm for solving
  nonsymmetric linear systems.
\newblock {\em {SIAM} J. Sci. Stat. Comput.}, 7:856--869, 1986.

\bibitem{loca}
A.~G. Salinger, M.~Bou-Rabee, R.~P. Pawlowski, E.~D. Wilkes, E.~A. Burroughs,
  R.~B. Lehoucq, and L.~A. Romero.
\newblock {LOCA}: A library of continuation algorithms - theory and
  implementation manual.
\newblock Technical report, Sandia National Laboratory, 2001.
\newblock \texttt{http://www.cs.sandia.gov/projects/loca/index.html}.

\bibitem{meschach}
David~E. Stewart and Zbigniew Leyk.
\newblock Meschach library.
\newblock Technical report, Australian National University, 1994.
\newblock \texttt{http://www.netlib.org/c/meschach/readme}.

\bibitem{takac98}
Peter Tak{\'a}{\v{c}}.
\newblock Bifurcations to invariant 2-tori for the complex {G}inzburg-{L}andau
  equation.
\newblock {\em Appl. Math. Comput.}, 89:241--257, 1998.

\bibitem{vanHecke02}
Martin van Hecke.
\newblock Coherent and incoherent structures in systems described by the {1D}
  {CGLE}: Experiments and identification.
\newblock {\em Physica~D}, 174:134--151, 2003.

\bibitem{vanSaarloos94}
W.~van Saarloos.
\newblock The complex {Ginzburg-Landau} equation for beginners.
\newblock In {\em Spatio-Temporal Patterns in Nonequilibrium Complex Systems},
  P. E. Cladis and P. Palffy-Muhoray, eds., Studies in the Sciences of
  Complexity, Proceedings XXI. Addison-Wesley, Reading, MA, 1994.

\bibitem{wulff06}
Claudia Wulff and Andreas Schebesch.
\newblock Numerical continuation of symmetric periodic orbits.
\newblock {\em {SIAM} J. Appl. Dyn. Syst.}, 5(3):435--475, 2006.

\end{thebibliography}
\end{document}